\documentclass[reqno]{amsart}
\usepackage[utf8x]{inputenc}
\usepackage{amsfonts,amsmath,amssymb,mathtools,enumitem,bm}
\usepackage{datetime}
\usepackage{hyperref}
\hypersetup{colorlinks=true, pdfstartview=FitV, linkcolor=BrickRed,citecolor=BrickRed, urlcolor=BrickRed, linktoc=page}

\usepackage{graphicx}
\usepackage[dvipsnames]{xcolor}
\usepackage{caption}
\usepackage{subcaption}

\definecolor{labelkey}{rgb}{0.6,0,0}
\usepackage{cancel}

\mathtoolsset{showonlyrefs=true}

\usepackage[margin=2.5cm]{geometry}
\numberwithin{equation}{section}

\def\dt{\partial_t}
\def\kap{\kappa}
\def\jpt{\langle t \rangle}
\def\U{\mathcal{U}_\pm}
\usepackage{dsfont}

\newtheorem{theorem}{Theorem}[section]
\newtheorem{lemma}[theorem]{Lemma}
\newtheorem{corollary}[theorem]{Corollary}
\newtheorem{proposition}[theorem]{Proposition}

\newtheorem{remark}[theorem]{Remark}

\title[Global solutions to Navier-Stokes-Coriolis equations]{\vspace*{-1cm}Global axisymmetric solutions for Navier-Stokes equation with rotation uniformly in the inviscid limit}

\author{Haram Ko}
\thanks{Department of Mathematics, Brown University, Providence, RI, USA, haram\_ko@brown.edu}

\begin{document}

\subjclass[2020]{35Q30, 76D03, 76D05, 76U05}
\keywords{Navier-Stokes equations, rotating fluids, nonlinear stability, dispersion}

\begin{abstract}
We prove that the solutions to the 3D Navier-Stokes equation with constant rotation exist globally for small axisymmetric initial data, where the smallness is uniform with respect to the viscosity $\nu \in [0,\infty)$. This expands the work by Guo, Pausader, and Widmayer \cite{GPW} which showed the global axisymmetric stability of rotation for 3D incompressible Euler's equation, to the viscous case, but for a single threshold that works for arbitrary viscosity.
This is achieved by suitably adapting the dispersive framework established in \cite{GPW} to the Navier-Stokes setting.
\end{abstract}

\setcounter{tocdepth}{1}
\maketitle
\tableofcontents
\vspace*{-.75cm}

\section{Introduction}

We consider the following 3D incompressible Navier-Stokes-Coriolis equation, which describes the evolution of an incompressible viscous fluid governed by the Navier-Stokes equation in a rotating environment:
\begin{equation}\label{eq:fullNSC}
	\begin{cases}
		\dt U - \nu \Delta U + (U\cdot\nabla)U + \delta^{-1} e_3 \times U + \nabla P & = 0, \\
		\text{div}(U) & = 0,
	\end{cases} \qquad (t,x) \in \mathbb{R}_+\times\mathbb{R}^3,
\end{equation}
where $\nu$ is the viscosity of the fluid and $\delta^{-1}$ is the speed of the rotating background. From hereon, whenever there is no source of confusion, we will simply write the Navier-Stokes(-Coriolis) equation to denote the incompressible Navier-Stokes(-Coriolis) equation.

The rotation of the background is known to be an important source of macroscopic properties of fluids, and is able to create global behaviors, e.g. the Gulf Stream. We refer to \cite{GFD, CDGG06} for more details on the geophysical importance of \eqref{eq:fullNSC}. Furthermore, there is a $1$-$1$ correspondence between solutions of \eqref{eq:NSC} and the solutions of Navier-Stokes equation, since the well-known steady-state solutions, constant speed rotations $U_{\delta}(x,y,z) = \delta^{-1}(-y,x,0)$, give rise to \eqref{eq:fullNSC} when we study the evolution of its axisymmetric perturbation. In this regard, our work can be thought of as constructing global solutions of the Navier-Stokes equation with prescribed behavior at infinity. \\

It is well known that solutions of the Navier-Stokes equations exist globally for small initial data in appropriate spaces \cite{FK, Ka, KT}. Indeed, the same is true for the Navier-Stokes-Coriolis equation, but even better results are true if we introduce rotation, by exploiting the extra degree of freedom given by the speed of rotation $\delta^{-1}$. Chemin, Desjardins, Gallagher, and Grenier \cite{CDGG02} proved the existence of global solutions of \eqref{eq:fullNSC} for any arbitrarily large, low regularity data of a certain form and small enough $\delta$ which depends on the initial data. We can interpret this as a first answer to the following question, using the concept of stability threshold similar to one used in \cite{BGM}:\\\\
\textbf{Question} : \textit{Given a norm $X$, let the threshold $\epsilon = \epsilon_X(\nu,\delta)$ be the largest value such that any initial data $u_0$ of $||u_0||_X < \epsilon$ leads to global solution of \eqref{eq:fullNSC}. How does $\epsilon_X(\nu,\delta)$ behave as $\nu,\delta$ vary?} \\\\
Several results can be interpreted using this framework. Koh, Lee, and Takada \cite{KLT} found an explicit function for the lower bound of $\epsilon(\nu,\delta)$ using the Strichartz estimates in Sobolev spaces. See also \cite{IT} and \cite{SYC} for the preceding or more recent results in different function spaces along a similar line. While these results mainly tried to include more initial data as $\delta \to 0$, in the other extreme, \cite{GIMS} and \cite{HS} have proved that as $\delta \to \infty$, the threshold becomes independent of $\delta$ in certain spaces. However, all of the aforementioned works have been focused on bounding $\epsilon(\nu,\delta)$ in terms of $\delta$ with fixed $\nu$, usually by setting $\nu=1$. As such, for the behavior of $\epsilon(\nu,\delta)$ as $\nu$ varies, there has been no better result than what we can get from scaling, which only gives the trivial claim that the threshold is not less than 0 when $\nu \to 0$. In summary, with less emphasis on the norms, the following have been established:
\begin{equation}\label{best-results}
	\lim_{\delta \to 0} \epsilon(\nu,\delta) = \infty, \quad
	\lim_{\nu \to \infty} \epsilon(\nu,\delta) = \infty, \quad
	\varliminf_{\delta \to \infty} \epsilon(\nu,\delta) \gtrsim \nu.
\end{equation}

At the heart of all analysis lies the fact that rotation brings dispersion to the Euler and Navier-Stokes systems. In \cite{GHPW, GPW}, Guo, Pausader, and Widmayer showed that this dispersive effect can give rise to global solutions to Euler-Coriolis equations for axisymmetric data. In our terminology, this amounts to saying
\begin{equation*}
	\epsilon(\nu=0, \delta=1) \geq \epsilon_0 > 0,
\end{equation*}
where $\epsilon_0$ is an absolute constant. This strongly indicates that we may be able to complement \eqref{best-results} with nontrivial lower bound for inviscid limit case.

The goal of this paper is to establish the existence of global solutions to \eqref{eq:fullNSC} with initial data smaller than a uniform threshold $\epsilon^*$, for continuum of viscosities $\nu \in \ [0,\infty)$. As $\nu \to 0$, we regain the Euler-Coriolis case \cite{GPW}, and $\nu=1$ represents the unit viscosity case, which has been extensively studied in the literature. To shift the role of major parameter from the speed of rotation to the viscosity, we scale variables to
\begin{equation}\label{scale-fullNSC-to-kappa}
	u(t,x) = \delta \cdot U(\delta t, x), \quad p(x) = \delta^2 \cdot P(x),
\end{equation}
to fix the rotation speed to 1 and work with the following equation with $\kap=\nu\delta$:
\begin{equation}\label{eq:NSC}
	\begin{cases}
		\dt u - \kap \Delta u + (u\cdot\nabla)u + e_3 \times u + \nabla p & = 0, \\
		\text{div}(u) & = 0.
	\end{cases}
\end{equation}

\begin{theorem}\label{Thm-abbrev}
	(a) There exists $\epsilon^* > 0, N_0 > 0$ and a norm $Y$, independent of $\kappa \in [0,\infty)$, such that for any axisymmetric function $u_0 \in H^{2N_0}(\mathbb{R}^3)$ with $||u_0||_{Y \cap H^{2N_0}} < \epsilon^*$, there exists a global solution to the Cauchy problem of \eqref{eq:NSC}, $u \in C([0,\infty); H^{2N_0})$.
	
	In particular, the solution $U(x,y,z) = (-y,x,0)$ of the Navier-Stokes equation
	\begin{equation}
		\begin{cases}
			\dt U - \nu \Delta U + (U\cdot\nabla)U + \nabla P & = 0, \\
			\text{\normalfont div}(U) & = 0,
		\end{cases} \qquad (t,x) \in \mathbb{R}_+\times\mathbb{R}^3,
	\end{equation}
	is globally stable under axisymmetric perturbation, uniformly for all $0 \leq \nu < \infty$.\\
	(b) Two solutions $U_1, U_2$ to \eqref{eq:fullNSC} of viscosities $\nu_1, \nu_2$, respectively, with the same axisymmetric initial data $U_0$ whose rescaling by \eqref{scale-fullNSC-to-kappa} satisfies the assumptions in (a) above, satisfy
	\begin{equation}\label{finite-time-convergence}
		||U_1(t) - U_2(t)||_{L^2}^2 \lesssim |\nu_1 - \nu_2| \cdot \delta^{-2-C\epsilon^*} (\epsilon^*)^2 t^{1 + C \epsilon^*}
	\end{equation}
	for small enough $\nu_1, \nu_2$ for some $C$.
\end{theorem}

The significance of Theorem \ref{Thm-abbrev}(a) lies in the fact that we have a single threshold $\epsilon^*$ and a norm that guarantee global existence for all small viscosities up to 0; see Theorem \ref{MainThmDetailed} for more detailed statement. The threshold only grows for large $\kappa$ by the standard method which allows us to have one threshold $\epsilon^*$ for all viscosities. Also, Theorem \ref{Thm-abbrev}(b) tells us that this uniformity implies continuous dependence of solutions in viscosity for finite times. In particular, the continuous dependence is true including $\nu = 0$, showing that the inviscid limit holds in finite time intervals. This provides theoretical justification and quantitative error bounds for approximating the Euler-Coriolis equation using the Navier-Stokes-Coriolis equation. Given that it is often numerically advantageous to simulate the Navier-Stokes equations than the Euler equation, we believe this could be useful to scientists and engineers who use solutions to 3D Navier-Stokes equations with small viscosity as alternate solutions to the 3D Euler equations in a rotating environment. \\

We note that this is not the first result to focus on the viscosity in Navier-Stokes-Coriolis system. In \cite{BMN97}, the authors obtained uniform (in $\nu$ not in $\delta$) bounds on the difference between the solutions of \eqref{eq:fullNSC} and their approximate solutions in generic torus, but their final results on $\epsilon(\nu,\delta)$ for torus again relied on $\nu$ for $\nu \in (0,1)$; see also \cite{BMN99, BMN01}. Chemin, Desjardins, Gallagher, and Grenier \cite{CDGG00} studied the fluid under anisotropic viscosity, where the viscosity in the direction parallel to the rotating axis, $\nu_V$, can be different from the horizontal viscosity $\nu$, and proved that the threshold size for small data global existence is indepndent in $\nu_V$. In this paper, we improved these for axisymmetric functions by obtaining independence of $\epsilon(\nu,\delta)$ in the vanishing limit of isotropic single $\nu$.

\subsection{Quantitative comparison with previous results}

The uniformity with respect to $\kap$ gives an improvement over the previous results. A concrete way to witness this fact is via scaling, and we give two examples how explicitly our result expands the known regime for global existence.

\begin{itemize}

\item In \cite{CDGG02}, Chemin, Desjardins, Gallagher, and Grenier proved the global existence of \eqref{eq:fullNSC}, and showed that the solution converges to the sum of bidimensional Navier-Stokes flow $\bar{U}$ and the free linear part of \eqref{eq:fullNSC}, $W_F^{\delta}$, when $\delta$ is small enough. However, both results implicity depend on $\nu$, and required $\nu^{-\frac{1}{2}} ||\overline{U_0}||_{L^2(\mathbb{R}^2)} \lesssim 1$ to close the estimate. Via scaling \eqref{scale-fullNSC-to-kappa}, this is equivalent to
\begin{equation*}
	||\overline{u_0}||_{L^2(\mathbb{R}^2)} \lesssim \delta \nu^{\frac{1}{2}},
\end{equation*}
which degenerates as $\nu \to 0$. Hence, Theorem \ref{Thm-abbrev} complements \cite{CDGG02} for the cases $||\overline{U_0}||_{L^2(\mathbb{R}^2)} \gtrsim \nu^{\frac{1}{2}}$ by rescaling to the model \eqref{eq:NSC} with small $\kap$.

\bigskip
\item In \cite{IT, KLT}, the authors already fixed the viscosity to $1$ and considered
\begin{equation*}
	\dt v - \Delta v + (v\cdot\nabla)v + \Omega e_3 \times v + \nabla q = 0,
\end{equation*}
in the more generic setting where there is no assumption on the form of initial data and it only needs to lie in a Sobolev space. This equation can be obtained from a different scaling
\begin{equation*}
	v(t,x) = \sqrt{\kap^{-1} \Omega} \cdot u(\Omega t, \sqrt{\kap\Omega} x), \quad q(x) = \kap^{-1}\Omega \cdot p(\sqrt{\kap\Omega} x).
\end{equation*}
For $s \in (\frac{1}{2}, \frac{9}{10})$, it was shown that if $||v_0||_{\dot{H}^s} \lesssim \Omega^{\frac{1}{2}(s-\frac{1}{2})}$, then $v(t,x)$ is globally defined. $\dot{H}^s$ norms of $u$ and $v$ are related by $\Vert v \Vert_{\dot{H}^s} = \kap^{\frac{1}{2}(s-\frac{5}{2})} \Omega^{\frac{1}{2}(s-\frac{1}{2})} \Vert u \Vert_{\dot{H}^s}$, and from this, one can observe that the Theorem \ref{Thm-abbrev} allows initial data with bigger $\dot{H}^s$ norm compared to \cite{KLT} as we include $\kap \ll 1$.

\end{itemize}

\medskip

The examples above clearly show there has been a lack of control on $\epsilon(\nu,\delta)$ as $\nu \to 0$. Via scaling \eqref{scale-fullNSC-to-kappa}, $u(t,x) = \delta U(\delta t, x)$, Theorem \ref{Thm-abbrev} indicates that initial data of size $||U_0|| \leq \delta^{-1} \epsilon^*$ generate global solution, no matter how close $\nu$ is to $0$. In the threshold terminology, we improved \eqref{best-results} and established
\begin{equation*}
	\varliminf_{\nu \to 0} \epsilon(\nu,\delta) \gtrsim \delta^{-1}.
\end{equation*}

We note that, as in many previous literatures, we imposed a structural assumption on the initial data. Here, the proof relies on the fact that the velocity is axisymmetric. The assumption is inherited from \cite{GHPW, GPW} as we have used similar framework. Recently, Ren and Tian \cite{RT} announced to have removed axisymmetry condition for Euler-Coriolis system which could be an interesting direction to apply to Navier-Stokes equation as well.

\subsection{Overview of the techniques and plan of the article}
As can be inferred from the history of the problem and known results, the most subtle part in the proof is to prevent $\kap$ from shrinking the allowed size of initial data. The dissipation term $\kap\Delta u$ has always been used to construct global solutions, but comes with the price of limiting the threshold in terms of $\kap$. As such, we will rely less on the decay properties of the heat semi-group, but more on the dispersion $\Lambda$ given by the Coriolis term.(see Section 2 for the definition of $\Lambda$) To this end, we inherit the frameworks and a number of techniques from the Euler-Coriolis case \cite{GPW}, where the authors exploited dispersive effect to a great extent. In this sense, this paper lies in the stream of works in fluid dynamics, for example, \cite{DIPP, GMS12, GMS15, GIP, GPW, IP15, IP18, PW}, that incorporated spacetime resonance method of combining vector field and normal form methods \cite{Germain} to analyze nonlinear dispersive equations.

Besides from the technical modifications, one notable difference between Euler-Coriolis and Navier-Stokes-Coriolis equations is the number of symmetries. Euler-Coriolis had two symmetries - scaling and rotation around $z$-axis - giving rise to $S$- and $\Omega$-vector fields respectively. However, Navier-Stokes-Coriolis equation loses the scaling symmetry due to the introduction of the Laplacian term. This is seemingly a big issue, as the role of $S$-vector field is crucial in \cite{GPW}. Having said that, using the fact that the commutator between $S$ and the Laplacian gives `lower-order' operator, we will show throughout the paper, for example in Section 4, how this lack of scaling symmetry can be managed.

Another difference arises in the normal form. From integration by parts in time, there are two sources where normal form can be beneficial: discrepancies between $\Lambda(\zeta)$'s or between $\kap|\zeta|^2$'s, $\zeta \in \{\xi,\xi-\eta,\eta\}$. It's interesting to see that, due to the different nature of dispersion and dissipation, they never compete, but always help. It's exactly this fact that enables us to use only the differences between $\Lambda(\cdot)$'s to prevent $\kap$ from being used in the estimates. Hence, by modifying the normal form in a suitable way (see Corollary \ref{Cor-normalform} or the entire Section 5.2), we will be able to exploit the lack of spacetime resonance in a way that does not involve $\kap$. \\

Having summarized the differences, we will use similar settings and techniques used in analyzing the dispersive-only version of \eqref{eq:NSC}, the Euler-Coriolis equation. Namely, two vector fields $S$ and $\Omega$ stemming from the symmetries of Euler-Coriolis, the use of variables $\U$ that diagonalize the linearized equation, localizations $P_{k,p,q}$ reflecting the property of dispersion relation $\Lambda$, and the choice of norms from dispersive linear estimates. We will see that the framework fits quite smoothly with some proper modifications.\\

In section 2, some preliminary facts and necessary notations are introduced. In section 3, we present the main theorem in full detail and lay out the summary of its proof. Dispersive decay for Euler-Coriolis equations inherited from \cite{GPW} and the energy estimates for Navier-Stokes-Coriolis equation will be given in Section 4, and several methods to bound bilinear terms are explained in Section 5. We use these settings to bootstrap the $B$ norm in section 6, and the $X$ norm in section 7.

\section{Preliminaries}

As stressed above, since we want a result uniform in $\kap$, in particular, as $\kap \to 0$, a lot of the framework from the study of inviscid case comes in handy. We list here the main tools we inherit from \cite{GPW}.

\subsection{Dispersive unknowns}
At the heart of the analysis lies the fact that the Coriolis force brings dispersion into the system \cite{Dut2005, GPW, KLTEuler} with dispersion relation
\begin{equation}\label{DefofLambda}
	\Lambda(\xi) = \frac{\xi_3}{|\xi|}.
\end{equation}
The dispersion relation becomes visible once we write the system in $A:=|\nabla_h|^{-1} \text{curl}_h(u)$, and $C:= |\nabla||\nabla_h|^{-1} u_3$ where $\nabla_h = (\partial_1, \partial_2, 0)$, $\text{curl}_h u = \partial_1 u_2 -\partial_2 u_1$. From the divergence-free condition, \eqref{eq:NSC} is equivalent to
\begin{equation*}
	\begin{split}
		& \dt A - \kap\Delta A - i\Lambda C = -|\nabla_h|^{-1} \partial_j \partial_n \in^{j k} \{u_n u_k\} - i\Lambda|\nabla| \in^{j k} |\nabla_h|^{-1} \partial_j \{u_3 u_k\}, \\
		& \dt C - \kap\Delta C - i\Lambda C = -i\Lambda |\nabla|\sqrt{1-\Lambda^2} \{ |\nabla_h|^{-2} \partial_j\partial_k \{u_j u_k\} + u_3^2\} - |\nabla||\nabla_h|^{-1} \partial_j (1-2\Lambda^2)\{u_3 u_j\}.
	\end{split}
\end{equation*}
We use the same dispersive unknowns $U_{\pm} = A \pm C$, and the profiles $\U := e^{\mp it\Lambda}U_{\pm}(t)$ as in \cite{GPW} to use the axisymmetry and simplify the equation. Then, the equations become
\begin{equation}\label{ProfileEquation}
	\begin{split}
		\dt \mathcal{U}_+ - \kap\Delta \mathcal{U}_+ = & \mathcal{Q}_{\mathfrak{m}_+^{++}}(\mathcal{U}_+, \mathcal{U}_+)
		+ \mathcal{Q}_{\mathfrak{m}_+^{+-}}(\mathcal{U}_+, \mathcal{U}_-)
		+ \mathcal{Q}_{\mathfrak{m}_+^{--}}(\mathcal{U}_-, \mathcal{U}_-), \\
		\dt \mathcal{U}_- - \kap\Delta \mathcal{U}_- = & \mathcal{Q}_{\mathfrak{m}_-^{++}}(\mathcal{U}_+, \mathcal{U}_+)
		+ \mathcal{Q}_{\mathfrak{m}_-^{+-}}(\mathcal{U}_+, \mathcal{U}_-)
		+ \mathcal{Q}_{\mathfrak{m}_-^{--}}(\mathcal{U}_-, \mathcal{U}_-),
	\end{split}
\end{equation}
where
\begin{gather}
	\mathcal{Q}_{\mathfrak{m}_{\mu}^{\mu_1, \mu_2}}(\mathcal{U}_{\mu_1}, \mathcal{U}_{\mu_2})(t) := \mathcal{F}^{-1} \left( \int_{\mathbb{R}^3} e^{it\Phi_{\mu}^{\mu_1, \mu_2}} \mathfrak{m}_{\mu}^{\mu_1, \mu_2}(\xi,\eta) \widehat{\mathcal{U}_{\mu_1}}(\xi-\eta) \widehat{\mathcal{U}_{\mu_2}}(\eta) d\eta \right),
	\quad \mu,\mu_1,\mu_2 \in \{+,-\}, \\
	\Phi_{\mu}^{\mu_1, \mu_2}(\xi,\eta) = \mu\Lambda(\xi) + \mu_1 \Lambda(\xi-\eta) + \mu_2 \Lambda(\eta), \quad \mu, \mu_1, \mu_2 \in \{+,-\}, \label{def-of-Phi}
\end{gather}
and the multipliers are of the form
\begin{equation}
	\begin{split}
		\mathfrak{m}_{\mu}^{\mu_1, \mu_2}(\xi,\eta) & = |\xi|\mathfrak{n}_{\mu}^{\mu_1, \mu_2}(\xi,\eta), \\
		\mathfrak{n}_{\mu}^{\mu_1, \mu_2}(\xi,\eta) & \in \text{span}\left\{ \Lambda(\zeta_1)\sqrt{1-\Lambda^2(\zeta_2)} \cdot \bar{\mathfrak{n}}(\xi,\eta), \zeta_j \in \{\xi,\xi-\eta,\eta\}, \bar{\mathfrak{n}} \in \bar{E} \right\}, \mu, \mu_1, \mu_2 \in \{+,-\}, \\
		\bar{E} & := \left\{ \Lambda(\zeta), \sqrt{1-\Lambda^2(\zeta)}, \frac{\xi_h \cdot \theta_h}{|\xi_h||\theta_h|}, \frac{\xi_h^\perp \cdot \theta_h}{|\xi_h||\theta_h|} : \zeta \in \{\xi,\xi-\eta,\eta\}, \theta \in \{\xi-\eta, \eta\} \right\}.
	\end{split}
\end{equation}
See \cite[Lemma 2.3]{GPW} for the derivation; although we have $\kappa\Delta u$ in addition, since this is a linear term, one can derive \eqref{ProfileEquation} by modifying the proof there in an obvious way. In practice, when the signs do not matter, we will omit $\mu, \mu_1, \mu_2$ and simply write $\Phi, \mathfrak{m}$ and $\mathcal{U}$ or $\U$. For example, the Duhamel's formula becomes
\begin{gather}
	\U(t) = e^{t\kap\Delta}\U(0) + \sum_{\mathfrak{m}} \mathcal{B}_\mathfrak{m}(\U, \U)(t), \\
	\mathcal{B}_\mathfrak{m}(f_1, f_2)(t) := \int_0^t e^{(t-s)\kap\Delta} \mathcal{Q}_\mathfrak{m}(f_1, f_2)(s) ds. \label{def-bilinear-term}
\end{gather}

\subsection{Vector fields}
Next, we introduce the vector fields:
\begin{gather*}
	S\mathbf{v} = (x\cdot\nabla)\mathbf{v} - \mathbf{v}, \quad Sf = x\cdot\nabla f - 2f, \\
	\Omega\mathbf{v} = (x_1\partial_2 - x_2\partial_1)\mathbf{v} - \mathbf{v}, \quad \Omega f = (x_1\partial_2 - x_2\partial_1)f, \\
	\Upsilon f = \partial_{\phi} f.
\end{gather*}
$S$ and $\Omega$ come from the symmetries of the Euler-Coriolis equation, and since $S = \rho\partial_\rho$ and $\Omega = \partial_\theta$ in spherical coordinates up to lower order terms, $\Upsilon$ is introduced to complete the space of derivatives. Note that $\Omega$ commutes with the Navier-Stokes-Coriolis equation and hence axisymmetry is preserved as the equation has rotational symmetry, but $S$ does not. However, it misses by a little: if we apply $S$ to \eqref{eq:NSC}, $u$ satisfies
\begin{equation*}
	\dt Su - \kap\Delta (S-2)u + (Su \cdot \nabla)u + (u\cdot\nabla)Su + e_3 \times Su + \nabla Sp = 0.
\end{equation*}
The commutator is a lower order term and does not create an issue in the analysis; see the energy estimates in Section 4 or normal form in Section 5.2, for example. Hence, the $S$-vector field will still be useful.

\subsection{Localizations}
Let $\psi \in C_0^{\infty}(\mathbb{R}; [0,1])$ be a standard bump function, i.e. an even function with $\text{supp } \psi \subset [-2,2]$, decreasing in $\mathbb{R}_+$, and $\psi|_{[-1,1]} \equiv 1$. Put $\varphi(x) = \psi(x) - \psi(2x)$ and define $P_{k}, P_{k,p}$, and $P_{k,p,q}$ as
\begin{equation*}
	\mathcal{F}(P_k f)(\xi) = \varphi_k(\xi)\hat{f}(\xi), \quad \mathcal{F}(P_{k,p} f)(\xi) = \varphi_{k,p}(\xi)\hat{f}(\xi), \quad \mathcal{F}(P_{k,p,q} f)(\xi) = \varphi_{k,p,q}(\xi)\hat{f}(\xi),
\end{equation*}
for $k \in \mathbb{Z}, p,q \in \mathbb{Z}_-, \mathbb{Z}_- = \{n\in\mathbb{Z}, n\leq0\}$ where
\begin{equation*}
	\varphi_k(\xi) = \varphi(2^{-k} |\xi|), \quad
	\varphi_{k,p}(\xi) = \varphi_k(\xi) \varphi(2^{-p} \sqrt{1 - \Lambda^2(\xi)}), \quad
	\varphi_{k,p,q}(\xi) = \varphi_{k,p}(\xi) \varphi(2^{-q} \Lambda(\xi)).
\end{equation*}
$P_k$ is the standard Littlewood-Paley localizations based on the size of frequency, and the additional localizations $P_{k,p}$ and $P_{k,p,q}$ are introduced to reflect the anisotropy present in the dispersion relation $\Lambda$.
Another localization to measure the smoothness in angular direction will be needed. Using the Legendre polynomial
\begin{equation*}
	\mathfrak{Z}_n(x) = \frac{2n+1}{4\pi}L_n(x), \quad L_n(x) = \frac{1}{2^n (n!)} \frac{d^n}{dz^n}[(z^2-1)^n],
\end{equation*}
we define the angular localizations $\bar{R}_{\leq l}$ and $\bar{R}_l$ by
\begin{equation}\label{Def-of-Rl}
	\begin{split}
		\bar{R}_{\leq l}f(x) & = \sum_{n \geq 0} \psi(2^{-l}n) \int_{\mathbb{S}^2} f(|x|\theta) \mathfrak{Z}_n(\langle \theta, \frac{x}{|x|} \rangle) dS(\theta), \\
		\bar{R}_{l}f(x) & = \sum_{n \geq 0} \varphi(2^{-l}n) \int_{\mathbb{S}^2} f(|x|\theta) \mathfrak{Z}_n(\langle \theta, \frac{x}{|x|} \rangle) dS(\theta),
	\end{split}
\end{equation}
for $l \geq 0$ where $dS$ is the usual measure on the sphere. We refer to \cite{GPW} or the Appendix for the properties of $R_l$'s, but note here that they measure the size of the angular derivatives in the sense that
\begin{equation}
	\sum_{1 \leq a < b \leq 3} ||\Omega_{a b} \bar{R}_l f ||_{L^r} \simeq 2^l ||\bar{R}_l f||_{L^r}, \quad \Omega_{a b} = x_a\partial_b - x_b\partial_a,
\end{equation}
and that due to the `uncertainty principle'
\begin{equation*}
	||[\Omega_{a 3}, P_{k,p}]||_{L^r \to L^r} \lesssim 2^{-p},
\end{equation*}
we only need to consider the case $l + p \geq 0$, and hence define
\begin{equation*}
	R_l = \begin{cases}
		0, & p+l < 0, \\
		\bar{R}_{\leq l}, & p+l = 0, \\
		\bar{R}_l, & p+l > 0.
	\end{cases}
\end{equation*}
The dependence on $p$ exists, but will be obvious in practice and hence we neglect it in the notation.\\
\\
We will often need three different localizations - two for the input functions of bilinear terms, and one for the output. Localization indices of the output functions will be denoted by $k,p,q$, and $l$, while those of the input functions will be denoted by $k_j,p_j,q_j$, and $l_j$, $j=1,2$. For example, typical localized bilinear term under analysis will take the form
\begin{equation*}
	P_{k,p} \mathcal{Q}_\mathfrak{m}(P_{k_1,p_1} f_1, P_{k_2,p_2} f_2) = \mathcal{F}^{-1}\left\{ \int_{\mathbb{R}^3} e^{is\Phi(\xi,\eta)} \mathfrak{m}(\xi,\eta) \varphi_{k,p}(\xi) \varphi_{k_1,p_1}(\xi-\eta) \varphi_{k_2,p_2}(\eta) \hat{f}_1(\xi-\eta) \hat{f}_2(\eta) d\eta \right\},
\end{equation*}
with possibly additional localizations in $q, q_j$ or $l, l_j$'s. For this reason, it will be convenient to introduce
\begin{equation*}
	\chi_h(\xi,\eta) := \varphi_{k,p}(\xi) \varphi_{k_1,p_1}(\xi-\eta) \varphi_{k_2,p_2}(\eta), \quad \chi(\xi,\eta) = \varphi_{k,p,q}(\xi) \varphi_{k_1,p_1,q_1}(\xi-\eta) \varphi_{k_2,p_2,q_2}(\eta).
\end{equation*}
It is important to note that the localizations are not independent. From the definition of $R_l$, we implicitly assumed that $p+l$ will be nonnegative whenever we localize in both sizes, but even earlier than that, one can check that $2^{2p} + 2^q \simeq 1$ by definition. The same holds for $p_j,q_j,l_j$'s, $j=1,2$. Another fact of great importance is that indices across input functions and output function also cannot be independent, simply due to $\xi = (\xi-\eta) + \eta$. This will be investigated in more detail in Section 5.3. As such, we use
\begin{equation*}
	k_{\max} := \max\{k, k_1, k_2\}, \quad k_{\min} := \min\{k, k_1, k_2\},
\end{equation*}
since we will frequently need to compare their sizes. The same type of notations will be used for $p,q$, and $l$'s.

\section{Precise main theorem}

The main norms that we will bootstrap are the following:
\begin{equation}
	\begin{split}
		& ||f||_B := \sup_{k,p,q \in\mathds{Z}, p,q \leq 0} 2^{3k^+ - \frac{1}{2}k^-} 2^{-p-\frac{q}{2}} ||P_{k,p,q}f||_2, \\
		& ||f||_X := \sup_{\substack{k,p \in \mathds{Z}, p \leq 0 \leq l, \\ p+l \geq 0}} 2^{3k^+} 2^{(1+\beta)l + \beta p} ||P_{k,p} R_l f||_2,
	\end{split}
\end{equation}
where $k^+ = \max\{k,0\}$ and $k^- = \min\{k,0\}$. Now, we can state Theorem \ref{Thm-abbrev}(a) rigorously. Recall from the comment following the theorem that it is enough to consider only the small viscosities.
\begin{theorem}\label{MainThmDetailed}
	There exist $M, N, N_0, \beta, \epsilon >0$ with $N_0 \gg M \gg N, \beta^{-1}$ and $\beta^{-1} \geq 60$, all of which are independent of $\kap \in [0,1]$, such that if the dispersive unknowns $U^{\pm}_0$ of $u_0$ satisfy
	\begin{equation}\label{InitialCondition}
		\begin{split}
			||U_0^{\pm} ||_{H^{2N_0} \cap \dot{H}^{-1}} + ||S^a U_0^{\pm}||_{L^2 \cap \dot{H}^{-1}} \leq \epsilon, \quad & 0 \leq a \leq M, \\
			||S^b U_0^{\pm}||_B + ||S^b U_0^{\pm}||_X \leq \epsilon, \quad & 0 \leq b \leq N,
		\end{split}
	\end{equation}
	then there exists a global solution $u \in C(\mathbb{R}, \mathbb{R}^3)$ to the Cauchy problem \eqref{eq:NSC}.
\end{theorem}

The main theorem will be proved by showing the following bootstrap argument.
\begin{theorem}\label{MainThmBootstrap}
	Suppose that $u$ is a solution to \eqref{eq:NSC} up to some time $T > 0$ with initial data satisfying \eqref{InitialCondition}. There exists $c > c' > 1$ such that if
	\begin{equation}\label{BootstrapAssumption}
		\begin{split}
			|| S^n \mathcal{U}_\pm(t) ||_{B} + ||S^n \mathcal{U}_\pm(t)||_X \leq c\epsilon, & \hspace{4mm} 0 \leq n \leq N,
		\end{split}
	\end{equation}
	holds for all $t \in [0,T]$, then actually
	\begin{equation*}
		|| S^n \mathcal{U}_\pm(t) ||_{B} + ||S^n \mathcal{U}_\pm(t)||_X \leq c'\epsilon, \hspace{4mm} 0 \leq n \leq N,
	\end{equation*}
	for $t \in [0,T]$. Moreover, the estimates are uniform for $0 \leq \kap \lesssim 1$.
\end{theorem}

\medskip

\begin{proof}[Proof of Theorem \ref{Thm-abbrev}(b)]
	We prove the inequality for solutions $u_1, u_2$ of \eqref{eq:NSC} with viscosities $\kappa_1, \kappa_2$, respectively. Subtracting two equations and taking inner product with $u_1-u_2$, we have
	\begin{equation}
		\frac{1}{2}\dt ||u_1-u_2||_2^2 + \kap_2 ||u_1-u_2||_{\dot{H}^1}^2 \leq |\kap_1-\kap_2| ||u_1-u_2||_{\dot{H}^1} ||u_1 ||_{\dot{H}^1} + ||\nabla u_1||_\infty ||u_1-u_2||_2^2.
	\end{equation}
	Putting $v = u_1-u_2$ and assuming $\kap_1 \leq \kap_2$ WLOG, we can absorb the ${\dot{H}^1}$ norm of $v$ on the right-hand side to deduce
	\begin{equation}
		\dt ||v||_2^2 \leq |\kap_1 - \kap_2| ||u_1||_{\dot{H}^1}^2 + 2||\nabla u_1||_\infty ||v||_2^2.
	\end{equation}
	Since $v(0) \equiv 0$, Gronwall's inequality gives us
	\begin{equation}
		||v(t)||_2^2 \leq \int_0^t |\kap_1 - \kap_2| ||u_1||_{\dot{H}^1}^2 \cdot e^{\int_s^t 2||\nabla u_1(\tau)||_\infty d\tau} ds.
	\end{equation}
	From the dispersive decay and the energy estimates in Corollary \ref{energy-estimate} and the dispersive decay Proposition \ref{dispersive-decay} and the bootstrap Thoerem \ref{MainThmBootstrap}, we get
	\begin{equation}
		||v(t)||_2^2 \lesssim |\kap_1 - \kap_2| \int_0^t \epsilon^2 \langle s \rangle^{2C\epsilon} \left(\frac{t}{s}\right)^{C\epsilon} ds \lesssim |\kap_1 - \kap_2| \epsilon^2 t^{1+C\epsilon}.
	\end{equation}
	Scaling back, we obtain \eqref{finite-time-convergence}.
\end{proof}

The bulk of the paper, namely the last two sections, is devoted to proving Theorem \ref{MainThmBootstrap}. The sections preceding them are linear estimates and the propagation of the energy type norms in the bootstrap assumption \eqref{InitialCondition} (Section 4), and the technical lemmas (Section 5) needed to bound the bilinear terms \eqref{def-bilinear-term} in Duhamel's formulas for the last two sections. Although the linear terms $e^{t\kap\Delta}\U(0)$ in Duhamel's formulas do not commute with the $S$-vector field since it does not commute with the equation, these are harmless, as
\begin{equation}\label{S-heat-commutator}
	[S, e^{t\kap\Delta}] = 2t\kap\Delta e^{t\kap\Delta}
\end{equation}
has an operator norm bounded by an absolute constant. One can check that more iterations of $S$-vector fields do not generate qualitative change, so that we only need to choose $c$ and $c'$ larger than sum of these types.

Hence, the challenge is to deal with the bilinear term. Details of the estimates will be provided later, but we describe here how we reduce the problem into estimating atoms. First, we cite \cite{GPW} to claim that
\begin{equation}
	S_\xi\mathcal{F}\left(\mathcal{Q}_{\mathfrak{m}}(\U,\U)\right)(\xi) = \mathcal{F}\left(-2\,\mathcal{Q}_{\mathfrak{m}}(\U,\U) +\mathcal{Q}_{\mathfrak{m}}(S\U,\U) + \mathcal{Q}_{\mathfrak{m}}(\U,S\U)\right)(\xi).
\end{equation}
Together with \eqref{S-heat-commutator}, the bilinear terms with some copies of $S$-vector field have the form
\begin{equation}
	S^n \mathcal{B}_\mathfrak{m}(\U, \U) = \sum_{a+b+c \leq n} \int_0^t h_a((t-s)\kap\Delta) e^{(t-s)\kap\Delta} \mathcal{Q}_\mathfrak{m}(S^b\U, S^c\U)(s) ds,
\end{equation}
where $h_a(\cdot)$ is a polynomial of order $a$. After writing this as a sum of localized pieces, we use energy estimates together with the bootstrap assumptions to reduce the problem into estimating finitely many localized pieces. In addition to localizations in space introduced already, we also use time localizations $\tau_m(s), m=0,1,\cdots, L+1$, $L = \lceil \log_2(t+2) \rceil$ where
\begin{equation}
	\text{supp }\tau_m \subset [2^{m-1}, 2^{m+1}], m=1, \cdots, L, \quad \text{supp }\tau_0 \subset [0,2], \quad \text{supp }\tau_{L+1} \subset [t-2,t], \quad \sum_{m=0}^{L+1} \tau_m(s) \equiv 1,
\end{equation}
and $\int_0^t \tau'm(s) ds \lesssim 1$. After these reductions, most of the time we will be working on $(\log \langle t\rangle)$-many pieces of
\begin{equation}
	P_{k,p(,q)} (R_l) \int_0^t \tau_m(s) h_a((t-s)\kap\Delta) e^{(t-s)\kap\Delta} \mathcal{Q}_\mathfrak{m}(P_{k_1,p_1(,q_1)} (R_{l_1}) S^b\U, P_{k_2,p_2(,q_2)} (R_{l_2}) S^c\U)(s) ds,
\end{equation}
and show that their sum is bounded by $c'\epsilon$. Parentheses around $q, q_j, R_l, R_{l_j}$s denote that we will use finer localizations only when they are needed, not from the outset. In Section 2.3, it was mentioned that the localizations are not independent. As such, we will divide cases according to their relations explained in Section 5.3, and each configuration will require different methods.

\section{Dispersive decay and energy estimates}

One of the most important estimates from having the dispersion relation is the dispersive decay estimate, and indeed we inherit and exploit this estimate proved in \cite[Proposition 4.1]{GPW}.
\begin{proposition}\label{dispersive-decay}
	Define the norm $D$ as $$ ||f||_D := \sup_{\substack{0 \leq a \leq 3 \\ 0 \leq b \leq 2}} (||S^a f||_B + ||S^b f||_X).$$ Then, we can split $$P_{k,p,q} e^{it\Lambda}f = I_{k,p,q}(f) + II_{k,p,q}f$$ such that for any $0 < \beta' < \beta$,
	\begin{equation}\label{linear-decay}
		\begin{split}
			||I_{k,p,q}(f)||_{\infty} & \lesssim 2^{\frac{3}{2}k - 3k^+} \cdot \min\{2^{2p+q}, 2^{-p-\frac{q}{2}} t^{-\frac{3}{2}}\} ||f||_D, \\
			||II_{k,p,q}(f)||_{2} & \lesssim 2^{-3k^+} t^{-1-\beta'} 2^{(-1-2\beta')p} \cdot \boldsymbol{1}_{2^{2p+q} \gtrsim t^{-1}} ||f||_D.
		\end{split}
	\end{equation}
\end{proposition}
\begin{proof}
	See \cite[Proposition 4.1]{GPW}.
\end{proof}

As a first application, we can prove slow energy growth Corollary \ref{energy-estimate} for \eqref{eq:NSC}.
\begin{lemma}
	Put $c_{n,k}$ as the $k$th coefficient of $(x+\frac{3}{2})(-1)^{n-k}(x+3)^{n-k-1}(x-2)^n$ when $k < n$, and $c_{n,n} = 1$. Then,
	\begin{equation}
		- \langle S^n u, S^n \Delta u \rangle = \sum_{k=0}^n c_{n,k} ||\nabla S^k u||^2_{L^2}
	\end{equation}
\end{lemma}

\begin{proof}
	From $S\Delta u = \Delta (S-2)u$, we can first put $\langle S^n u, S^n \Delta u \rangle = \langle S^n u, \Delta (S-2)^n u \rangle$. Using the fact that $S^* = -S-5$ for vectors,
	\begin{equation*}
		\langle S^{n+1} u, \Delta S^ku\rangle = \langle S^nu, (-S-5)\Delta S^ku \rangle = -\langle S^nu, \Delta((S-2)S^ku + 5S^ku) \rangle = -\langle S^nu, \Delta (S+3)S^ku \rangle
	\end{equation*}
	Also, we have $\langle Sv, \Delta Sv \rangle = - || \nabla Sv||^2_{L^2}$ and
	\begin{equation*}
		\begin{split}
			\langle Sv, \Delta v \rangle &= \sum_{i,j,k} \int x_j\partial_j v_i \partial_k\partial_k v_i dx - \int v \cdot \Delta v dx \\
			& = - \sum_{i,j,k} \int [ x_j \partial_j \partial_k v_i \partial_k v_i + \partial_k v_i \partial_k v_i ]dx - \int v \cdot \Delta v dx \\
			& = - \sum_{i,j} \int x_j \partial_j |\nabla v_i|^2/2 dx = \frac{3}{2}||\nabla v||^2_{L^2}
		\end{split}
	\end{equation*}
	Combining these 3 facts, we can write $\langle S^n u, \Delta(S-2)^n u \rangle$ as a linear combination of $||\nabla S^ku||^2_{L^2}, 0\leq k\leq n$. Except for $||\nabla S^n u||^2_{L^2}$ and $||\nabla S^{n-1} u||^2_{L^2}$,  the coefficients will be determined by $\langle S^k u, \Delta S^k u \rangle = - || \nabla S^k u||^2_{L^2}$ and $\langle S^{k+1} u, \Delta S^k u \rangle = \frac{3}{2} || \nabla S^k u||^2_{L^2}$ after sending $S$ several times from left to right in the inner product $\langle S^n u, \Delta S^a u\rangle$. Using these 3 rules with the fact that starting coefficients are $\dbinom{n}{k}(-2)^k$, the lemma is proved.
\end{proof}

Similarly, $\langle S^n u, (S-2)^n u\rangle$ can also be written as a linear combination $\overset{n}{\underset{k=0}{\sum}} c'_{n,k} ||S^k u||_2^2$ for some $c'_{n,k}$'s where $c'_{n,n} = 1$. Using these facts, we can get bounds for $S^n u$'s.

\begin{lemma}
	There exist $a_{n,k}, a'_{n,k} \geq 0, a_{n,n}=a'_{n,n}=1$ such that
	\begin{equation}
		\sum_{k=0}^n a_{n,k} (||S^k u(t)||^2_2 - ||S^k u_0||^2_2) + 2\kappa \int_0^t ||\nabla S^n u(s)||^2_2 ds \lesssim \int_0^t (||u(s)||_\infty + ||\nabla u(s)||_\infty)(\sum_{k=0}^n ||S^k u(s)||_2)^2 ds
	\end{equation}
	\begin{equation}
		\begin{split}
			\sum_{k=0}^n a'_{n,k} (|| |\nabla|^{-1} S^k u(t)||^2_2 - || |\nabla|^{-1} S^k u_0||^2_2) & + 2\kappa \int_0^t || S^n u(s)||^2_2 ds \\ \lesssim & \int_0^t ||u(s)||_\infty || |\nabla|^{-1} S^n u(s)||_2 (\sum_{k=0}^n || S^k u(s)||_2) ds
		\end{split}
	\end{equation}
\end{lemma}

\begin{proof}
	Since $S$ commutes with the Euler-Coriolis equation, iterating $S$, we obtain
	\begin{equation}\label{energy1}
		\dt S^nu - \kappa S^n \Delta u + e_3 \times S^nu + \partial_j \{S^n\{u_i u_j\}\} + \nabla S^n p = 0.
	\end{equation}
	Taking inner product with $S^n u$ and arranging in the same way as \cite{GHPW},
	\begin{equation*}
		\frac{1}{2}\dt||S^n u||^2_{L^2} - \kappa \langle S^n u, S^n \Delta u \rangle \lesssim ||\nabla u||_{L^\infty}||S^n u||^2_{L^2} + ||S^n u||_{L^2} ||u||_{L^\infty} (\sum_{a=0}^{n}||S^a u||_{L^2})
	\end{equation*}
	We now focus on the dissipation term. In particular, the highest order term will always be $\kappa ||\nabla S^n u||^2_{L^2}$. So we induct on $n$. For $n=0$, the usual energy estimate gives,
	\begin{equation*}
		||u(t)||^2_{L^2} - ||u_0||^2_{L^2} + 2\kappa\int_0^t ||\nabla u(s)||^2_{L^2} ds = 0
	\end{equation*}
	Suppose the lemma holds until $n-1$. Apply Lemma 4.2 and put $d_{n,k} = \min\{0, c_{n,k}\}$ for $k < n$ and $d_{n,n} = 1$. Then,
	\begin{equation*}
		\begin{split}
			||S^n u(t)||^2_2 - ||S^n u_0||^2_2 &-2\kappa \int_0^t \langle S^n u, S^n \Delta u \rangle ds = ||S^n u(t)||^2_2 - ||S^n u_0||^2_2 +2\kappa \int_0^t \sum_{k=0}^n c_{n,k} ||\nabla S^k u(s)||^2_2 ds \\
			& \geq ||S^n u(t)||^2_2 - ||S^n u_0||^2_2 +2\kappa \sum_{k=0}^n d_{n,k} \int_0^t ||\nabla S^k u(s)||^2_2 ds \\
			& = ||S^n u(t)||^2_2 - ||S^n u_0||^2_2 +2\kappa\int_0^t ||\nabla S^n u(s)||^2_2 ds \\
			& \hspace{2mm} - \sum_{k=0}^{n-1} d_{n,k} \left[ \sum_{j=0}^k a_{k,j} (||S^j u(t)||^2_2 - ||S^k u_0||^2_2) + \int_0^t \langle S^k u(s), S^k( (u\cdot\nabla)u(s) ) \rangle ds \right] \\
			& = ||S^n u(t)||^2_2 - ||S^n u_0||^2_2 +2\kappa\int_0^t ||\nabla S^n u(s)||^2_2 ds \\
			& \hspace{2mm} + \sum_{j=0}^{n-1}\sum_{k=j}^{n-1} |d_{n,k}|a_{k,j} (||S^j u(t)||^2_2 - ||S^k u_0||^2_2) - \sum_{k=0}^{n-1} d_{n,k} \int_0^t \langle S^k u(s), S^k( (u\cdot\nabla)u(s) ) \rangle ds
		\end{split}
	\end{equation*}
	Hence, by defining $a_{n, j} = \sum_{k=j}^{n-1} |d_{n,k}|a_{k,j}$ for $j < n$ and dealing with the nonlinearity as in \cite{GHPW}, we obtain the first inequality for $n$, completing the induction.
	
	Proof for the $\dot{H}^{-1}$ norms works the same, since from \eqref{energy1}, we can also get
	\begin{equation*}
		\frac{1}{2}\dt|| |\nabla|^{-1} S^n u||^2_{L^2} + \kappa \langle S^n u, (S-2)^n u \rangle \lesssim || |\nabla|^{-1} S^n u||_{L^2} ||S^n(u \otimes u)||_2 \lesssim || |\nabla|^{-1} S^n u||_{L^2} ||u||_{L^\infty} \sum_{a=0}^{n}||S^a u||_{L^2}
	\end{equation*}
	Hence, using the similar induction process as before, by putting $a'_{n,k} = \sum_{j=k}^{n-1} |\min\{c'_{n,j}, 0\}| a'_{k,j}$, we get the second inequality.
\end{proof}

\begin{lemma}
	If $u$ is a solution to \eqref{eq:NSC}, it satisfies
	\begin{equation}
		||u(t)||^2_{H^m} + 2\kappa\int_0^t ||\nabla u(s)||^2_{H^m} ds - ||u_0||^2_{H^m} \lesssim \int_0^t ||\nabla u(s)||_{L^\infty} ||u(s)||^2_{H^m} ds
	\end{equation}
\end{lemma}

\begin{proof}
	The proof is almost identical to the one in \cite{GHPW}. Taking derivatives and then having dot product to \eqref{eq:NSC}
	\begin{equation*}
		\begin{split}
			\frac{1}{2}\dt||\partial^{\mu}u||^2_2 - \kappa\int\Delta\partial^{\mu}u^{\alpha} \partial^{\mu}u^{\alpha} dx &= - \int \partial^{\mu}u^{\alpha} \partial^{\mu} \partial^{\beta}(u^{\alpha}u^{\beta}) dx \\
			& = - \int \partial^{\mu}u^{\alpha} \partial^{\mu} ((\partial^{\beta}u^{\alpha})u^{\beta}) dx \hspace{8mm} (\because \partial^{\beta}u^{\beta} = 0) \\
			& = - \int \partial^{\mu}u^{\alpha} (\partial^{\mu-e_i} \partial^{\beta}u^{\alpha}) \partial_i u^{\beta} dx \hspace{3mm} (\because  \int \partial^{\mu}u^{\alpha} (\partial^{\mu} \partial^{\beta}u^{\alpha})u^{\beta} dx = 0) \\
			& \lesssim ||\nabla u||_{L^\infty} ||u||_{\dot{H}^{|\mu|}}^2
		\end{split}
	\end{equation*}
	Hence, integrating with respect to $t$ and adding up the result for $|\mu| \leq m$, gives the lemma.
\end{proof}

\begin{corollary}\label{energy-estimate}
	Suppose $u$ is a solution to \eqref{eq:NSC} with initial data satisfying \eqref{InitialCondition}. If $u$ satisfies \eqref{BootstrapAssumption}, up to some time $T > 0$. Then,
	\begin{equation*}
		||U^\pm(t)||_{H^{2N_0} \cap \dot{H}^{-1}} + ||S^n U^\pm(t)||_{L^2 \cap \dot{H}^{-1}} \leq \epsilon \jpt^{C\epsilon}, \hspace{10mm} 0 \leq n \leq M.
	\end{equation*}
\end{corollary}

\begin{proof}
	Using the previous lemmas, the proof is complete once we show
	\begin{equation*}
		(1+t)||u(t)||_\infty + (1+t)||\nabla u(t)||_\infty \lesssim C\epsilon,
	\end{equation*}
	which can be achieved in exactly the same way as in the proof of \cite[Corollary 7.2]{GPW}, as we again imposed initial conditions and bootstrap assumptions on the dispersive unknowns and profiles.
\end{proof}

\section{Bounds on bilinear terms}

Two major ways of controlling the bilinear terms are integration by parts and normal form. We first set the framework for integration by parts along vector fields. The crucial finding in \cite{GPW} showing that there's no spacetime resonance was through the quantity
\begin{equation}
	\bar{\sigma}(\xi,\eta) := \xi_3 \eta_h - \eta_3 \xi_h = -(\xi \times \eta)_h^{\perp}
\end{equation}
defined in the Fourier space.
\begin{proposition}\label{no-spacetime-resonance}
	Assume that $|\Phi| \leq 2^{q_{\max}-10}$. Then $2^{p_{\max}} \sim 1$, and $|\bar{\sigma}| \gtrsim 2^{q_{\max}} 2^{k_{\max} + k_{\min}}$. Also,
	\begin{equation}
		|S_\eta \Phi| + |\Omega_\eta \Phi| \sim \frac{|\xi_h -\eta_h|}{|\xi-\eta|^3} |\bar{\sigma}(\xi,\eta)|.
	\end{equation}
\end{proposition}
$\Phi = \Phi(\xi,\eta)$ is defined at \eqref{def-of-Phi}. Recall the comment following the definition that we will not specify the signs when they are irrelevant.
\begin{proof}
	See \cite[Lemma 5.1, Proposition 5.2]{GPW}.
\end{proof}
This proposition implies that when we don't perform normal form due to the lack of lower bound on the phase, we can integrate by parts along either $S$ or $\Omega$ vector field thanks to the following lemma.

\medskip

\subsection{Bounds for repeated integration by parts}

The followings are the bounds for repeated integration by parts along vector fields which will be used frequently.
\begin{lemma}\label{IbyP}
	(1) Assume the localization parameters are such that $|\bar{\sigma}| \cdot \chi_h \gtrsim 2^{k_{\max}+k_{\min}} 2^{p_{\max}} =: L_1$. Then for any $N \in \mathbb{N}$,
	\begin{equation*}
		\begin{split}
			|| \mathcal{F}\{ \mathcal{Q}_{\mathfrak{m} \chi_h}(R_{l_1}f_1, R_{l_2}f_2) \}||_\infty \lesssim & \hspace{1mm} ||\mathfrak{m} \chi_h||_\infty \cdot (s^{-1} 2^{-p_1 + 2k_1} L_1^{-1} [1 + 2^{k_2 - k_1 + l_1}])^N \\ & \cdot ||P_{k_1,p_1} R_{l_1} (1,S)^N f_1||_2 \cdot ||P_{k_2,p_2} R_{l_2} (1,S)^N f_2||_2.
		\end{split}
	\end{equation*}
	
	(2) Assume the localization parameters are such that $|\bar{\sigma}| \cdot \chi \gtrsim 2^{k_{\max}+k_{\min}} 2^{p_{\max} + q_{\max}} =: L_2$. Then for any $N \in \mathbb{N}$,
	\begin{equation*}
		\begin{split}
			|| \mathcal{F}\{ \mathcal{Q}_{\mathfrak{m} \chi}(R_{l_1}f_1, R_{l_2}f_2) \}||_\infty \lesssim & \hspace{1mm} ||\mathfrak{m} \chi||_\infty \cdot (s^{-1} 2^{-p_1 + 2k_1} L_2^{-1} [1 + 2^{k_2 - k_1}(2^{q_2 - q_1} + 2^{l_1})])^N \\ & \cdot ||P_{k_1,p_1,q_1} R_{l_1} (1,S)^N f_1||_2 \cdot ||P_{k_2,p_2,q_2} R_{l_2} (1,S)^N f_2||_2.
		\end{split}
	\end{equation*}
	Both estimates hold for the indices of $f_1$ and $f_2$ swapped by symmetry.
\end{lemma}
\begin{proof}
	See \cite[Lemma 5.7]{GPW}.
\end{proof}

In addition to this main integration by parts bounds, we will have occasions where we will need to use integrate by parts along the vertical direction, using
\begin{equation*}
	D_3^{\eta} := |\eta|\partial_{\eta_3}.
\end{equation*}
\begin{lemma}\label{vertical-variant}
	Assume that the localization parameters are such that $p_1,p_2 \ll 0$, $k_1 \sim k_2$, and $|D_3^{\eta} \Phi| \cdot \chi_h \gtrsim L_3$. Then, for any $N \in \mathbb{N}$,
	\begin{equation*}
		\begin{split}
			|| \mathcal{F}\{ \mathcal{Q}_{\mathfrak{m} \chi_h}(R_{l_1}f_1, R_{l_2}f_2) \}||_\infty \lesssim & \hspace{1mm} ||\mathfrak{m} \chi_h||_\infty \cdot (s^{-1} L_3^{-1} [1 + 2^{l_1 + p_1} + 2^{l_2 + p_2}])^N \\ & \cdot ||P_{k_1,p_1} R_{l_1} (1,S)^N f_1||_2 \cdot ||P_{k_2,p_2} R_{l_2} (1,S)^N f_2||_2.
		\end{split}
	\end{equation*}
\end{lemma}
\begin{proof}
	See \cite[Section 5.3.3]{GPW}.
\end{proof}

\medskip

\subsection{Normal Form}
When we have a good lower bound on the phase, which we can realize using Proposition \ref{no-spacetime-resonance}, we have another way to estimate the quadratic nonlinearities.
\begin{lemma}\label{normalform-base}
	For the solution $\mathcal{U}$ of \eqref{eq:NSC}, define $\mathcal{U}^m_+ (t), \mathcal{U}^m_- (t)$ as
	\begin{equation*}
		\widehat{\mathcal{U}}^m_\mu(t) = \sum_{\mu_1,\mu_2 \in \{+,-\}} \int_0^t \tau_m(s) e^{-(t-s)\kap|\xi|^2} \int_{\mathds{R}^3} e^{is\Phi^{\mu_1,\mu_2}_\mu} \mathfrak{m}^{\mu_1,\mu_2}_\mu(\xi,\eta) \widehat{\mathcal{U}}_{\mu_1}(s,\eta) \widehat{\mathcal{U}}_{\mu_2}(s,\xi-\eta) d\eta ds.
	\end{equation*}
	Then, there holds
	\begin{align*}
		\mathcal{U}^m_+(t) = \sum_{\mu_1,\mu_2 \in \{+,-\}} & \left\{ - \delta_{m,L} \mathcal{Q}_{\mathfrak{n}_+^{\mu_1,\mu_2}}(\mathcal{U}_{\mu_1},\mathcal{U}_{\mu_2})(t) + \delta_{m,0} e^{t\kap\Delta}[ \mathcal{U}^m_+(0) + \right. \mathcal{Q}_{\mathfrak{m}_+^{\mu_1,\mu_2}}(\mathcal{U}_{\mu_1},\mathcal{U}_{\mu_2})(0) ]  \\
		& + \int_0^t e^{(t-s)\kap\Delta} \tau_m(s) \mathcal{Q}_{\mathfrak{n}_+^{\mu_1,\mu_2}}((\dt-\kap\Delta)\mathcal{U}_{\mu_1},\mathcal{U}_{\mu_2})(s)ds \\
		& + \int_0^t e^{(t-s)\kap\Delta} \tau_m(s) \mathcal{Q}_{\mathfrak{n}_+^{\mu_1,\mu_2}}(\mathcal{U}_{\mu_1}, (\dt-\kap\Delta)\mathcal{U}_{\mu_2})(s)ds \\
		& + \left.\int_0^t e^{(t-s)\kap\Delta} \tau_m'(s) \mathcal{Q}_{\mathfrak{n}_+^{\mu_1,\mu_2}}(\mathcal{U}_{\mu_1},\mathcal{U}_{\mu_2})(s)ds \right\},
	\end{align*}
	\begin{equation*}
		\mathfrak{n}_+^{\mu_1,\mu_2} = -\frac{\mathfrak{m}_+^{\mu_1,\mu_2}}{i\Phi^{\mu_1,\mu_2}_\mu + 2\kap(\xi-\eta)\cdot\eta},
	\end{equation*}
	and similarly for $\mathcal{U}^m_-$.
\end{lemma}

\begin{proof}
	By computation,
	\begin{align*}
		& (\dt+\kap|\xi|^2) \int \tau_m(t) e^{it\Phi^{\mu_1,\mu_2}_\mu} \mathfrak{n}_+^{\mu_1,\mu_2}(\xi,\eta) \widehat{\mathcal{U}}_{\mu_1}(t,\eta) \widehat{\mathcal{U}}_{\mu_2}(t,\xi-\eta) d\eta \\
		= & \int \tau'_m(t) e^{it\Phi^{\mu_1,\mu_2}_\mu} \mathfrak{n}_+^{\mu_1,\mu_2}(\xi,\eta) \widehat{\mathcal{U}}_{\mu_1}(t,\eta) \widehat{\mathcal{U}}_{\mu_2}(t,\xi-\eta) d\eta \\
		+ & \int \tau_m(t) e^{it\Phi^{\mu_1,\mu_2}_\mu} \mathfrak{n}_+^{\mu_1,\mu_2}(\xi,\eta) \hspace{1mm} (\dt+\kap|\eta|^2) \widehat{\mathcal{U}}_{\mu_1}(t,\eta) \hspace{1mm} \widehat{\mathcal{U}}_{\mu_2}(t,\xi-\eta) d\eta \\
		+ & \int \tau_m(t) e^{it\Phi^{\mu_1,\mu_2}_\mu} \mathfrak{n}_+^{\mu_1,\mu_2}(\xi,\eta) \hspace{1mm} \widehat{\mathcal{U}}_{\mu_1}(t,\eta) \hspace{1mm} (\dt+\kap|\xi-\eta|^2)\widehat{\mathcal{U}}_{\mu_2}(t,\xi-\eta) d\eta \\
		+ & \int [i\Phi^{\mu_1,\mu_2}_\mu + \kap|\xi|^2 - \kap|\eta|^2 - \kap|\xi-\eta|^2] \tau_m(t) e^{it\Phi^{\mu_1,\mu_2}_\mu} \mathfrak{n}_+^{\mu_1,\mu_2}(\xi,\eta) \widehat{\mathcal{U}}_{\mu_1}(t,\eta) \widehat{\mathcal{U}}_{\mu_2}(t,\xi-\eta) d\eta
	\end{align*}
	where the last line is simply the Fourier transform of $-\mathcal{Q}^{(m)}_{\mathfrak{m}_+^{\mu_1,\mu_2}} (\mathcal{U}_{\mu_1}, \mathcal{U}_{\mu_2})(t) = -(\dt-\kap\Delta) \mathcal{U}_+^m(t)$. Hence, merging this term with left-hand side, we get a heat equation for $\mathcal{U}_+^m(t) + \mathcal{Q}_{\mathfrak{n}_+^{\mu_1,\mu_2}}(\mathcal{U}_{\mu_1}, \mathcal{U}_{\mu_2})$, and solving the heat equation gives the desired result.
\end{proof}

The normal form will be helpful when we can use the discrepancy between the sizes of $\Lambda$ in $\xi, \xi-\eta,$ and $\eta$. Since the normal form produces ``time derivative terms", $(\dt-\kap\Delta)f$, we need a way to bound these using the quantities we can control. This is the content of following lemma, which is a counterpart of \cite[Lemma 6.1]{GPW}.

\begin{lemma}\label{Lemmadt-delta}
	Let $\U$ be a solution of \eqref{eq:NSC} satisfying the bootstrap assumptions \eqref{BootstrapAssumption} in Theorem \ref{MainThmBootstrap}. For $m\in\mathbb{N}$ and $t \in [2^m, 2^{m+1}]$, we have
	\begin{equation}\label{dt-deltadecay}
		\Vert P_k S^n (\dt - \kap\Delta) \U(t) \Vert_{L^2} \lesssim 2^{\frac{k}{2} - k^+} \cdot 2^{-\frac{3}{2}m + \gamma m} \cdot \epsilon^2,
	\end{equation}
	for some $0 < \gamma \ll \beta$, and any $0 \leq n \leq N$.
\end{lemma}

\begin{proof}
	The proof follows that of \cite[Lemma 6.1]{GPW}. The same method works because applying $\dt - \kap\Delta$ removes the main flow in the Duhamel's formula, and only leaves the dispersive bilinear term, i.e. the same terms as in the time derivative of Euler-Coriolis flow. In other words,
	\begin{align*}
		(\dt-\kap\Delta)\U = & (\dt-\kap\Delta) e^{t\kap\Delta}\U(0) + (\dt-\kap\Delta) \sum_{\mathfrak{m}} \int_0^t e^{(t-s)\kap\Delta} \mathcal{Q}_\mathfrak{m} (\U, \U)(s) ds \\
		= & \sum_{\mathfrak{m}} \mathcal{Q}_\mathfrak{m} (\U, \U)(t).
	\end{align*}
	Hence, \eqref{dt-deltadecay} follows in the same way, as we have same dispersive decay and proved similar energy estimates as in \cite{GPW}.
\end{proof}

\begin{corollary}\label{Cor-normalform}
	For the solution $\mathcal{U}$ of \eqref{eq:NSC}, define $\mathcal{U}^{m,a,b}_+ (t), \mathcal{U}^{m,a,b}_- (t)$ as
	\begin{equation*}
		\widehat{\mathcal{U}}^{m,a,b}_\mu(t) = \sum_{\mu_1,\mu_2 \in \{+,-\}} \int_0^t \tau_m(s) e^{(t-s)\kap\Delta} \int_{\mathds{R}^3} e^{is\Phi^{\mu_1,\mu_2}_\mu} \mathfrak{m}_\mu^{\mu_1,\mu_2}(\xi,\eta) \widehat{S^a \mathcal{U}}_{\mu_1}(s,\eta) \widehat{S^b \mathcal{U}}_{\mu_2}(s,\xi-\eta) d\eta ds.
	\end{equation*}
	Then, there holds
	\begin{equation}\label{normalform}
		\begin{split}
			\mathcal{U}^{m,a,b}_+(t) = \sum_{\mu_1,\mu_2 \in \{+,-\}} & \left\{ - \delta_{m,L} \mathcal{Q}_{\mathfrak{n}_+^{\mu_1,\mu_2,a,b}}(S^a \mathcal{U}_{\mu_1},S^b \mathcal{U}_{\mu_2})(t) + \delta_{m,0} e^{t\kap\Delta}[ \mathcal{U}^{m,a,b}_+(0) + \right. \mathcal{Q}_{\mathfrak{n}_+^{\mu_1,\mu_2,a,b}}(S^a \mathcal{U}_{\mu_1},S^b \mathcal{U}_{\mu_2})(0) ] \\
			& + \int_0^t e^{(t-s)\kap\Delta} \tau_m(s) \mathcal{Q}_{\mathfrak{n}_+^{\mu_1,\mu_2,a,b}}(S^a(\dt-\kap\Delta)\mathcal{U}_{\mu_1}, S^b \mathcal{U}_{\mu_2})(s)ds \\
			& + \int_0^t e^{(t-s)\kap\Delta} \tau_m(s) \mathcal{Q}_{\mathfrak{n}_+^{\mu_1,\mu_2,a,b}}(S^a \mathcal{U}_{\mu_1}, S^b(\dt-\kap\Delta)\mathcal{U}_{\mu_2})(s)ds \\
			& + \left. \int_0^t e^{(t-s)\kap\Delta} \tau_m'(s) \mathcal{Q}_{\mathfrak{n}_+^{\mu_1,\mu_2,a,b}}(S^a\mathcal{U}_{\mu_1}, S^b\mathcal{U}_{\mu_2})(s)ds \right\},
		\end{split}
	\end{equation}
	\begin{equation*}
		\mathfrak{n}_+^{\mu_1,\mu_2,a,b} = -\frac{\mathfrak{m}_+^{\mu_1,\mu_2}}{i\Phi_{\mu}^{\mu_1,\mu_2} + \kap|\xi|^2 - (2a+1)\kap|\eta|^2 - (2b+1)\kap|\xi-\eta|^2},
	\end{equation*}
	and similarly for $\mathcal{U}^{m,a,b}_-$.
\end{corollary}

The terms multiplied by $\kap$ in the denominator of $\mathfrak{n}$ will be ignored in estimates using the normal form, because we want to show estimates that are uniform with respect to $\kap$, in particular as $\kap \to 0$. In other words, qualitatively, this corollary will be used in the same way as Lemma \ref{normalform-base}, as we will only use the imaginary part of the denominator $i\Phi_{\mu}^{\mu_1,\mu_2}$ and ignore terms with $a$ and $b$, whenever we have changed the multipliers. It suffices to note that we can use the normal form even when the $S$-vector fields are applied to $\U$, even though \eqref{eq:NSC} does not commute with $S$ vector field. Together with Lemma \ref{Lemmadt-delta}, this gives another way to estimate bilinear terms.

\begin{proof}
	We modify the proof of Lemma \ref{normalform-base}. Again, by computation,
	{\allowdisplaybreaks
	\begin{align*}
		& (\dt+\kap|\xi|^2) \int \tau_m(t) e^{it\Phi_{\mu}^{\mu_1,\mu_2}} \mathfrak{n}_+^{\mu_1,\mu_2,a,b}(\xi,\eta) \widehat{S^a \mathcal{U}}_{\mu_1}(t,\eta) \widehat{S^b \mathcal{U}}_{\mu_2}(t,\xi-\eta) d\eta \\
		= \vphantom{+} & \int \tau'_m(t) e^{it\Phi_{\mu}^{\mu_1,\mu_2}} \mathfrak{n}_+^{\mu_1,\mu_2,a,b}(\xi,\eta) \widehat{S^a \mathcal{U}}_{\mu_1}(t,\eta) \widehat{S^b \mathcal{U}}_{\mu_2}(t,\xi-\eta) d\eta \\
		+ & \int \tau_m(t) e^{it\Phi_{\mu}^{\mu_1,\mu_2}} \mathfrak{n}_+^{\mu_1,\mu_2,a,b}(\xi,\eta) \hspace{1mm} (\dt+\kap|\eta|^2 + 2a\kap|\eta|^2) \widehat{S^a \mathcal{U}}_{\mu_1}(t,\eta) \hspace{1mm} \widehat{S^b \mathcal{U}}_{\mu_2}(t,\xi-\eta) d\eta \\
		+ & \int \tau_m(t) e^{it\Phi_{\mu}^{\mu_1,\mu_2}} \mathfrak{n}_+^{\mu_1,\mu_2,a,b}(\xi,\eta) \hspace{1mm} \widehat{S^a \mathcal{U}}_{\mu_1}(t,\eta) \hspace{1mm} (\dt+\kap|\xi-\eta|^2 + 2b\kap|\xi-\eta|^2)\widehat{S^b \mathcal{U}}_{\mu_2}(t,\xi-\eta) d\eta \\
		+ & \int [i\Phi_{\mu}^{\mu_1,\mu_2} + \kap|\xi|^2 - (2a+1)\kap|\eta|^2 - (2b+1)\kap|\xi-\eta|^2] \tau_m(t) e^{it\Phi_{\mu}^{\mu_1,\mu_2}} \mathfrak{n}_+^{\mu_1,\mu_2,a,b}(\xi,\eta) \widehat{S^a \mathcal{U}}_{\mu_1}(t,\eta) \widehat{S^b \mathcal{U}}_{\mu_2}(t,\xi-\eta) d\eta \\
		= \vphantom{+} & \int \tau'_m(t) e^{it\Phi_{\mu}^{\mu_1,\mu_2}} \mathfrak{n}_+^{\mu_1,\mu_2,a,b}(\xi,\eta) \widehat{S^a \mathcal{U}}_{\mu_1}(t,\eta) \widehat{S^b \mathcal{U}}_{\mu_2}(t,\xi-\eta) d\eta \\
		+ & \int \tau_m(t) e^{it\Phi_{\mu}^{\mu_1,\mu_2}} \mathfrak{n}_+^{\mu_1,\mu_2,a,b}(\xi,\eta) \hspace{1mm} \mathcal{F}\{S^a (\dt-\kap\Delta)\mathcal{U}_{\mu_1}\}(t,\eta) \hspace{1mm} \widehat{S^b \mathcal{U}}_{\mu_2}(t,\xi-\eta) d\eta \\
		+ & \int \tau_m(t) e^{it\Phi_{\mu}^{\mu_1,\mu_2}} \mathfrak{n}_+^{\mu_1,\mu_2,a,b}(\xi,\eta) \hspace{1mm} \widehat{S^a \mathcal{U}}_{\mu_1}(t,\eta) \hspace{1mm} \mathcal{F}\{S^b (\dt-\kap\Delta)\mathcal{U}_{\mu_2}\}(t,\xi-\eta) d\eta \\
		- & \int \tau_m(t) e^{it\Phi_{\mu}^{\mu_1,\mu_2}} \mathfrak{m}_+^{\mu_1,\mu_2}(\xi,\eta) \widehat{S^a \mathcal{U}}_{\mu_1}(t,\eta) \widehat{S^b \mathcal{U}}_{\mu_2}(t,\xi-\eta) d\eta
	\end{align*}}
	Observing that we recovered $\mathcal{U}^{m,a,b}_\lambda$ in the last line, the same argument as in the proof of Lemma \ref{normalform-base} gives us the desired result.
\end{proof}

\begin{remark}
	Since we now know how to use the normal form and Lemma \ref{Lemmadt-delta} together to obtain additional decay, for notational convenience, when we are already in the framework of considering functions $f_j = S^{n_j} \mathcal{U}_j$, we will write $(\dt - \kap\Delta) f_j$'s instead of reverting to $S^{n_j} (\dt - \kap\Delta) \mathcal{U}_j$, and declare we have used normal form.
\end{remark}

Also, whenever it's clear what we are dealing with, we will simply write $\mathfrak{n}$ for the altered multiplier instead of writing full notation $\mathfrak{n}_{\pm}^{\mu_1,\mu_2,a,b}$ to keep notations simpler. \\

There are occasions where we need to perform normal form using arbitrary cutoff $\lambda > 0$. In these cases, we use the notation
\begin{equation*}
	\mathfrak{m} = \psi(\lambda^{-1}\Phi)\mathfrak{m} + (1-\psi(\lambda^{-1}\Phi))\mathfrak{m} =: \mathfrak{m}^{res} + \mathfrak{m}^{nr}
\end{equation*}
to denote the resonant and nonresonant decomposition via multiplier. The bilinear terms can be expressed likewise, by writing
\begin{equation*}
	\mathcal{Q}_\mathfrak{m}(f_1,f_2) = \mathcal{Q}_{\mathfrak{m}^{res}}(f_1,f_2) + \mathcal{Q}_{\mathfrak{m}^{nr}}(f_1,f_2), \quad \mathcal{B}_\mathfrak{m}(g_1,g_2) = \mathcal{B}_{\mathfrak{m}^{res}}(g_1,g_2) + \mathcal{B}_{\mathfrak{m}^{nr}}(g_1,g_2),
\end{equation*}
and we have the following estimates.

\begin{lemma}\label{arbitrary-normal-form}
	Let $\lambda > 0$ and $F_j = P_{k_j, p_j} F_j$.
	(1) The non-resonant part satisfies
	\begin{equation}\label{normal-form-nr-1}
		||P_{k,p} \mathcal{Q}_{\Phi^{-1} \mathfrak{m}^{nr}}(F_1, F_2)||_2 \lesssim 2^{k+p_{\max}} \lambda^{-1} \cdot \min\{2^{k+p}, 2^{k_1+p_1}, 2^{k_2+p_2}\} ||F_1||_2 ||F_2||_2.
	\end{equation}
	(2) If we can choose $\lambda > 0$ such that $|\Phi \chi_h| \geq \lambda \gtrsim 1$, then we have that $\mathfrak{m}^{res} = 0$ and thus $\mathfrak{m} = \mathfrak{n}^{nr}$, and in addition to \eqref{normal-form-nr-1}, we also have
	\begin{equation}\label{normal-form-nr-2}
		||P_{k,p} \mathcal{Q}_{\Phi^{-1} \mathfrak{m}^{nr}}(F_1, F_2)||_2 \lesssim 2^{k+p_{\max}} \cdot \min\{||e^{it\Lambda}F_1||_\infty ||F_2||_2, \hspace{2mm} ||F_1||_2||e^{it\Lambda}F_2||_\infty \}.
	\end{equation}
	(3) If there holds that $|\partial_{\eta_3} \Phi \chi_h| \gtrsim L > 0$, then we can use the set size estimates for both resonant and nonresonant terms:
	\begin{gather}
		||P_{k,p} \mathcal{Q}_{\mathfrak{m}^{res}}(F_1, F_2)||_2 \lesssim 2^{k+p_{\max}} \lambda^{1/2} L^{-1/2} \cdot \min\{2^{k_1+p_1}, 2^{k_2+p_2}\} ||F_1||_2 ||F_2||_2, \\
		||P_{k,p} \mathcal{Q}_{\Phi^{-1} \mathfrak{m}^{nr}}(F_1, F_2)||_2 \lesssim 2^{k+p_{\max}} |\log \lambda| \lambda^{-1/2} L^{-1/2} \cdot \min\{2^{k_1+p_1}, 2^{k_2+p_2}\} ||F_1||_2 ||F_2||_2.
	\end{gather}
\end{lemma}
\begin{proof}
	See \cite[Lemma 5.10]{GPW}.
\end{proof}

Set size estimates can be used more generally, and will be used a lot in the analysis. See appendix A.2.

\subsection{Geometry of vectors}

The configuration of vectors $\xi, \xi-\eta, \eta$ will determine whether we can use integration by parts or normal form. As we can already see, the three vectors are related, and hence we can group the configurations into some cases amenable for our analysis.
\begin{lemma}\label{two_big_sizes}
	Assume that $p \leq \min\{p_1, p_2\} - 10$. Then on the support of $\chi_h$ there holds that $p+k < p_1 + k_1 - 4$, and thus $p_2 + k_2 - 2 \leq p_1 + k_1 \leq p_2 + k_2 + 2$. Moreover, either one of the following options holds:
	\begin{enumerate}
		\item $|k_1 - k_2| \leq 4$, and thus also $|p_1 - p_2| \leq 6$,
		\item $k_2 < k_1 - 4$ then $|k - k_1| \leq 2$ and $p_1 \leq p_2 - 2$, so that $p \leq p_1 - 10 \leq p_2 - 12$,
		\item  $k_1 < k_2 - 4$ then $|k - k_2| \leq 2$ and $p_2 \leq p_1 - 2$, so that $p \leq p_2 - 10 \leq p_1 - 12$.
	\end{enumerate}
\end{lemma}
\begin{proof}
	See \cite[Lemma 5.8]{GPW}.
\end{proof}

Similar results hold with $p,p_1,p_2$ permuted, and also with $q,q_1,q_2$ instead of $p,p_1,p_2$, as the proof only uses how the vectors are located. We add one more type of case where we have two small $p$'s and one large $p$ (or similarly with $q$'s).

\begin{lemma}\label{two_small_sizes}
	Assume that $p_1, p_2 \leq p - 15$ and $|p_1 - p_2| \leq 2$. Then on the support of $\chi_h$, there holds that $k_1, k_2 \geq k + 6$, $|k_1 - k_2| \leq 2$, and $|k_1+p_1 - k_2-p_2| \leq 4$, $k_2+p_2 \geq k+p - 2$.
	
	Similarly, assume that $q_1, q_2 \leq q - 15$ and $|q_1 - q_2| \leq 2$. Then on the support of $\chi$, there holds that $k_1, k_2 \geq k + 6$, $|k_1 - k_2| \leq 2$, and $|k_1+q_1 - k_2-q_2| \leq 4$, $k_2+q_2 \geq k+q - 2$.
\end{lemma}
\begin{proof}
	Suppose, for the sake of contradiction, that $k_1 < k+6$. From $\xi_h = (\xi-\eta)_h + \eta_h$, we have $2^{k+p-4} - 2^{k_1+p_1} \leq 2^{k_2+p_2}$. Since $(k+p) - (k_1+p_1) > 9$, this gives $2^{k+p-5} \leq 2^{k_2+p_2}$, hence $k - k_2 \leq p_2 - p + 5 \leq -10$. However, this leads to  $k \leq k_2 - 10$ and $k_1 \leq k_2 -4$ which contradicts $\xi = (\xi-\eta) + \eta$. Therefore, $k_1 \geq k+6$. Symmetric arguments give $k_2 \geq k+6$, and as the two largest indices, $k_1$ and $k_2$ should have similar sizes.
	
	$|k_1+p_1 - k_2-p_2| \leq 4$ is now immediate, and hence, it should be the case $(\xi-\eta)_h \sim \eta_h \gtrsim \xi_h$, which gives $k_2+p_2 \geq k+p - 2$.
	
	Replacing $(p,p_1,p_2)$ to $(q,q_1,q_2)$ and $(\xi_h, (\xi-\eta)_h, \eta_h)$ to $(\xi_3, (\xi-\eta)_3, \eta_3)$ and applying the same arguments gives the results between sizes of $q$'s and $k$'s.
\end{proof}

\begin{remark}
	Again, permuting the roles of $p, p_j$'s give analogous results. We will now use the notation $\ll$ instead of, for example, saying $p_1 \leq p - 15$, since the results follow from the configuration of vectors.
\end{remark}

\section{$B$-norm bounds}

In this section, we prove the first half of the bootstrap argument.
\begin{theorem}
	If $u$ is the solution to \eqref{eq:NSC} with $\kap \lesssim 1$ and the initial data satisfying \eqref{InitialCondition}, then assuming \eqref{BootstrapAssumption}, we have that
	\begin{equation*}
		||S^n \U(t)||_B \lesssim \epsilon + \epsilon^2
	\end{equation*}
\end{theorem}

\begin{proof}
	Control of the linear term $||S^n(e^{t\kap\Delta} \U(0))||_B = ||\overset{n}{\underset{a=0}{\sum}} h_a(t\kap\Delta) e^{t\kap\Delta} S^{n-a} \U(0)||_B \lesssim \epsilon$ was already mentioned in Section 3, and so were the steps to deal with the nonlinear term. We localize in time $t \in [2^m, 2^{m+1}]$ first, then in space by $k,p,l$, and in $q$, reducing the indices every time by energy estimates or bootstrap assumptions. After we are restricted to certain ranges of indices, we use repeated integration by parts using vector fields or perform normal forms, depending on the geometry of vectors in the restricted localized regions. Recall that
	\begin{equation*}
		S^n \mathcal{B}_\mathfrak{m}(\U, \U) = \sum_{a+b+c \leq n} \int_0^t \tau_m(s) h_a((t-s)\kap\Delta) e^{(t-s)\kap\Delta} \mathcal{Q}_\mathfrak{m}(S^b\U, S^c\U)(s) ds.
	\end{equation*}
	Since $||h_a(-(t-s)\kap|\xi|^2) e^{-(t-s)\kap|\xi|^2}||_\infty \lesssim 1$, it's enough to show
	\begin{equation*}
		|| \mathcal{B}_\mathfrak{m}(F_1, F_2) ||_B \lesssim 2^{-\delta m} \epsilon^2
	\end{equation*}
	for $F_j = S^{b_j} \U$. Also, we will often prove stronger bound
	\begin{equation*}
		2^{k+4k^+} || \mathcal{F}\{\mathcal{B}_\mathfrak{m} (F_1, F_2)\} ||_\infty \lesssim 2^{-\delta m} \epsilon^2.
	\end{equation*}
	
	First of all, using the energy bounds proved on previous section,
	\begin{align*}
		2^{k+4k^+} || \mathcal{F}\{P_k \mathcal{B}_\mathfrak{m} (P_{k_1} F_1, P_{k_2} F_2)\} ||_\infty
		\lesssim 2^{2k+4k^+} \sum_{k_1,k_2} 2^m & \cdot \min\{ 2^{-N_0 k_1} ||F_1||_{H^{N_0}}, 2^{k_1} ||F_1||_{\dot{H}^{-1}}\} \\
		& \cdot \min\{ 2^{-N_0 k_2} ||F_2||_{H^{N_0}}, 2^{k_2} ||F_2||_{\dot{H}^{-1}}\}.
	\end{align*}
	Since the two biggest indices of $\{k, k_1, k_2\}$ should have similar sizes, the summations to either infinity gives some trivial decay:
	\begin{align*}
		& \sum_{k_1 \geq \delta_0 m, k_2 \in \mathds{Z}} 2^{2k+4k^+} \cdot \min\{ 2^{-N_0 k_1} ||F_1||_{H^{N_0}}, 2^{k_1}||F_1||_{\dot{H}^{-1}}\}
		\cdot \min\{ 2^{-N_0 k_2} ||F_2||_{H^{N_0}}, 2^{k_2} ||F_2||_{\dot{H}^{-1}}\} \\
		\leq & \sum_{k_1 \geq \delta_0 m} 2^{-(N_0-6)k_1} \epsilon \cdot \sum_{k_2} \min\{ 2^{-N_0 k_2}, 2^{k_2}\} \epsilon \tag*{(when $k \simeq k_1 \gg k_2$)} \\
		+ & \sum_{k_1 \geq \delta_0 m} 2^{-N_0 k_1} \epsilon \cdot \sum_{k_2} \min\{ 2^{-(N_0-6) k_2}, 2^{3k_2}\} \epsilon \tag*{(when $k \simeq k_2 \gg k_1$)} \\
		+ & \sum_{k_1 \geq \delta_0 m} 2^{-(2N_0 - 6)k_1} \epsilon^2 \tag*{(when $k_1 \simeq k_2 \gg k$)} \\
		\simeq & \hspace{1mm} 2^{-(N_0-6)\delta_0 m} \epsilon^2 = 2^{(-2+\frac{12}{N_0})m} \epsilon^2
	\end{align*}
	for $\delta_0 = 2/N_0$ gives acceptable contribution as long as $N_0 > 12$. Also,
	\begin{align*}
		& \sum_{k_1 \leq -2m, k_2 \in \mathds{Z}} 2^{2k+4k^+} \cdot \min\{ 2^{-N_0 k_1} ||F_1||_{H^{N_0}}, 2^{k_1}||F_1||_{\dot{H}^{-1}}\}
		\cdot \min\{ 2^{-N_0 k_2} ||F_2||_{H^{N_0}}, 2^{k_2} ||F_2||_{\dot{H}^{-1}}\} \\
		\leq & \sum_{k_1 \leq -2m} 2^{3k_1} \epsilon \cdot \sum_{k_2} \min\{ 2^{-N_0 k_2}, 2^{k_2}\} \epsilon \tag*{(when $k \simeq k_1 \gg k_2$)} \\
		+ & \sum_{k_1 \leq -2m} 2^{k_1} \epsilon \cdot \sum_{k_2} \min\{ 2^{-(N_0-6) k_2}, 2^{3k_2}\} \epsilon \tag*{(when $k \simeq k_2 \gg k_1$)} \\
		+ & \sum_{k_1 \leq -2m} 2^{2k_1} \epsilon^2 \tag*{(when $k_1 \simeq k_2 \gg k$)} \\
		\simeq & \hspace{1mm} 2^{-2m} \epsilon^2,
	\end{align*}
	which is again an acceptable summation. Similar calculations apply for summations over $k_1 \in \mathds{Z}, k_2 \notin (-2m, \delta_0 m)$, and from $2^k \lesssim 2^{k_1} + 2^{k_2}$, this also restricts $k$ to be less than $\delta_0 m$. Lastly, when $k \leq -2m$,
	\begin{align*}
		& \sum_{k_1, k_2 \in \mathds{Z}} 2^{2k+4k^+} \cdot \min\{ 2^{-N_0 k_1} ||F_1||_{H^{N_0}}, 2^{k_1}||F_1||_{\dot{H}^{-1}}\}
		\cdot \min\{ 2^{-N_0 k_2} ||F_2||_{H^{N_0}}, 2^{k_2} ||F_2||_{\dot{H}^{-1}}\} \\
		\leq & \hspace{1mm} 2^{-2m} \sum_{k_1 \in \mathds{Z}} \min\{ 2^{-N_0 k_1}, 2^{k_1}\} \epsilon \cdot \sum_{k_2 \in \mathds{Z}} \min\{ 2^{-N_0 k_2}, 2^{k_2}\} \epsilon \simeq 2^{-2m} \epsilon^2,
	\end{align*}
	so we can restrict to the case where $-2m \leq k, k_1, k_2 \leq \delta_0 m$. Now localizing further in $p_j, l_j, j=1,2$, with denoting $f_j = P_{k_j,p_j} R_{l_j} F_j$,
	\begin{equation*}
		2^{k+4k^+} ||\mathcal{F}\{ P_k \mathcal{B}_\mathfrak{m}(f_1,f_2) \} ||_\infty
		\lesssim 2^{k+4k^+} \cdot 2^m \cdot 2^{k+p_{\max}} \cdot 2^{-3k_1^+ - l_1} ||f_1||_X \cdot 2^{-3k_2^+ - l_2} ||f_2||_X \lesssim 2^{(1+3\delta_0)m -l_1-l_2} \epsilon^2
	\end{equation*}
	so that summation over either $l_1 \geq 2m$ or $l_2 \geq 2m$ give acceptable bounds. Along with the natural assumption made before that $p_j + l_j \geq 0$, we can only consider the indices
	\begin{equation*}
		-2m \leq k, k_1, k_2 \leq \delta_0 m, \hspace{4mm} -2m \leq p_1, p_2 \leq 0, \hspace{4mm} -p_1 \leq l_1 \leq 2m, -p_2 \leq l_2 \leq 2m.
	\end{equation*}
	Since these are only $(\log \langle t\rangle)^6$ many cases ($k$ is fixed), now it's enough to show
	\begin{equation}\label{Bnorm_goal}
		2^{3k^+ - \frac{1}{2}k^- - p -\frac{q}{2}} ||P_{k,p,q} \mathcal{B}_\mathfrak{m}(f_1,f_2) ||_2 \lesssim 2^{-\delta m} \epsilon^2
	\end{equation}
	for fixed $k,k_j,p_j,q_j, j=1,2$ in the above ranges. Also, in many cases, we will show instead
	\begin{equation*}
		2^{k+4k^+} ||\mathcal{F}\{P_{k,p,q} \mathcal{Q}_\mathfrak{m}(f_1,f_2) \}||_\infty \lesssim 2^{-(1+\delta)m} \epsilon^2.
	\end{equation*}
	
	\noindent\textbf{1. Gap in $p$ with $p_{\max} \sim 0$.}\\ From $|\bar{\sigma}| \sim 2^{k_{\min} + k_{\max} + p_{\max}}$, we can conduct integration by parts in a controllable way in this region. We use Lemma \ref{IbyP} to conduct repeated integration by parts along either vector fields, and assume without loss of generality that $p_1 \leq p_2$. \\
	\textcircled{1} $p \ll p_1,p_2$ : By Lemma \ref{two_big_sizes}, we have three possible cases. \\
	
	\noindent 1-1 : $2^{k_1} \sim 2^{k_2}, \hspace{2mm} 2^{p_1} \sim 2^{p_2} \sim 1$. First, integration by parts gives,
	\begin{align*}
		2^{k+4k^+} ||\mathcal{F}\{ \mathcal{Q}_\mathfrak{m}(f_1,f_2) \}||_\infty
		& \lesssim 2^{k+4k^+} \cdot 2^{k+p_{\max}} \cdot (s^{-1} \cdot 2^{-p_1 + 2k_1 - k_{\max} - k_{\min}} [1+2^{k_2 - k_1 + l_1}] )^N ||(1,S^N) f_1||_2 ||(1,S^N) f_2||_2 \\
		& \lesssim 2^{2k+4k^+} \cdot (s^{-1} \cdot 2^{k_1 - k + l_1})^N \epsilon^2
	\end{align*}
	so that when $l_1 \leq k - k_1 + (1-\delta)m$, we are done. By symmetry, we are okay when $l_2 \leq k - k_2 + (1-\delta)m$. Else,
	\begin{align*}
		2^{k+4k^+} ||\mathcal{F}\{ \mathcal{Q}_\mathfrak{m}(f_1,f_2) \}||_\infty
		& \lesssim 2^{2k+4k^+} 2^{-3k_1^+ - l_1} ||f_1||_X 2^{-3k_2^+ - l_2} ||f_2||_X \\
		& \leq 2^{k_1 + k_2 - 2(1-\delta)m} \epsilon^2 \leq 2^{-(1+\delta)m} \epsilon^2.
	\end{align*}
	
	\noindent 1-2 : $2^{k_2} \ll 2^{k_1} \sim 2^k, \hspace{2mm} p \ll p_1 \ll p_2 \sim 0$. In this case, having $2^{-p_2}$ is favorable. Hence, integrating by parts in $V_\eta$,
	\begin{align*}
		2^{k+4k^+} ||\mathcal{F}\{ \mathcal{Q}_\mathfrak{m}(f_1,f_2) \}||_\infty
		& \lesssim 2^{k+4k^+} \cdot 2^{k} \cdot (s^{-1} \cdot 2^{-p_2 + 2k_2 - k_{\max} - k_{\min}} [1+2^{k_1 - k_2 + l_2}] )^N ||(1,S^N) f_1||_2 ||(1,S^N) f_2||_2 \\
		& \lesssim 2^{2k+4k^+} \cdot (s^{-1} \cdot 2^{l_2})^N \epsilon^2
	\end{align*}
	so that if $l_2 \leq (1-\delta)m$ we are done. Else,
	\begin{align*}
		2^{k+4k^+} ||\mathcal{F}\{ \mathcal{Q}_\mathfrak{m}(f_1,f_2) \}||_\infty
		& \lesssim 2^{2k+4k^+} ||f_1||_2 \cdot 2^{-3k_2^+ - (1+\beta)l_2 - \beta p_2} ||f_2||_X \\
		& \lesssim 2^{2k+4k^+} 2^{-(1+\beta)(1-\delta)m} \epsilon^2 \leq 2^{-(1+\delta)m} \epsilon^2.
	\end{align*}
	
	\noindent 1-3 : $2^{k_1} \ll 2^{k_2} \sim 2^k, \hspace{2mm} p \ll p_2 \ll p_1 \sim 0$. This is against the assumption $p_1 \leq p_2$, hence is excluded. \\
	
	\noindent\textcircled{2} $p_1 \ll p_2, p$ : By Lemma \ref{two_big_sizes}, we again have three possible cases.\\
	2-1 : $2^{k} \sim 2^{k_2}, \hspace{2mm} 2^{p} \sim 2^{p_2} \sim 1$. Via the same integration by parts as in the case 1-1,
	\begin{align*}
		2^{k+4k^+} ||\mathcal{F}\{ \mathcal{Q}_\mathfrak{m}(f_1,f_2) \}||_\infty
		& \lesssim 2^{2k+4k^+} \cdot (s^{-1} 2^{-p_1 + l_1})^N \epsilon^2 \\
		\text{or, } & \lesssim 2^{2k+4k^+} \cdot (s^{-1} [2^{k_2 - k_1} + 2^{l_2}])^N \epsilon^2
	\end{align*}
	so that if $l_1 \leq p_1 + (1-\delta)m$ or $\max\{k_2 - k_1, l_2\} \leq (1-\delta)m$, we are done. Else, we use the set size estimate Lemma \ref{set-size-estimate}:
	\begin{align*}
		2^{3k^+ - \frac{1}{2}k^- -p -\frac{q}{2}} ||P_{k,p,q} \mathcal{Q}_\mathfrak{m}(f_1,f_2) ||_2
		& \lesssim 2^{3k^+ - \frac{1}{2}k^- -\frac{q}{2}} 2^k 2^{k_1+p_1 + \frac{k+q}{2}} 2^{-l_1}||f_1||_X \cdot 2^{-3k_2^+ - l_2} ||f_2||_X \\
		& \lesssim 2^{\frac{k^+}{2} + 2k + (k_1-k_2) - l_2 + (p_1 - l_1)} \epsilon^2 \leq 2^{\frac{k^+}{2} + 2k - 2(1-\delta)m} \epsilon^2 \leq 2^{-\delta m} \epsilon^2.
	\end{align*}
	
	\noindent 2-2 : $2^{k_2} \ll 2^k \sim 2^{k_1}, \hspace{2mm} p_1 \ll p \ll p_2 \sim 0$. The relations between $k$'s and the size of $p_2$ are the same as those in the case 1-2, and thus can be dealt with exactly the same way as those were enough to close the estimate.
	
	\noindent 2-3 : $2^k \ll 2^{k_2} \sim 2^{k_1}, \hspace{2mm} p_1 \ll p_2 \ll p \sim 0$. Integrating by parts along $V_\eta$,
	\begin{align*}
		2^{k+4k^+} ||\mathcal{F}\{ \mathcal{Q}_\mathfrak{m}(f_1,f_2) \}||_\infty
		& \lesssim 2^{2k+4k^+} \cdot (s^{-1} \cdot 2^{-p_2 + 2k_2 - k_{\max} - k_{\min}} [1+2^{k_1 - k_2 + l_2}] )^N ||(1,S^N) f_1||_2 ||(1,S^N) f_2||_2 \\
		& \lesssim 2^{2k+ 4k^+} (s^{-1} 2^{-p_2 + k_2 - k + l_2})^N \epsilon^2
	\end{align*}
	so that when $ l_2 - p_2 + k_2 - k \leq (1-\delta)m$, we are done. Else, we also use $2^{k+p} \sim 2^{k_2+p_2}$, so that $2^{\beta(k_2 - k)} \sim 2^{-\beta p_2}$, and
	\begin{align*}
		2^{3k^+ - \frac{1}{2}k^- -p -\frac{q}{2}} ||P_{k,p,q} \mathcal{Q}_\mathfrak{m}(f_1,f_2) ||_2
		& \lesssim 2^{3k^+ - \frac{1}{2}k^- -\frac{q}{2}} 2^k 2^{k_2 + p_2 + \frac{k+q}{2}} \cdot 2^{p_1} ||f_1||_2 \cdot 2^{-3k_2^+ - (1+\beta)l_2 - \beta p_2} ||f_2||_X \\
		& \lesssim 2^{\frac{1}{2}k^+ + k + k_2} 2^{p_1 + p_2} 2^{-\beta p_2} 2^{-(1+\beta)p_2 + (1+\beta)(k_2 - k) - (1+\beta)(1- \delta)m} \epsilon^2 \\
		& \lesssim 2^{\frac{1}{2}k^+ + 2k_2} 2^{p_1 - 3\beta p_2} 2^{-(1 + \beta/2)m} \epsilon^2 \leq 2^{-(1+\delta)m} \epsilon^2
	\end{align*}
	from $p_1 \leq p_2$ and $\beta \leq 1/3$. \\
	
	\noindent \textcircled{3} $p_2 \ll p_1, p$ : This case is excluded from our assumption $p_1 \leq p_2$. \\
	
	\noindent \textcircled{4} $p_1 \sim p_2 \ll p \sim 0$ : By Lemma \ref{two_small_sizes}, $k_1 \sim k_2 \gg k$ holds, and using integration by parts first,
	\begin{align*}
		2^{k+4k^+} ||\mathcal{F}\{ \mathcal{Q}_\mathfrak{m}(f_1,f_2) \}||_\infty
		& \lesssim 2^{2k+4k^+} \cdot (s^{-1} 2^{-p_1 + 2k_1 - k - k_1} \cdot [1+ 2^{k_2 - k_1 + l_1}])^N ||(1,S)^N f_1||_2 ||(1,S)^N f_2||_2 \\
		& \lesssim 2^{2k+4k^+} \cdot (s^{-1} 2^{-p_1 + k_1 - k + l_1})^N \epsilon^2,
	\end{align*}
	so that if $l_1 - p_1 + k_1 - k \leq (1-\delta)m$, we are done. Else, we also use $2^{k_1 + p_1} \sim 2^{k_2 + p_2} \gtrsim 2^{k+p} \sim 2^k \leftrightarrow -p_1 - k_1 +k \lesssim 0$, so that
	\begin{align*}
		2^{k+4k^+} ||\mathcal{F}\{ \mathcal{Q}_\mathfrak{m}(f_1,f_2) \}||_\infty
		& \lesssim 2^{2k+4k^+} \cdot 2^{-3k_1^+ -(1+\beta)l_1 - \beta p_1} ||f_1||_X \cdot 2^{-3k_2^+ + p_2} ||f_2||_B \\
		& \lesssim 2^{2k} 2^{(1+\beta)(-p_1 + k_1 - k) - \beta p_1} 2^{-(1+\beta)(1-\delta)m} \cdot 2^{p_2} \epsilon^2 \\
		& \lesssim 2^{2k + (1+3\beta)(k_1 - k)} 2^{p_2 - p_1} 2^{-(1+\beta/2)m} \epsilon^2 \lesssim 2^{-(1+\delta)m} \epsilon^2
	\end{align*}
	since $p_1 \sim p_2$ and $\beta \leq 1/3$. \\
	
	\noindent \textcircled{5} $p \sim p_1 \ll p_2 \sim 0$ : Lemma \ref{two_small_sizes} gives $k \sim k_1 \gg k_2$. From integration by parts,
	\begin{align*}
		2^{k+4k^+} ||\mathcal{F}\{ \mathcal{Q}_\mathfrak{m}(f_1,f_2) \}||_\infty
		& \lesssim 2^{2k+4k^+} \cdot (s^{-1} 2^{-p_2 + 2k_2 - k_1-k_2}[1+ 2^{k_1 - k_2 + l_2}])^N ||(1,S)^N f_1||_2 ||(1,S)^N f_2||_2 \\
		& \lesssim 2^{2k+4k^+} \cdot (s^{-1} 2^{l_2})^N \epsilon^2
	\end{align*}
	so that when $l_2 \leq (1-\delta)m$, we are done. Else,
	\begin{align*}
		2^{k+4k^+} ||\mathcal{F}\{ \mathcal{Q}_\mathfrak{m}(f_1,f_2) \}||_\infty
		& \lesssim 2^{2k+4k^+} ||f_1||_2 \cdot 2^{-3k_2^+ - (1+\beta)l_2 - \beta p_2} ||f_2||_X \\
		& \lesssim2^{2k + 4k^+} 2^{-(1+\beta)(1-\delta)m} \epsilon^2 \lesssim 2^{(-(1+\delta)m)} \epsilon^2.
	\end{align*}
	
	\noindent \textcircled{6} $p \sim p_2 \ll p_1 \sim 0$ : This case is excluded by our assumption $p_1 \leq p_2$. \\
	
	\noindent\textbf{2. Gap in $p$ with $p_{\max} \ll 0$.} In this case, every $\Lambda(\zeta) \geq \frac{1}{2}, \zeta \in \{\xi,\xi-\eta,\eta\}$, so that we have $|\Phi| \geq \frac{1}{2}$, enabling us to perform normal form with $\mathfrak{n} \lesssim 2^k$ every time. Compared to \eqref{normalform}, we have $h_a((t-s)\kappa\Delta)$'s in addition, but one can immediately observe that this can be absorbed by the exponential again. The first two terms with time integration are bounded easily from
	\begin{align}\label{Bnorm-normalform}
		2^{k+4k^+} ||\mathcal{F}\{ \mathcal{B}_\mathfrak{n}((\dt-\kap\Delta)f_1,f_2) \}||_\infty & \lesssim 2^{k+4k^+} 2^m 2^{k+p_{\max}} ||(\dt-\kap\Delta)f_1||_2 ||f_2||_2 \\
		& \lesssim 2^{2k+4k^+} 2^m 2^{-\frac{3}{2}m + \gamma m} \epsilon^2 \cdot \epsilon \lesssim 2^{-\frac{1}{4}m} \epsilon^3
	\end{align}
	and the bound for the another term works symmetrically.
	
	We now turn to the first boundary term. This term is evaluated exactly at time $t$, so that we can use the linear decay estimates. When $\min\{p_1,p_2\} = p_1 \lesssim p$, set size estimate Lemma \ref{set-size-estimate} gives
	\begin{align*}
		2^{3k^+ - \frac{1}{2}k^-} 2^{-p-\frac{q}{2}} ||\mathcal{Q}_\mathfrak{n}(f_1,f_2)(t)||_2
		& \lesssim 2^{3k^+ - \frac{1}{2}k^- - p} 2^{k+p_{\max}} 2^{p_1} ||f_1||_B ||e^{it\Lambda} f_2||_\infty \\
		& \lesssim 2^{4k^+ + \frac{1}{2}k} 2^{-(1-\delta)m} \epsilon^2 \lesssim 2^{-\delta m} \epsilon^2,
	\end{align*}
	and similarly when $\min\{p_1,p_2\} = p_2 \lesssim p$. Hence we can now assume $p_{\min} = p$. If furthermore, $p_{\max} \leq -\delta m$, then
	\begin{align*}
		2^{k+4k^+} ||\mathcal{F}\{ \mathcal{Q}_\mathfrak{n}(f_1,f_2) \}||_\infty
		\lesssim 2^{2k+4k^+} 2^{p_{\max}} \cdot 2^{p_1}||f_1||_B \cdot 2^{p_2}||f_2||_B \lesssim 2^{-\delta m} \epsilon^2
	\end{align*}
	Lastly, if $p_{\max} \geq -\delta m$, then from $|\bar{\sigma}| \gtrsim 2^{-\delta m}2^{k_{\max} + k_{\min}}$ and
	\begin{align*}
		\left\vert \frac{1}{s \cdot V_\eta \Phi} V_\eta \left(\frac{1}{i\Phi + 2\kap(\xi-\eta)\cdot \eta}\right) \right\vert
		= \left\vert \frac{1}{s \cdot V_\eta \Phi} \frac{i V_\eta\Phi + 2\kap(\xi-\eta)\cdot\eta - 2\kap|\eta|^2}{[i\Phi + 2\kap(\xi-\eta)\cdot \eta]^2} \right\vert \lesssim s^{-1}|V_\eta \Phi|^{-1} (1 + \kap 2^{2k_2}),
	\end{align*}
	we can use repeated integration by parts again as in the Case 1.
	
	The second boundary term only exists when $m=0$, and hence the restrictions on the indices become fixed numbers. Also, $||e^{t\kap\Delta}\mathcal{Q}_\mathfrak{n}(\U,\U)(0)||_2 \lesssim \epsilon^2$, so it's acceptable.
	
	Lastly, from the above results,
	\begin{align*}
		2^{3k^+ - \frac{1}{2}k^-} 2^{-p-\frac{q}{2}} \left\Vert \int_0^t e^{(t-s)\kap\Delta} \tau_m'(s) \mathcal{Q}_\mathfrak{n}(f_1,f_2)(s)ds \right\Vert_2 & \leq \int_0^t |\tau_m'(s)| \cdot 2^{3k^+ - \frac{1}{2}k^-} 2^{-p-\frac{q}{2}} ||\mathcal{Q}_\mathfrak{n}(f_1,f_2)(s)||_2 ds \\
		& \lesssim 2^{-\delta m} \epsilon^2 \int_0^t |\tau_m'(s)| ds \lesssim 2^{-\delta m} \epsilon^2
	\end{align*}
	hence we are done. \\
	The proof above shows that the boundary term for $m=0 (t=0)$ and the term with $\tau_m'(s)$ will never be an issue when the other terms can be bounded. Hence, we omit them from the future use of normal forms.
	
	From hereon, we assume $p \sim p_1 \sim p_2$.\\
	
	\noindent \textbf{3. Gap in q.} We localize further in $q_1, q_2$, and reduce the cases first. For notational convenience, we keep using $f_j$'s for $P_{k_j,p_j,q_j}R_{l_j} F_j$'s. From
	\begin{equation*}
		2^{k+4k^+} ||\mathcal{F}\{ \mathcal{B}_\mathfrak{m}(f_1,f_2) \}||_\infty
		\lesssim 2^{2k+4k^+} 2^m \cdot 2^{\frac{q_1}{2}} ||f_1||_B \cdot 2^{\frac{q_2}{2}} ||f_2||_B,
	\end{equation*}
	we can see that it's enough to consider the cases when $ -3m \leq q_1,q_2 \leq 0$. Also, the assumption $q_{\min} \ll q_{\max} \leq 0$ makes $p_{\max} \sim 0$, so that all $2^p, 2^{p_1}, 2^{p_2}$ are comparable to 0 now. Without loss of generality, we assume $q_1 \leq q_2$.
	
	\noindent \textcircled{1} $q \ll q_1, q_2$ : By Lemma \ref{two_big_sizes}, we have three possible cases. \\
	3-1 : $2^{k_1} \sim 2^{k_2}, \hspace{2mm} 2^{q_1} \sim 2^{q_2}$. The integration by parts gives
	\begin{align*}
		2^{k+4k^+} ||\mathcal{F}\{ \mathcal{Q}_\mathfrak{m}(f_1,f_2) \}||_\infty
		& \lesssim 2^{2k+4k^+} (s^{-1} 2^{2k_1 - k_{\max} - k_{\min} - q_{\max}} \cdot [1+ 2^{k_2-k_1} (2^{q_2-q_1} + 2^{l_1}) ])^N ||(1,S)^N f_1||_2 ||(1,S)^N f_2||_2 \\
		& \lesssim 2^{2k+4k^+} \cdot (s^{-1} 2^{k_1 - k - q_{\max} + l_1})^N \epsilon^2,
	\end{align*}
	so that if $k_1 -k - q_{\max} + l_1 \leq (1-\delta)m$, we are done. Else,
	\begin{align*}
		2^{k+4k^+} ||\mathcal{F}\{ \mathcal{Q}_\mathfrak{m}(f_1,f_2) \}||_\infty
		& \lesssim 2^{k+4k^+} 2^{k+q_{\max}} \cdot 2^{-3k_1^+ + \frac{q_1}{2} - (1+\beta)l_1} ||f_1||_X ||f_2||_2 \\
		& \lesssim 2^{2k+k^+} 2^{\frac{3}{2}q_1} 2^{(1+\beta)(k_2-k) - (1+\beta)q_1 - (1+\beta)(1-\delta)m} \epsilon^2 \\
		& \lesssim 2^{-\delta m} \epsilon^2
	\end{align*}
	
	\noindent 3-2 : $2^{k_2} \ll 2^{k_1}, \hspace{2mm} 2^{q_1 - q_2} \sim 2^{k_2 - k_1} \ll 1$. Here, it's advantageous to see $-q_2$ from integration by parts, so
	\begin{align*}
		2^{k+4k^+} ||\mathcal{F}\{ \mathcal{Q}_\mathfrak{m}(f_1,f_2) \}||_\infty
		& \lesssim 2^{2k+4k^+} \cdot (s^{-1} 2^{2k_2 - k_1 - k_2 - q_2} \cdot [1 + 2^{k_1-k_2} (2^{q_1 - q_2} + 2^{l_2})] )^N \epsilon^2 \\
		& \lesssim 2^{2k+4k^+} \cdot (s^{-1} 2^{-q_2 + l_2})^N \epsilon^2
	\end{align*}
	so that if $-q_2 + l_2 \leq (1-\delta)m$, we are done. Else,
	\begin{align*}
		2^{k+4k^+} ||\mathcal{F}\{ \mathcal{Q}_\mathfrak{m}(f_1,f_2) \}||_\infty
		& \lesssim 2^{k+4k^+} 2^{k+q_{\max}} 2^{\frac{q_1}{2}} ||f_1||_B 2^{-3k_2^ - (1+\beta)l_2} ||f_2||_X \\
		& \lesssim 2^{2k+k^+} 2^{q_2 + \frac{1}{2}q_1} 2^{-(1+\beta)q_2 - (1+\beta)(1-\delta)m} \epsilon^2 \lesssim 2^{-\delta m} \epsilon^2
	\end{align*}
	
	\noindent 3-3 : $2^{k_1} \ll 2^{k_2}, \hspace{2mm} 2^{q_2 - q_1} \ll 1$. This is against the assumption $q_1 \leq q_2$, so we will not consider. \\
	
	\noindent \textcircled{2} $q_1 \ll q, q_2$ : Again by Lemma \ref{two_big_sizes}, we have three scenarios.
	
	\noindent 4-1 : $2^{k} \sim 2^{k_2}, \hspace{2mm} 2^q \sim 2^{q_2}$. Integration by parts gives
	\begin{align*}
		2^{k+4k^+} ||\mathcal{F}\{ \mathcal{Q}_\mathfrak{m}(f_1,f_2) \}||_\infty
		& \lesssim 2^{2k+4k^+} \cdot (s^{-1} 2^{2k_2 - k_{\max} - k_{\min} - q_{\max}} [1+2^{k_1-k_2} (2^{q_1 - q_2} + 2^{l_2})])^N \epsilon^2 \\
		& \lesssim 2^{2k+4k^+} \cdot (s^{-1} [2^{k_2 - k_1 - q_2} + 2^{-q_2 + l_2}])^N \epsilon^2
	\end{align*}
	so that if $k_2 - k_1 - q_2 \leq (1-\delta)m$ and $-q_2 + l_2 \leq (1-\delta)m$, we are done. Also,
	\begin{align*}
		2^{k+4k^+} ||\mathcal{F}\{ \mathcal{Q}_\mathfrak{m}(f_1,f_2) \}||_\infty
		& \lesssim 2^{2k+4k^+} \cdot (s^{-1} 2^{2k_1 - k_{\max} - k_{\min} - q_{\max}} [1+2^{k_2-k_1} (2^{q_2 - q_1} + 2^{l_1})])^N \epsilon^2 \\
		& \lesssim 2^{2k+4k^+} \cdot (s^{-1} [2^{-q_1} + 2^{-q_2 + l_1}])^N \epsilon^2
	\end{align*}
	so that we are also okay if $-q_1, -q_2 + l_1 \leq (1-\delta)m$. If $q_1 \leq -(1-\delta)m$,
	\begin{align*}
		2^{3k^+ - \frac{1}{2}k^- -\frac{q}{2}} ||\mathcal{Q}_\mathfrak{m}(f_1,f_2)||_2
		& \lesssim 2^{3k^+ - \frac{1}{2}k^- -\frac{q}{2}} 2^{k + \frac{k_1 + q_1}{2}} 2^{k+q_2} \cdot 2^{\frac{q_1}{2}} ||f_1||_B \cdot 2^{-(1+\beta)l_2} ||f_2||_X \\
		& \lesssim2^{\frac{3}{2}k + 4k^+} 2^{q_1 + \frac{q_2 + k_1}{2}} 2^{-(1+\beta)l_2} \epsilon^2
	\end{align*}
	since either $q_2 + k_1 \leq k_2 - (1-\delta)m$ or $q_2 - l_2 \leq -(1-\delta)m$, we have $-\frac{3}{2}(1-\delta)m$ combined with the bounde from $q_1$ and we are done. Lastly, if $q_2 - l_1 \leq -(1-\delta)m$,
	\begin{align*}
		2^{k+4k^+} ||\mathcal{F}\{ \mathcal{Q}_\mathfrak{m}(f_1,f_2) \}||_\infty
		\lesssim 2^{k+4k^+} 2^{k+q_2} 2^{-(1+\beta)l_1} ||f_1||_X 2^{\frac{q_2}{2}}||f_2||_B \lesssim 2^{-(1+\delta)m} \epsilon^2.
	\end{align*}
	
	\noindent 4-2 : $2^k \sim 2^{k_1}, \hspace{2mm} 2^{q-q_2} \sim 2^{k_2-k} \ll 1$. From
	\begin{align*}
		2^{k+4k^+} ||\mathcal{F}\{ \mathcal{Q}_\mathfrak{m}(f_1,f_2) \}||_\infty
		& \lesssim 2^{2k+4k^+} \cdot (s^{-1} 2^{2k_2 - k_{\max} - k_{\min} - q_{\max}} [1+2^{k_1-k_2} (2^{q_1 - q_2} + 2^{l_2})])^N \epsilon^2 \\
		& \lesssim 2^{2k+4k^+} \cdot (s^{-1} 2^{-q_2 + l_2})^N \epsilon^2
	\end{align*}
	so that if $-q_2 + l_2 \leq (1-\delta)m$, we are done. Else,
	\begin{align*}
		2^{k+4k^+} ||\mathcal{F}\{ \mathcal{Q}_\mathfrak{m}(f_1,f_2) \}||_\infty
		& \lesssim 2^{k+4k^+} \cdot 2^{k+q_2} \cdot 2^{q_1} ||f_1||_B \cdot 2^{-(1+\beta)l_2} ||f_2||_B \\
		& \lesssim 2^{2k + 4k^+} 2^{(1-\beta)q_1} 2^{(1+\beta)(q_2-l_2)} \epsilon^2 \lesssim 2^{-\delta m} \epsilon^2
	\end{align*}
	since $q_1 \leq q_2$.
	
	\noindent 4-3 : $2^{k_1} \sim 2^{k_2}, \hspace{2mm} 2^{q_2-q} \sim 2^{k-k_2} \ll 1$. Again, integration by parts gives,
	\begin{align*}
		2^{k+4k^+} ||\mathcal{F}\{ \mathcal{Q}_\mathfrak{m}(f_1,f_2) \}||_\infty
		& \lesssim 2^{2k+4k^+} \cdot (s^{-1} 2^{2k_2 - k_{\max} - k_{\min} - q_{\max}} [1+2^{k_1-k_2} (2^{q_1 - q_2} + 2^{l_2})])^N \epsilon^2 \\
		& \lesssim 2^{2k+4k^+} \cdot (s^{-1} 2^{k_2 - k -q + l_2})^N \epsilon^2
	\end{align*}
	so that if $k_2 - k -q + l_2 \leq (1-\delta)m$, we are done. Else,
	\begin{align*}
		2^{k+4k^+} ||\mathcal{F}\{ \mathcal{Q}_\mathfrak{m}(f_1,f_2) \}||_\infty
		& \lesssim 2^{k+4k^+} \cdot 2^{k+q} \cdot 2^{\frac{q_1}{2}} ||f_1||_B \cdot 2^{-(1+\beta)l_2} ||f_2||_X \\
		& \lesssim 2^{(1-\beta)k + 4k^+ + (\frac{1}{2} - \beta)q_1} 2^{(1+\beta)(k+q-l_2)} \epsilon^2 \lesssim 2^{-\delta m} \epsilon^2
	\end{align*}
	since $q_1 \leq q$. \\
	
	\noindent \textcircled{3} $q_2 \ll q, q_1$ : This case is against our assumption $q_1 \leq q_2$, hence excluded. \\
	
	\noindent\textcircled{4} $q_1 \sim q_2 \ll q$ : By Lemma \ref{two_small_sizes}, we get $k_1 \sim k_2 \gg k$. In this case, we have $|\Phi| \gtrsim 2^q$, so that we can use the normal form. One can see that the terms with $\dt-\kap\Delta$ can be treated in exactly the same way as in \eqref{Bnorm-normalform}. The boundary term is even easier to bound by
	\begin{align*}
		2^{3k^+ - \frac{1}{2}k^- -\frac{q}{2}} ||\mathcal{F}\{ \mathcal{Q}_\mathfrak{n}(f_1,f_2) \}||_2
		\lesssim 2^{3k^+ - \frac{1}{2}k^- -\frac{q}{2}} ||\mathcal{F} \mathfrak{n}||_1 ||e^{it\Lambda} f_1 ||_\infty 2^{\frac{q_2}{2}} ||f_2||_B \lesssim 2^{-\delta m} \epsilon^3
	\end{align*}
	
	\noindent\textcircled{5} $ q \sim q_1 \ll q_2$ : Lemma \ref{two_small_sizes} gives $k \sim k_1 \gg k_2$, and $|\Phi| \gtrsim 2^{q_2}$, so similarly as above,
	\begin{align*}
		2^{3k^+ - \frac{1}{2}k^- -\frac{q}{2}} ||\mathcal{F}\{ \mathcal{Q}_\mathfrak{n}(f_1,f_2) \}||_2
		\lesssim 2^{3k^+ - \frac{1}{2}k^- -\frac{q}{2}} ||\mathcal{F} \mathfrak{n}||_1 2^{\frac{q_1}{2}} ||f_1||_B ||e^{it\Lambda} f_2 ||_\infty \lesssim 2^{-\delta m} \epsilon^3
	\end{align*}
	and the other terms follow by the same reasoning.
	
	\noindent\textcircled{6} $q\sim q_2 \ll q_1$ : This case is excluded from our assumption $q_1 \leq q_2$.
	
	Hence, from hereon, we further assume $q \sim q_1 \sim q_2$.\\
	
	\noindent\textbf{4. No gaps.} Assume without loss of generality that $f_1$ has at most $N/2$ copies of vector field $S$. Using Proposition \ref{dispersive-decay}, we can decompose $P_{k_1,p_1,q_1} e^{it\Lambda}f_1$ into two terms such that
	\begin{equation*}
		||I||_\infty \lesssim 2^{-p_1-\frac{q_1}{2}}t^{-\frac{3}{2}} 2^{\frac{3}{2}k_1 - 3k_1^+} \epsilon, \hspace{6mm} ||II||_2 \lesssim 2^{3k_1^+} t^{-1/2} \epsilon,
	\end{equation*}
	for each $k_1,p_1,q_1$ that are leftover, and we have
	\begin{align*}
		2^{3k^+ - \frac{1}{2}k^- -\frac{q}{2}} ||\mathcal{F}\{ \mathcal{Q}_\mathfrak{m}(f_1,f_2) \}||_2 & \lesssim 2^{3k^+ - \frac{1}{2}k^- -\frac{q}{2}} 2^{k+p+q} [||I||_\infty 2^{p+\frac{q}{2}} ||f_2||_B + ||II||_2 ||e^{it\Lambda} f_2||_\infty] \\
		& \lesssim 2^{\frac{k}{2}+4k^+ +\frac{q}{2}} [2^{-\frac{3}{2}m}\epsilon^2 + 2^{3k_1^+} 2^{-\frac{1}{2}m} 2^{-\frac{2}{3}m} \epsilon^2] \lesssim 2^{-\frac{13}{12}m} \epsilon^2,
	\end{align*}
	hence we are finished.
	
\end{proof}

\section{$X$-norm bounds}

To bound the $X$-norm, first note that again from $S^n(e^{t\kap\Delta}\U(0)) = \underset{0 \leq a \leq n}{\sum} h_a(t\kap\Delta) e^{t\kap\Delta} S^{n-a}\U(0)$, the bound of the linear term is straightforward from the initial condition. Hence, the entire section is devoted to controlling the nonlinear term in the Duhamel's formula.

\subsection{Large $l$ case}

We first bound the $X$ norm in the case when the parameter $l$ is bigger compared to $m$. To be precise, we prove the following.

\begin{proposition}
	Assume the bootstrap condition \eqref{BootstrapAssumption} holds, and let $\delta = 2M^{-1/2} \ll \beta = \frac{1}{60}$. Then, for $F_j = S^{b_j}\U, 0 \leq b_1 + b_2 \leq N, j=1,2$, we have
	\begin{equation}\label{Xnorm-large-l}
		\sup_{k\in\mathds{Z}, l+p \geq 0, l \geq (1+\delta)m} 2^{3k_+}2^{(1+\beta)l}2^{\beta p} \Vert P_{k,p}R_l \mathcal{B}_\mathfrak{m}(F_1, F_2) \Vert_{L^2} \lesssim 2^{-\delta^2 m}\epsilon^2.
	\end{equation}
\end{proposition}

The rest of the section 7.1 will be the proof of this proposition. As in Section 6, $\mathcal{B}_\mathfrak{m}$ should be modified so that we have $h_n((t-s)\kappa\Delta)e^{(t-s)\kap\Delta}$'s instead of $e^{(t-s)\kap\Delta}$, but it will be immaterial as they mostly have same qualitative effect, and in particular have bounded operator norms. They will manifest their presence only when we are performing normal forms, where the effects were already captured in Corollary \ref{Cor-normalform}.\\

\noindent \textbf{7.1.1 When $l + p \leq \delta m$} \\
We first restrict the ranges of indices involved as in the previous section. By the energy estimates,
\begin{align*}
	& 2^{3k^+}2^{(1+\beta)l}2^{\beta p} \Vert P_{k,p}R_l S^n \mathcal{B}_\mathfrak{m}(f_1, f_2) \Vert_{L^2} \\
	\lesssim & \sum_{\substack{a+b+c \leq n \\ k_1,k_2}} 2^{3k^+}2^{(1+\beta)l} 2^{\beta p} \cdot 2^m \cdot ||h_a((t-s)\kap|\xi|^2) e^{-(t-s)\kap|\xi|^2} ||_{L^{\infty}_\xi} \cdot 2^{\frac{3}{2}k + p} \\
	& \hspace{11mm} \cdot \min\{2^{-N_0 k_1} ||S^b f_1||_{H^{N_0}}, 2^{k_1} ||S^b f_1||_{\dot{H}^{-1}} \} \cdot \min\{2^{-N_0 k_2} ||S^c f_2||_{H^{N_0}}, 2^{k_2} ||S^c f_2||_{\dot{H}^{-1}} \} \\
	\lesssim & \hspace{2mm} 2^{(1+\beta)\delta m + m} 2^{3k^+ + \frac{3}{2}k} \cdot \min\{2^{-N_0 k_1}, 2^{k_1} \} \cdot \min\{2^{-N_0 k_2} , 2^{k_2} \} \epsilon^2.
\end{align*}
Recalling that the biggest parameters in $k, k_1, k_2$ always come at least as a pair, we see that the summation gives \eqref{Xnorm-large-l} unless $-2m \leq k, k_1, k_2 \leq \delta_0 m$, where $\delta_0 = 2N_0^{-1}$ again. Similarly,
\begin{align*}
	2^{3k^+}2^{(1+\beta)l}2^{\beta p} \Vert P_{k,p}R_l S^n \mathcal{B}_\mathfrak{m}(f_1, f_2) \Vert_{L^2}
	& \lesssim \sum_{\substack{a+b+c \leq n \\ p_1,p_2 \leq 0 \leq p_1+l_1, p_2+l_2}} 2^{3k^+}2^{(1+\beta)l} 2^{\beta p} \cdot 2^m \cdot ||h_a((t-s)\kap|\xi|^2) e^{-(t-s)\kap|\xi|^2} ||_{L^{\infty}_\xi} \\
	& \cdot 2^{\frac{3}{2}k + p} \cdot \min\{2^{p_1} ||S^b f_1||_{B}, 2^{-l_1} ||S^b f_1||_{X} \} \cdot \min\{2^{p_2} ||S^c f_2||_{B}, 2^{-l_2} ||S^c f_2||_{X} \} \\
	\lesssim & \hspace{2mm} 2^{(1+\beta)\delta m + m + 5\delta_0 m} \cdot \min\{2^{p_1}, 2^{-l_1} \} \cdot \min\{2^{p_2}, 2^{-l_2} \} \epsilon^2,
\end{align*}
gives \eqref{Xnorm-large-l} if we limit the summation to the complement of $-2m \leq p_1, p_2 \leq 0$ and $-p_j \leq l_j \leq 2m, j=1,2$. Hence, from hereon, we assume
\begin{equation*}
	-2m \leq k, k_j \leq \delta_0 m, \hspace{4mm} -2m \leq p_j \leq 0, \hspace{4mm} -p_j \leq l_j \leq 2m, \hspace{3mm} j=1,2,
\end{equation*}
and hence ignore $2^{3k^+}$ in front. We divided into a few cases according to the indices.

\noindent
\textcircled{1} $p_1, p_2 \ll 0$: Since we already have $p \leq \delta m - l \leq -m$, this means $q \sim q_1 \sim q_2 \sim 0$, and hence we can use normal form. Terms with time derivative can be bounded as
\begin{align*}
	2^{(1+\beta)l}2^{\beta p} \Vert P_{k,p}R_l \mathcal{B}_\mathfrak{n}((\dt-\kap\Delta)f_1, f_2) \Vert_{L^2}
	& \lesssim 2^{(1+\beta)l} 2^{\beta p} \cdot 2^m \cdot 2^k \cdot 2^{\frac{3}{2}k + p} \cdot 2^{-\frac{3}{2}m + \gamma m} \epsilon^2 \cdot \epsilon \\
	& \lesssim 2^{(-\frac{1}{2} + \gamma + 2\delta + 3\delta_0)m} \epsilon^3,
\end{align*}
and symmetrically for the term with $(\dt - \kap\Delta)f_2$. The boundary term needs more care. WLOG assume $p_2 \leq p_1$. When $p_2 \lesssim p$,
\begin{align*}
	2^{(1+\beta)l}2^{\beta p} \Vert P_{k,p}R_l \mathcal{Q}_{\mathfrak{n}^{nr}}(f_1, f_2) \Vert_{L^2}
	& \lesssim 2^{(1+\beta)(l+p)} 2^{-p} \cdot 2^k \cdot ||e^{it\Lambda}f_1||_\infty \cdot 2^{p_2} ||f_2||_B \\
	& \lesssim 2^{\delta m + \delta_0 m} \cdot 2^{-m + \gamma m} \epsilon \cdot 2^{p_2 - p} \epsilon \\
	& \lesssim 2^{(-1 + \gamma + \delta + \delta_0)m} \epsilon^2.
\end{align*}
So now we assume $p \ll p_2 \leq p_1$. We integrate by parts using the fact that $|\bar{\sigma}| \sim 2^{p_1 + k_{\max} + k_{\min}}$ in this region. This will be enough if
\begin{equation*}
	2^{-2p_1+2k_1 - k_{\max} - k_{\min}}(1 + 2^{k_2 - k_1 + l_1}) \leq 2^{(1-\delta)m},
\end{equation*}
so we look at the complement cases. If $k = k_{\min}$, then $k_1 \sim k_2$, so this means we can assume $2p_1 + k - k_1 - l_1 \leq -(1-\delta)m$. Hence,
\begin{align*}
	2^{(1+\beta)l}2^{\beta p} \Vert P_{k,p}R_l \mathcal{Q}_{\mathfrak{n}^{nr}}(f_1, f_2) \Vert_{L^2}
	& \lesssim 2^{(1+\beta)l} 2^{\beta p} \cdot 2^k \cdot 2^{\frac{3}{2}k + p} \cdot 2^{-l_1} ||f_1||_X \cdot 2^{p_2} ||f_2||_B \\
	& \lesssim 2^{(1+\beta)\delta m} \cdot 2^{\frac{1}{2}(2p_1 + k - l_1)} 2^{2k - \frac{l_1}{2}} \epsilon^2 \\
	& \lesssim 2^{(-\frac{1}{2} + 2\delta + 2\delta_0)m} \epsilon^2.
\end{align*}
If $k_{\min} = k_1 \text{ or } k_2$, this means we can assume $2p_1 - l_1 \leq -(1-\delta)m$ or $2p_1 + k_2 - k_1 \leq -(1-\delta)m$ depending on the minimum $k_j$. By changing the set size estimate from $2^{\frac{3}{2}k + p}$ to $2^{k+p + \frac{k_{\min}}{2}}$, we can close the case in exactly the same way as above. \\

\noindent
\textcircled{2} $p_1 \sim 0$: The other case $p_2 \sim 0$ will follow by symmetric argument. Here, we can benefit from integration by parts as now we have $|\bar{\sigma}| \sim 2^{k_{\max} + k_{\min}}$. Hence, integration by parts will be enough if
\begin{equation*}
	2^{-p_1+2k_1 - k_{\max} - k_{\min}}(1 + 2^{k_2 - k_1 + l_1}) \leq 2^{(1-\delta)m},
\end{equation*}
so we look at the complement cases again. If $k = k_{\min}$, then $k_1 \sim k_2$, so this means we can assume $k - k_1 - l_1 \leq -(1-\delta)m$. Hence,
\begin{align*}
	2^{(1+\beta)l}2^{\beta p} \Vert P_{k,p}R_l \mathcal{Q}_{\mathfrak{n}^{nr}}(f_1, f_2) \Vert_{L^2}
	& \lesssim 2^{(1+\beta)l} 2^{\beta p} \cdot 2^k \cdot 2^{\frac{3}{2}k + p} \cdot 2^{-l_1} ||f_1||_X \cdot ||f_2||_2 \\
	& \lesssim 2^{(1+\beta)\delta m} \cdot 2^{k - l_1} 2^{\frac{3}{2}k} \epsilon^2 \\
	& \lesssim 2^{(-1 + 2\delta + 2\delta_0)m} \epsilon^2.
\end{align*}
If $k_{\min} = k_1 \text{ or } k_2$, this means we can assume $- l_1 \leq -(1-\delta)m$ or (depending on the minimum $k_j$) $k_2 - k_1 \leq -(1-\delta)m$. The additional $k_2$ (if needed) can be obtained from using $\dot{H}^{-1}$ norm for $f_2$, and we can close the case in exactly the same way as above.

\noindent \textbf{7.1.2 When $l+p \geq \delta m$}

In this case, we have to deal with the largest parameter $l$ directly. Hence, the reduction of indices will be done in terms of $l$. From the energy estimates, we have
\begin{align*}
	& 2^{3k^+} 2^{(1+\beta)l}2^{\beta p} \Vert P_{k,p}R_l \mathcal{B}_\mathfrak{m}(f_1, f_2) \Vert_{L^2} \\
	\lesssim & \sum_{k_1,k_2} 2^{3k^+} 2^{(1+\beta)(l+p)} \cdot 2^m \cdot 2^{\frac{5}{2}k} \cdot \min\{2^{-N_0 k_1} ||f_1||_{H^{N_0}}, 2^{k_1} ||f_1||_{\dot{H}^{-1}} \} \cdot \min\{2^{-N_0 k_2} ||f_2||_{H^{N_0}}, 2^{k_2} ||f_2||_{\dot{H}^{-1}} \} \\
	\lesssim & \hspace{2mm} 2^{(1+\beta)(l+p) + m} 2^{3k^+ + \frac{3}{2}k} \cdot \min\{2^{-N_0 k_1}, 2^{k_1} \} \cdot \min\{2^{-N_0 k_2} , 2^{k_2} \} \epsilon^2.
\end{align*}
By the same reasoning as in the case \textbf{7.1.1}, this already gives the bound \eqref{Xnorm-large-l} unless $-2l \leq k,k_1,k_2 \leq \delta_0 l$. Similar repetitions with the use of $B$ and $X$ norms enable us to reduce the ranges to
\begin{equation*}
	-2l \leq k, k_j \leq \delta_0 l, \hspace{4mm} -2l \leq p_j \leq 0, \hspace{4mm} -p_j \leq l_j \leq 2l, \hspace{3mm} j=1,2.
\end{equation*}

From a brute estimate
\begin{align*}
	2^{3k^+} 2^{(1+\beta)l}2^{\beta p} \Vert P_{k,p}R_l \mathcal{B}_\mathfrak{m}(f_1, f_2) \Vert_{L^2}
	& \lesssim 2^{3k^+} 2^{(1+\beta)l} 2^{\beta p} \cdot 2^m \cdot 2^{k+p_{\max}} \cdot 2^{\frac{k_1+p_1}{3} + \frac{2(k_2+p_2)}{3} + \frac{1}{2}k_2} \cdot ||f_1||_2 ||f_2||_2 \\
	& \lesssim 2^m 2^{(1+\beta)(l - l_1 - l_2)} 2^{\beta p + (\frac{1}{3}-\beta)p_1 + (\frac{2}{3}-\beta)p_2 + p_{\max}} 2^{\frac{4}{3}k_{\max} + \frac{7}{6}k_2} \epsilon^2,
\end{align*}
if
\begin{equation}\label{large-l-cond1}
	l_1+l_2 \geq (1+\delta^2)l + (1-\frac{\beta}{2})m + k_2,
\end{equation}
we have the desired bound. Now we have to establish subtly different reasoning to extract the larger parameter $l$ instead of $m$ from integration by parts. Luckily, the estimates are not so different from the Euler-Coriolis case, since when we are using
\begin{equation*}
	R_l^{(j)} \mathcal{Q}_\mathfrak{m}(f_1,f_2) = 2^{-2l} R_l^{(j+1)} \Omega_{a3}^2 \mathcal{Q}_\mathfrak{m}(f_1,f_2).
\end{equation*}
$\Omega_{a3}$ acts to the bilinear term in the same way as in the Euler-Coriolis case thanks to $\Omega^\xi_{a3} |\xi| = 0$, so that $\Omega^\xi_{a3} (h_n(-(t-s)\kap|\xi|^2) e^{-(t-s)\kap|\xi|^2}) = 0$. Here, $R_l^{(j)}$s are all qualitatively similar localizations of $R_l = R_l^{(0)}$. Hence, we only need to estimate the terms where $\Omega^\xi_{a3}$ hits the phase, multiplier, or the input functions, which are the same terms as in the Euler-Coriolis case. To be more explicit, we will have terms of the form
\begin{equation}\label{omegaA3onDuhamel}
	\int_0^t (h_n(-(t-s)\kap|\xi|^2) e^{-(t-s)\kap|\xi|^2} \int_{\mathds{R}^3_\eta} (\Omega_{a3}^\xi)^{n_1}P_{k,p}(\xi) \cdot (\Omega_{a3}^\xi)^{n_2} \mathfrak{m} \cdot (\Omega_{a3}^\xi)^{n_3} e^{is\Phi} \cdot (\Omega_{a3}^\xi)^{n_4} P_{k_1, p_1}(\xi-\eta) \cdot (\Omega_{a3}^\xi)^{n_5} \hat{f_1}(\xi-\eta) \cdot \hat{f_2}(\eta) d\eta,
\end{equation}
where $n_1+n_2+n_3+n_4+n_5 \leq K$ for some $K$. To change derivatives in $\xi$ to those of $\xi-\eta$, we introduce
\begin{gather*}
	\Omega_{a3}^\xi = \Gamma_S S^{\xi-\eta} + \Gamma_\Omega \Omega_{a3}^{\xi-\eta}, \\
	\Gamma_S := \frac{\xi_a(\xi-\eta)_3 - \xi_3(\xi-\eta)_a}{|\xi-\eta|^2}, \hspace{3mm} \Gamma_\Omega := \frac{(\xi-\eta)_a \xi_a + (\xi-\eta)_3 \xi_3}{|\xi-\eta|^2}
\end{gather*}
which satisfy $|\Gamma_S|, |\Gamma_\Omega| \lesssim 2^{k - k_1}$. Direct computations to the `elementary components' give
\begin{align*}
	\Omega_{a3}^\xi \frac{(\xi-\eta)_a}{|\xi-\eta|} & = -\Gamma_\Omega \frac{(\xi-\eta)_3}{|\xi-\eta|}, \\
	\Omega_{a3}^\xi \frac{(\xi-\eta)_3}{|\xi-\eta|} & = \Gamma_\Omega \frac{(\xi-\eta)_a}{|\xi-\eta|} \\
	\Omega_{a3}^\xi \frac{\xi_a}{|\xi-\eta|} & = -\frac{\xi_3}{|\xi-\eta|} - \Gamma_S \frac{\xi_a}{|\xi-\eta|} \\
	\Omega_{a3}^\xi \frac{\xi_3}{|\xi-\eta|} & = \frac{\xi_a}{|\xi-\eta|} - \Gamma_S \frac{\xi_3}{|\xi-\eta|},
\end{align*}
and hence,
\begin{align*}
	\Omega_{a3}^\xi \Gamma_S & = -\Gamma_\Omega - \Gamma_S^2 + \Gamma_\Omega^2 \\
	\Omega_{a3}^\xi \Gamma_\Omega & = \Gamma_S (1 - 2\Gamma_\Omega),
\end{align*}
so that they form an algebra under the operation $\Omega_{a3}^\xi$, where each operation increases the order by 1. We also need how $\Lambda(\xi), \sqrt{1 - \Lambda^2(\xi)}, \frac{\xi_h \cdot \theta_h}{|\xi_h| |\theta_h|}, \frac{\xi_h^{\perp} \cdot \theta_h}{|\xi_h| |\theta_h|} , \theta \in \{\xi-\eta, \eta\}$ change under $\Omega_{a3}^\xi$. Again the action of $\Omega_{a3}$ on `elementary components' changes them into
\begin{align*}
	\Omega_{a3}^\xi \Lambda(\xi) & = \frac{\xi_a}{|\xi|}, \\
	\Omega_{a3}^\xi \frac{\xi_a}{|\xi|} & = - \Lambda(\xi), \\
	\Omega_{a3}^\xi \sqrt{1 - \Lambda^2(\xi)} & = - (1 - \Lambda^2(\xi))^{-\frac{1}{2}} \frac{\xi_a}{|\xi|}, \\
	\Omega_{a3}^\xi \frac{\xi_a}{|\xi_h|} & = -\frac{\xi_3}{|\xi_h|} + \frac{\xi_3 \xi_a^2}{|\xi_h|^3}, \\
	\Omega_{a3}^\xi \frac{\xi_b}{|\xi_h|} & = \frac{\xi_a \xi_b \xi_3}{|\xi_h|^3}, \\
	\Omega_{a3}^\xi \frac{\xi_3}{|\xi_h|} & = \frac{\xi_a}{|\xi_h|} + \frac{\xi_a \xi_3^2}{|\xi_h|^3},
\end{align*}
so that these form an algebra again. Using these to estimate the terms in \eqref{omegaA3onDuhamel}, we have
\begin{align*}
	|(\Omega_{a3}^\xi)^{n_1}P_{k,p}(\xi)| & \lesssim (2^{-p})^{n_1}, \\
	|(\Omega_{a3}^\xi)^{n_2} \mathfrak{m}| & \lesssim 2^{k + p_{\max}} (1 + 2^{-p} + 2^{k - k_1 - p_1})^{n_2} \\
	|(\Omega_{a3}^\xi)^{n_4} P_{k_1, p_1}(\xi-\eta)| & \lesssim (2^{k - k_1 - p_1})^{n_4} \\
	|(\Omega_{a3}^\xi)^{n_5} \hat{f_1}(\xi-\eta)| & \lesssim (2^{k - k_1})^{n_5} S^{j_1} \Omega_{a3}^{j_2} \hat{f_1}(\xi-\eta), \hspace{5mm} j_1 + j_2 \leq n_5.
\end{align*}
The phase term is more complicated, as it changes in different pattern depending on whether the $\Omega_{a3}^\xi$ acts on $e^{is\Phi}$ or on $(\Omega_{a3}^\xi)^{n} \Phi$. Instead of almost exact estimate, we use
\begin{equation*}
	\Omega_{a3}^\xi \Phi \lesssim 2^p + 2^{k-k_1+p_1}, \hspace{3mm} (\Omega_{a3}^\xi)^{n} \Phi \lesssim 1 + (2^{k-k_1})^n
\end{equation*}
to obtain
\begin{align*}
	(\Omega_{a3}^\xi)^{n_3} e^{is\Phi} & \lesssim \sum_{j=1}^{n_3} [s(2^p + 2^{k-k_1+p_1})]^j (1 + (2^{k-k_1})^{n_3 - j}) + \sum_{j=1}^{n_3/2} s^j (1 + (2^{k-k_1})^{n_3}) \\
	& \leq [s(2^p + 2^{k-k_1+p_1}) + 1 + 2^{k-k_1}]^{n_3} + s^{n_3/2} (1 + 2^{k-k_1})^{n_3}.
\end{align*}
Combining all these, modulo lower order terms, we finally achieve
\begin{align*}
	||R_l^{(n)} \mathcal{Q}_\mathfrak{m}(f_1, f_2) ||_2 & \lesssim 2^{\frac{3}{2}k_{\min}} 2^{k+p_{\max}} \sum_{j=0}^{1} 2^{-Kl} [2^{-p} + 2^m(2^p + 2^{k-k_1+p_1}) + 2^{m/2}(1 + 2^{k-k_1}) + 2^{k-k_1+l_1}]^K ||S^j f_1||_2 ||f_2||_2 \\
	+ & 2^{\frac{3}{2}k_{\min}} 2^{k+p_{\max}} 2^{-2l} (2^{-l}[2^{-p} + 2^m(2^p + 2^{k-k_1+p_1}) + 2^{m/2}(1 + 2^{k-k_1}) + 2^{k-k_1+l_1}])^{K'-2} ||S^2 f_1||_2 ||f_2||_2,
\end{align*}
by stopping integration by parts if $f_1$ gets two $S$ vector field applied.

This shows integration by parts along $\Omega_{a3}$ is enough when $2^{-p} + 2^m(2^p + 2^{k-k_1+p_1}) + 2^{m/2}(1 + 2^{k-k_1}) + 2^{k-k_1+l_1} < (1-\frac{\delta}{2})l$. In the complement range, we first remove the case when $-p$ is the dominant index. In this case,
\begin{align*}
	2^{3k^+} 2^{(1+\beta)l}2^{\beta p} \Vert P_{k,p}R_l \mathcal{B}_\mathfrak{m}(f_1, f_2) \Vert_{L^2}
	&\lesssim \sum_{\substack{k_1,p_1,l_1, \\ k_2,p_2,l_2}} 2^{3k^+ + \frac{5}{2}k} 2^{(1+\beta-K)(l+p)} \cdot 2^m \cdot \min\{2^{-N_0 k_2}, 2^{k_2}, 2^{p_2}, 2^{-l_1}\} \\
	& \hspace{13mm} \cdot \min\{ 2^{-N_0 k_2}, 2^{k_2}, 2^{p_2}, 2^{-l_2}\} \epsilon^2 \\
	& \lesssim 2^{-(1+\beta - K)\delta m + m} \epsilon^2,
\end{align*}
which can give acceptable contribution when $K$ is taken large enough. Hence, we drop $2^{-p}$ in the front from hereon. \\
After obtaining the analogous result for both $\xi-\eta$ and $\eta$, and hence assuming $k_1 \geq k_2$ without loss of generality, we are left with the cases when
\begin{equation*}
	l_1 \geq (1-\delta^2)l, \hspace{4mm} k - k_2 + \max\{m + p_2, \frac{m}{2}, l_2\} \geq (1-\delta^2)l.
\end{equation*}
If $l_2$ is the biggest, then \eqref{large-l-cond1} is satisfied, and we are done. If $m + p_2$ is the biggest, we have $k-k_2+p_2 \geq \delta m/2$, and integration by parts is enough if
\begin{equation*}
	2^{-p_2+2k_2 - p_{\max} - k_{\max} - k_{\min}}(1 + 2^{k_1 - k_2 + l_2}) \sim 2^{-p_2 - p_{\max} + k_2 - k_{\min} + l_2} \leq 2^{(1-3\delta)m},
\end{equation*}
since we now have a bound on $l$ in terms of $m$, $l \leq (1+2\delta^2)(m + p_2 + k - k_2)$, so that we can bound $2^{(1+\beta)l + m}$ through iterated integration by parts. Finally, we assume $p_2 + p_{\max} + k_{\min} - k_2 - l_2 \leq -(1-3\delta)m$. Harder case is when $k_{\min} = k$ so that $k_2 - k_{\min}$ does not vanish. In this case,
\begin{align*}
	& 2^{3k^+} 2^{(1+\beta)l}2^{\beta p} \Vert P_{k,p}R_l \mathcal{B}_\mathfrak{m}(f_1, f_2) \Vert_{L^2} \\
	& \lesssim 2^{3k^+} 2^{(1+\beta)l} 2^{\beta p} \cdot 2^m \cdot 2^{k+p_{\max}} \cdot 2^{k_1+p_1 + \frac{1}{2}k} \cdot 2^{-3k_1^+ -(1+\beta)l_1 - \beta p_1} ||f_1||_X \cdot 2^{-(1+\beta)l_2 - \beta p_2} ||f_2||_X \\
	& \lesssim 2^m \cdot 2^{k_1 + \frac{3}{2}k} 2^{(1+\beta)(l-l_1)} 2^{(1+\beta)(-p_2 - p_{\max} + k_2 - k - (1-3\delta)m)} 2^{\beta(p-p_1-p_2) + p_1 + p_{\max}} \epsilon^2 \\
	& \lesssim 2^m \cdot 2^{k_2 + \frac{3}{2}k} 2^{\delta^2 l} 2^{(1+\beta)(k_2 - k) - (1+\beta)m + 6\delta m} 2^{\beta(p - p_{\max}) - (1+2\beta)p_2} \epsilon^2 \\
	& \lesssim 2^{(-\beta + 6\delta)m} \cdot 2^{2\delta^2(m + k - k_2)} \cdot 2^{\frac{5}{2}k_2} 2^{(1+2\beta)(k - k_2 - \delta m/2)} \epsilon^2 \\
	& \lesssim 2^{(-\beta + 6\delta)m} \epsilon^2.
\end{align*}
Lastly, when $\frac{m}{2}$ is the biggest, we have $k - k_2 \geq \frac{1}{2}l$. Thus,
\begin{align*}
	2^{3k^+} 2^{(1+\beta)l}2^{\beta p} \Vert P_{k,p}R_l \mathcal{B}_\mathfrak{m}(f_1, f_2) \Vert_{L^2}
	& \lesssim 2^{3k^+} 2^{(1+\beta)l} 2^{\beta p} \cdot 2^m \cdot 2^{k+p_{\max}} \cdot 2^{\frac{3}{2}k_2 + p_2} \cdot 2^{-3k_1^+ - l_1} ||f_1||_X \cdot 2^{k_2} ||f_2||_{\dot{H}^{-1}} \\
	& \lesssim 2^m \cdot 2^{\beta l + (l-l_1)} 2^{k + \frac{5}{2}k_2} \epsilon^2 \\
	& \lesssim 2^m \cdot 2^{\beta l + \delta^2 l} \cdot 2^{\delta_0 l - \frac{5}{4}l + \frac{5}{2}k} \epsilon^2 \\
	& \lesssim 2^{(-\frac{1}{4} + \beta + \delta^2 + 4\delta_0)l} \epsilon^2.
\end{align*}
This completes the case with large $l$.

\subsection{Small $l$ case}

We investigate the case when $l \leq (1+\delta) m$. Now $m$ plays the role of large parameter again, and we can always use $(1+\beta)l \leq (1+\beta+2\delta)m$ when estimating $X$ norm.\\

\begin{proposition}
	Assume the bootstrap condition \eqref{BootstrapAssumption} holds, and let $\delta = 2M^{-1/2} \ll \beta \leq \frac{1}{60}$. Then for $F_j = S^{b_j} \mathcal{U}_\pm$, $0\leq b_1+b_2 \leq N, j=1,2$, we have
	\begin{equation}
		\sup_{k\in\mathds{Z}, l+p \geq 0, l \leq (1+\delta)m} 2^{3k_+}2^{(1+\beta)l}2^{\beta p} \Vert P_{k,p}R_l \mathcal{B}_\mathfrak{m}(F_1, F_2) \Vert_{L^2} \lesssim 2^{-\delta^2 m}\epsilon^2.
	\end{equation}
\end{proposition}
The rest of the section will be the proof of this proposition. The proof will follow quite similar pattern of the proof of $B$-norm bounds. The difference is that the estimates have to be more delicate as the $2^{(1+\beta)l}$ weight is harder to control. \\

Through the usual reduction process, we can assume
\begin{equation*}
	-2m \leq k, k_j \leq \delta_0 m, \hspace{4mm} -2m \leq p_j \leq 0, \hspace{4mm} -p_j \leq l_j \leq 2m, \hspace{3mm} j=1,2
\end{equation*}
Now we subdivide the cases according to the sizes of and relations between the indices as in Section 6.\\

\noindent\textbf{7.2.1 Gap in $p$ with $p_{\max} \sim 0$} \\
We consider the case when $p_{\min} \ll p_{\max} \sim 0$. WLOG, we may assume $p_1 \leq p_2$, which leave us with 4 different cases. \\
\textcircled{1} $p \ll p_1, p_2$: We invoke Lemma \ref{two_big_sizes} again.\\
\noindent
1-1) $k_1\sim k_2$ and $p_1 \sim p_2 \sim 0$ : Iterated integration by parts gives the result if $2^{k_1-k+l_1} \lesssim 2^{(1-\delta)m}$, or $2^{k_2-k+l_2} \lesssim 2^{(1-\delta)m}$. In other cases,
\begin{align*}
	2^{(1+\beta)l}2^{\beta p} \Vert P_{k,p}R_l \mathcal{B}_\mathfrak{m}(f_1, f_2) \Vert_{L^2}
	& \lesssim 2^{3k_+}2^{(1+\beta)l}2^{\beta p} \cdot 2^m \cdot 2^k \cdot 2^{\frac{3}{2}k+p} \cdot 2^{-(1+\beta)l_1} \Vert f_1\Vert_X \cdot 2^{-(1+\beta)l_2} \Vert_X \\
	& \lesssim 2^{(2+\beta+2\delta)m} \cdot 2^{\frac{5}{2}k} \cdot 2^{-2(1+\beta)(k_1 - k - (1-\delta)m)} \epsilon^2 \\
	& \lesssim 2^{-(\beta - 6\delta)m} \epsilon^2,
\end{align*}
which is acceptable since $\delta \ll \beta$.\\

\noindent
1-2) $k_2 \ll k_1 \sim k$, and $p \ll p_1 \ll p_2 \sim 0$ : Iterated integration by parts is enough if either $l_2 \leq (1-\delta)m$ or $-p_1 + \max\{k_1-k_2, l_1\} \leq (1-\delta)m$. In the opposite case, if $k_1-k_2 \geq l_1$,
\begin{align*}
	2^{(1+\beta)l}2^{\beta p} \Vert P_{k,p}R_l \mathcal{B}_\mathfrak{m}(f_1, f_2) \Vert_{L^2}
	& \lesssim 2^{(1+\beta)l}2^{\beta p} \cdot 2^m \cdot 2^k \cdot 2^{\frac{3}{2}k_2 + p_2} \cdot 2^{p_1} ||f_1||_B \cdot 2^{-(1+\beta)l_2} ||f_2||_X \\
	& \lesssim 2^{(2+\beta+2\delta)m} \cdot 2^{k+\frac{3}{2}k_2} \cdot 2^{(1+\beta)p_1} \cdot 2^{-(1+\beta)l_2} \epsilon^2 \\
	& \lesssim 2^{(2+\beta+2\delta)m} 2^{k + (1+\beta)k_1 + (\frac{1}{2}-\beta)k_2} 2^{-(1+\beta)(1-\delta)m} 2^{-(1+\beta)(1-\delta)m} \epsilon^2 \\
	& \lesssim 2^{-(\beta - 6\delta)m} \epsilon^2.
\end{align*}
If $k_1-k_2 \leq l_1$, using $2^{k_2} \sim 2^{k_1+p_1}$ in addition,
\begin{align*}
	2^{(1+\beta)l}2^{\beta p} \Vert P_{k,p}R_l \mathcal{B}_\mathfrak{m}(f_1, f_2) \Vert_{L^2}
	& \lesssim 2^{(1+\beta)l}2^{\beta p} \cdot 2^m \cdot 2^k \cdot 2^{\frac{1}{2}k_2 + k_1 + p_1} \cdot 2^{-(1+\beta)l_1 - \beta p_1} ||f_1||_X \cdot 2^{-(1+\beta)l_2} ||f_2||_X \\
	& \lesssim 2^{(2+\beta+2\delta)m} \cdot 2^{k+\frac{3}{2}k_1 + \frac{3}{2}p_1} \cdot 2^{-(1+\beta)l_1} 2^{-(1+\beta)l_2} \epsilon^2 \\
	& \lesssim 2^{(2+\beta+2\delta)m} \cdot 2^{3\delta_0 m} 2^{(\frac{3}{2} - \beta)p_1} 2^{-2(1+\beta)(1-\delta)m} \epsilon^2 \\
	& \lesssim 2^{-(\beta - 6\delta)m} \epsilon^2.
\end{align*}
\noindent 1-3) $p \ll p_2 \ll p_1$ : This case is excluded from the assumption $p_1 \leq p_2$. \\

\noindent
\textcircled{2} $p_1 \ll p_2, p$. Again we use Lemma \ref{two_big_sizes} to subdivide cases.\\
2-1) $k \sim k_2$, and $p \sim p_2 \sim 0$ : Integration by parts is enough if $-p_1+l_1 \leq (1-\delta)m$ or $\max\{k_2 - k_1, l_2\} \leq (1-\delta)m$. In other cases, when $k_2-k_1 \geq l_2$,
\begin{align*}
	2^{(1+\beta)l}2^{\beta p} \Vert P_{k,p}R_l \mathcal{B}_\mathfrak{m}(f_1, f_2) \Vert_{L^2}
	& \lesssim 2^{(1+\beta)l} \cdot 2^m \cdot 2^k \cdot 2^{\frac{3}{2}k_1 +p_1} \cdot 2^{-(1+\beta)l_1 - \beta p_1} ||f_1||_X \cdot \Vert f_2 \Vert_{L^2} \\
	& \lesssim 2^{(2+\beta+2\delta)m} 2^{\frac{5}{2}k_2 - 2\beta l_1} 2^{-(1-\beta)(1-\delta)m} 2^{-\frac{3}{2}(1-\delta)m} \epsilon^2 \\
	& \lesssim 2^{(-\frac{1}{2} - 2\beta -5\delta)m} \epsilon^2,
\end{align*}
which gives an acceptable bound. When $k_2 - k_1 \leq l_2$, we need to divide cases even further. If $2p_1 +k \gtrsim k_1$, we can get the extra $p_1$ from $k_1$:
\begin{align*}
	2^{(1+\beta)l}2^{\beta p} \Vert P_{k,p}R_l \mathcal{B}_\mathfrak{m}(f_1, f_2) \Vert_{L^2}
	& \lesssim 2^{(1+\beta)l} \cdot 2^m \cdot 2^k \cdot 2^{\frac{3}{2}k_1 + p_1} \cdot 2^{-(1+\beta)l_1 - \beta p_1} ||f_1||_X \cdot 2^{-(1+\beta)l_2} ||f_2||_X \\
	& \lesssim 2^{(2+\beta+2\delta)m} 2^{\frac{3}{2}k+k_1} \cdot 2^{-(1+\beta)l_1 + (2-\beta)p_1} 2^{-(1+\beta)l_2} \epsilon^2 \\
	& \lesssim 2^{(2+\beta+2\delta + 3\delta_0)m} 2^{-2(1+\beta)(1-\delta)m} \epsilon^2 \\
	& \lesssim 2^{-(\beta - 6\delta)m} \epsilon^2.
\end{align*}
Otherwise, we use normal form. From $\partial_{\eta_3} \Phi = \pm \frac{|\eta_h|^2}{|\eta|^3} \pm \frac{|\xi_h-\eta_h|^2}{|\xi-\eta|^3}$ and $\frac{|\eta_h|^2}{|\eta|^3} \sim 2^{2p_2-k_2} \sim 2^{-k_2}$, $\frac{|\xi_h-\eta_h|^2}{|\xi-\eta|^3} \sim 2^{2p_1-k_1}$, we have that $|\partial_{\eta_3} \Phi| \gtrsim 2^{-k_2}$. We choose $\lambda = 2^{-10\delta m}$ as the cutoff for resonant term, and use Lemma \ref{arbitrary-normal-form}:
\begin{align*}
	2^{(1+\beta)l}2^{\beta p} \Vert P_{k,p}R_l \mathcal{B}_{\mathfrak{m}^{res}}(f_1, f_2) \Vert_{L^2}
	& \lesssim 2^{(1+\beta)l} \cdot 2^m \cdot 2^k \cdot 2^{k_1+p_1} \cdot (2^{-10\delta m}2^{k_2})^{1/2} \cdot 2^{-l_1} ||f_1||_X \cdot 2^{-(1+\beta)l_2} ||f_2||_X \\
	& \lesssim 2^{(2+\beta+2\delta)m} 2^{\frac{3}{2}k + k_1} 2^{-5\delta m} 2^{p_1-l_1} 2^{-(1+\beta)l_2} \epsilon^2 \\
	& \leq 2^{[(4+\beta)\delta + 3\delta_0 - 5\delta]m} \epsilon^2
\end{align*}
which gives an acceptable contribution. For the nonresonant term, we can bound each term after using the normal form. We use another version of Lemma \ref{arbitrary-normal-form} for the boundary term:
\begin{align*}
	2^{(1+\beta)l}2^{\beta p} \Vert P_{k,p}R_l \mathcal{Q}_{\mathfrak{n}^{nr}}(f_1, f_2) \Vert_{L^2}
	& \lesssim 2^{(1+\beta)l} |\log(2^{-10\delta m})| \cdot 2^k \cdot 2^{k_1+p_1} \cdot (2^{-10\delta m}2^{-k_2})^{-1/2} \cdot 2^{-l_1} ||f_1||_X \cdot 2^{-(1+\beta)l_2} ||f_2||_X \\
	& \lesssim 2^{(1+\beta+2\delta)m} \cdot 10\delta m \cdot 2^{\frac{3}{2}k + k_1} 2^{5\delta m} 2^{p_1-l_1} 2^{-(1+\beta)l_2} \epsilon^2 \\
	& \lesssim 2^{-(1 - 20\delta)m} \epsilon^2.
\end{align*}
The other terms can be treated using Lemma \ref{Lemmadt-delta} :
\begin{align*}
	2^{(1+\beta)l}2^{\beta p} \Vert P_{k,p}R_l \mathcal{B}_{\mathfrak{n}^{nr}}((\dt-\kap\Delta)f_1, f_2) \Vert_{L^2}
	& \lesssim 2^{(2+\beta+2\delta)m} \cdot 2^{k+10\delta m} \cdot 2^{\frac{3}{2}k_1+p_1} \cdot 2^{-\frac{3}{2}m +\gamma m} \epsilon^2 \cdot 2^{-(1+\beta)l_2} ||f_2||_X \\
	& \lesssim 2^{-(\frac{1}{2} - \gamma - 16\delta)m} \epsilon^3,
\end{align*}
\begin{align*}
	2^{(1+\beta)l}2^{\beta p} \Vert P_{k,p}R_l \mathcal{B}_{\mathfrak{n}^{nr}}(f_1, (\dt-\kap\Delta)f_2) \Vert_{L^2}
	& \lesssim 2^{(2+\beta+2\delta)m} \cdot 2^{k+10\delta m} \cdot 2^{\frac{3}{2}k_1+p_1} \cdot 2^{-l_1} ||f_1||_X \cdot 2^{-\frac{3}{2}m +\gamma m} \epsilon^2 \\
	& \lesssim 2^{-(\frac{1}{2} - \beta - \gamma - 16\delta)m} \epsilon^3
\end{align*}
is enough since $\frac{1}{3} > \beta \gg \gamma$.\\

\noindent
2-2) $k_2 \ll k_1 \sim k$, and $p_1 \ll p \ll p_2 \sim 0$: Again, the integration by parts is enough when
\begin{gather*}
	2^{-p_1+2k_1 - k_{\max} - k_{\min}}(1 + 2^{k_2 - k_1 + l_1}) \sim 2^{-p_1}(2^{k_1 - k_2} + 2^{l_1}) \leq 2^{(1-\delta)m}, \quad \text{or,} \\
	2^{-p_2+2k_2 - k_{\max} - k_{\min}}(1 + 2^{k_1 - k_2 + l_2}) \sim 2^{-p_2 + l_2} \leq 2^{(1-\delta)m},
\end{gather*}
i.e. when $l_2 \leq (1-\delta)m$ or $-p_1+\max\{k_1-k_2, l_1\} \leq (1-\delta)m$. Hence, we look at the opposite case.

First we check the case $k_1-k_2 \geq l_1$. If $k_2 \leq -20\delta m$ in addition,
\begin{align*}
	2^{(1+\beta)l}2^{\beta p} \Vert P_{k,p}R_l \mathcal{B}_\mathfrak{m}(f_1, f_2) \Vert_{L^2}
	& \lesssim 2^{(2+\beta + 2\delta)m} 2^{\beta p} \cdot 2^k \cdot 2^{\frac{3}{2}k_2 + p_2} \cdot 2^{p_1} ||f_1||_B \cdot 2^{-(1+\beta)l_2} ||f_2||_X \\
	& \lesssim 2^{(2+\beta + 2\delta)m} \cdot 2^{k + \frac{k_2}{2}} \cdot 2^{k_2 + p_1} \cdot 2^{-(1+\beta)l_2} \epsilon^2 \\
	& \lesssim 2^{(2+\beta + 2\delta)m} \cdot 2^{2k_1 -10\delta m} \cdot 2^{-(1-\delta)m} \cdot 2^{-(1+\beta)(1-\delta)m} \epsilon^2 \\
	& \leq 2^{-5\delta m + 2\delta_0 m} \epsilon^2,
\end{align*}
which is acceptable since $\delta_0 \ll \delta$. Otherwise, if $k_2 \geq -20\delta m$,
\begin{align*}
	2^{(1+\beta)l}2^{\beta p} \Vert P_{k,p}R_l \mathcal{B}_\mathfrak{m}(f_1, f_2) \Vert_{L^2}
	& \lesssim 2^{(2+\beta + 2\delta)m} 2^{\beta p} \cdot 2^k \cdot 2^{\frac{3}{2}k_1 + p_1} \cdot 2^{p_1} ||f_1||_B \cdot 2^{-(1+\beta)l_2} ||f_2||_X \\
	& \lesssim 2^{(2+\beta + 2\delta)m} \cdot 2^{k + \frac{3}{2}k_1 + 2(k_2 + 20\delta m)} \cdot 2^{2p_1 - (1+\beta)l_2} \epsilon^2 \\
	& \lesssim 2^{(2+\beta + 42\delta)m} \cdot 2^{\frac{9}{2}k_1} \cdot 2^{-2(1-\delta)m - (1+\beta)(1-\delta)m} \epsilon^2 \\
	& \lesssim 2^{(-1 + 46\delta + 5\delta_0)m} \epsilon^2,
\end{align*}
which also gives enough contribution.

Now, we assume $l_1 \geq k_1 - k_2$, which gives us $p_1 - l_1 \leq -(1-\delta)m$. We use the set size estimate Lemma \ref{arbitrary-normal-form}(3) along with the normal form using the fact that $|\partial_3 \Phi| \sim |2^{2p_1-k_1} \pm 2^{2p_2-k_2}| \gtrsim 2^{2p_2 - k_2} \sim 2^{-k_2}$ is bounded below by $2^{-\delta_0 m}$. Dividing the terms into resonant and nonresonant terms according to $\lambda = 2^{-10\delta m}$,
\begin{align*}
	2^{(1+\beta)l}2^{\beta p} \Vert P_{k,p}R_l \mathcal{B}_{\mathfrak{m}^{res}}(f_1, f_2) \Vert_{L^2}
	& \lesssim 2^{(2+\beta + 2\delta)m} \cdot 2^k \cdot (\lambda 2^{\delta_0 m})^{\frac{1}{2}} \cdot 2^{k_1+p_1} \cdot 2^{-l_1} ||f_1||_X \cdot 2^{-(1+\beta)l_2} ||f_2||_X \\
	& \lesssim 2^{(2+\beta + 2\delta)m} \cdot 2^{k+k_1} \cdot 2^{(-5\delta + \frac{\delta_0}{2})m} \cdot 2^{-(1-\delta)m - (1+\beta)(1-\delta)m} \epsilon^2 \\
	& \lesssim 2^{(-(1-\beta)\delta + 3\delta_0)m} \epsilon^2,
\end{align*}
\begin{align*}
	2^{(1+\beta)l}2^{\beta p} \Vert P_{k,p}R_l \mathcal{Q}_{\mathfrak{n}^{nr}}(f_1, f_2) \Vert_{L^2}
	& \lesssim 2^{(1+\beta + 2\delta)m} \cdot 2^k \lambda^{-1} \cdot 2^{k_1+p_1} \cdot 2^{-l_1} ||f_1||_X \cdot 2^{-(1+\beta)l_2} ||f_2||_X \\
	& \lesssim 2^{(1+\beta + 2\delta)m} \cdot 2^{k+k_1+10\delta m} \cdot 2^{-(1-\delta)m - (1+\beta)(1-\delta)m} \epsilon^2 \\
	& \lesssim 2^{(-1 + 15\delta + 2\delta_0)m} \epsilon^2,
\end{align*}
\begin{align*}
	2^{(1+\beta)l}2^{\beta p} \Vert P_{k,p}R_l \mathcal{B}_{\mathfrak{n}^{nr}}((\dt-\kap\Delta)f_1, f_2) \Vert_{L^2}
	& \lesssim 2^{(2+\beta + 2\delta)m} \cdot 2^k \lambda^{-1} 2^{k_1+p_1} ||(\dt-\kap\Delta)f_1||_2 \cdot 2^{-(1+\beta)l_2} ||f_2||_X \\
	& \lesssim 2^{(2+\beta + 2\delta)m} \cdot 2^{k+k_1+10\delta m} \cdot 2^{-\frac{3}{2}m + \gamma m} \epsilon^2 \cdot 2^{-(1+\beta)(1-\delta)m} \epsilon \\
	& \lesssim 2^{(-\frac{1}{2} + \gamma + 15\delta + 2\delta_0)m} \epsilon^3,
\end{align*}
\begin{align*}
	2^{(1+\beta)l}2^{\beta p} \Vert P_{k,p}R_l \mathcal{B}_{\mathfrak{n}^{nr}}(f_1, (\dt-\kap\Delta)f_2) \Vert_{L^2}
	& \lesssim 2^{(2+\beta + 2\delta)m} \cdot 2^k \lambda^{-1} 2^{k_1+p_1} \cdot 2^{-l_1} ||f_1||_X || \cdot (\dt-\kap\Delta)f_2||_2 \\
	& \lesssim 2^{(2+\beta + 2\delta)m} \cdot 2^{k+k_1+10\delta m} \cdot 2^{-(1-\delta)m} \epsilon \cdot 2^{-\frac{3}{2}m + \gamma m} \epsilon^2 \\
	& \lesssim 2^{(-\frac{1}{2} + \beta + \gamma + 15\delta + 2\delta_0)m} \epsilon^3.
\end{align*}
Hence, we are done.\\

\noindent
2-3) $k \ll k_1 \sim k_2$, and $p_1 \ll p_2 \ll p \sim 0$: Again, integration by parts is enough if
\begin{align*}
	& 2^{-p_1+2k_1 - k_{\max} - k_{\min}}(1 + 2^{k_2 - k_1 + l_1}) \sim 2^{-p_1 + k_1 - k + l_1} \leq 2^{(1-\delta)m}, \quad \text{or,} \\
	& 2^{-p_2+2k_2 - k_{\max} - k_{\min}}(1 + 2^{k_1 - k_2 + l_2}) \sim 2^{-p_2 + k_2 - k + l_2} \leq 2^{(1-\delta)m}.
\end{align*}
So we assume $p_1 + k - l_1, p_2 + k - l_2 \leq k_1 - (1-\delta)m$. The issue is that it's hard to obtain all $p_1, p_2, -l_1, -l_2$ at the same time. In the case when $p_1 \geq -10\delta m$, this issue is resolved, and we have
\begin{align*}
	2^{(1+\beta)l}2^{\beta p} \Vert P_{k,p}R_l \mathcal{B}_\mathfrak{m}(f_1, f_2) \Vert_{L^2}
	& \lesssim 2^{(2+\beta + 2\delta)m} \cdot 2^k \cdot 2^{\frac{3}{2}k + p} \cdot 2^{(1+\beta)(p_1+p_2 + 20\delta m)} \cdot 2^{-(1+\beta)l_1} ||f_1||_X \cdot 2^{-(1+\beta)l_2} ||f_2||_X \\
	& \lesssim 2^{(2+\beta + 2\delta)m} \cdot 2^{(\frac{1}{2} - 2\beta)k + (1+\beta)20\delta m} \cdot 2^{(1+\beta)(p_1+k - l_1 + p_2 + k - l_2)} \epsilon^2 \\
	& \lesssim 2^{(-\beta + 30\delta + 3\delta_0)m} \epsilon^2,
\end{align*}
which is an enough contribution. Hence, we now assume $p_1 \leq -10\delta m$, and use normal form with $\lambda = 2^{-20}(2^q + 2^{2p_2})$. First, the nonresonant terms are bounded by
\begin{align*}
	2^{(1+\beta)l}2^{\beta p} \Vert P_{k,p}R_l \mathcal{Q}_{\mathfrak{n}^{nr}}(f_1, f_2) \Vert_{L^2}
	& \lesssim 2^{(1+\beta + 2\delta)m} \cdot 2^k (2^q + 2^{2p_2})^{-1} \cdot 2^{k_1+p_1 + \frac{k+q}{2}} \cdot 2^{-l_1} ||f_1||_X \cdot 2^{-(1+\beta)l_2 - \beta p_2} ||f_2||_X \\
	& \lesssim 2^{(1+\beta + 2\delta)m} \cdot \frac{2^{q/2}2^{p_2}}{2^q + 2^{2p_2}} \cdot 2^{\frac{3}{2}k + k_1} \cdot 2^{k_1-k-(1-\delta)m} 2^{-(1+\beta)(l_2+p_2)} \epsilon^2 \\
	& \leq 2^{(\beta + 3\delta)m} \cdot 2^{2k_1 + \frac{k}{2}} 2^{-(1+\beta)(l_2 + p_2)} \epsilon^2,
\end{align*}
which gives enough contribution if either $k - k_1 \leq -3\beta m$ or $l_2 + p_2 \geq 2\beta m$. However, if none of them hold
\begin{equation*}
	2\beta m \geq l_2 + p_2 \geq 2p_2 + k - k_1 + (1-\delta)m \gtrsim 3(k-k_1) + (1-\delta)m \geq (1-\delta - 9\beta)m
\end{equation*}
gives contradiction, where we used $2^{-p_2} \sim 2^{k_2-k} \sim 2^{k_1-k}$ in the current case. Hence, the boundary term is done, and
\begin{align*}
	2^{(1+\beta)l}2^{\beta p} \Vert P_{k,p}R_l \mathcal{B}_{\mathfrak{n}^{nr}}((\dt-\kap\Delta)f_1, f_2) \Vert_{L^2}
	& \lesssim 2^{(2+\beta + 2\delta)m} \cdot 2^k (2^q + 2^{2p_2})^{-1} \cdot 2^{k_2+p_2 + \frac{k+q}{2}} \cdot ||(\dt-\kap\Delta) f_1||_2 \cdot ||f_2||_2.
\end{align*}
From here, if $p_2 \leq -\frac{m}{3}$,
\begin{align*}
	& \lesssim 2^{(2+\beta + 2\delta)m} \cdot 2^{k+k_2} \cdot 2^{-\frac{3}{2}m + \gamma m} \epsilon \cdot ||f_2||_2 \\
	& \sim 2^{(\frac{1}{2} + \beta + \gamma + 2\delta)m} \cdot 2^{2k} 2^{2p_2} ||f_2||_B \epsilon^2 \\
	& \leq 2^{(-\frac{1}{6} + \beta + \gamma + 2\delta + 2\delta_0)m} \epsilon^3,
\end{align*}
closes the bound using $2^{p_2} \sim 2^{k-k_2}$, and if $p_2 \geq -\frac{m}{3}$,
\begin{align*}
	& \lesssim 2^{(2+\beta + 2\delta)m} \cdot 2^{k+k_2} \cdot 2^{-\frac{3}{2}m + \gamma m} \epsilon^2 \cdot 2^{-l_2} ||f_2||_X \cdot 2^{p_2 + \frac{m}{3}}\\
	& \lesssim 2^{(-\frac{1}{6} + \beta + \gamma + 3\delta + 2\delta_0)m} \epsilon^3.
\end{align*}
\begin{align*}
	2^{(1+\beta)l}2^{\beta p} \Vert P_{k,p}R_l \mathcal{Q}_{\mathfrak{n}^{nr}}(f_1, (\dt-\kap\Delta)f_2) \Vert_{L^2}
	& \lesssim 2^{(2+\beta + 2\delta)m} \cdot 2^k (2^q + 2^{2p_2})^{-1} 2^{k_1+p_1+ \frac{k+q}{2}} \cdot ||f_1||_2 \cdot 2^{-\frac{3}{2}m + \gamma m} \epsilon^2 \\
	& \lesssim 2^{(\frac{1}{2} + \beta + \gamma + 2\delta)m} \cdot 2^{2k_1 + \frac{k}{2}} \cdot 2^{2p_1} ||f_1||_B \epsilon^2
\end{align*}
finishes this term if $p_1 \leq -\frac{m}{3}$, and else,
\begin{align*}
	2^{(1+\beta)l}2^{\beta p} \Vert P_{k,p}R_l \mathcal{Q}_{\mathfrak{n}^{nr}}(f_1, (\dt-\kap\Delta)f_2) \Vert_{L^2}
	& \lesssim 2^{(2+\beta + 2\delta)m} \cdot 2^k (2^q + 2^{2p_2})^{-1} 2^{k_2+p_2+ \frac{k+q}{2}} \cdot ||f_1||_2 \cdot 2^{-\frac{3}{2}m + \gamma m} \epsilon^2 \\
	& \lesssim 2^{(\frac{1}{2} + \beta + \gamma + 2\delta)m} \cdot 2^{k_2 + \frac{3}{2}k} \cdot 2^{p_1 + \frac{m}{3}} 2^{-l_1}||f_1||_X \epsilon^2 \\
	& \lesssim 2^{(-\frac{1}{6} + \beta + \gamma + 3\delta + 3\delta_0)m} \epsilon^3,
\end{align*}
and this completes the nonresonant term.

For the resonant term, we use the estimate using the vertical variant Lemma \ref{vertical-variant}, since $|\partial_3 \Phi| \sim 2^{2p_2 - k_2}$. Then, integration by parts in $|\eta|\partial_{\eta_3}$ gives enough result if
\begin{equation*}
	2^{-2p_2}[1+2^{l_1+p_1} + 2^{l_2+p_2}] \leq 2^{(1-\delta)m},
\end{equation*}
as $|\eta|\partial{\eta_3}$ commutes with $e^{(t-s)\kap\Delta}$ so that we can use the same proof. Also, note that we have $q \sim 2p_2$ here, as
\begin{equation*}
	|\Lambda(\xi-\eta)| - |\Lambda(\eta)| = \frac{\sqrt{1 - \Lambda^2(\eta)}^2 - \sqrt{1 - \Lambda^2(\xi-\eta)}^2}{|\Lambda(\xi-\eta)| + |\Lambda(\eta)|} \in [2^{-6} 2^{2p_2}, 2^3 2^{2p_2}],
\end{equation*}
so that $|\Phi|$ cannot be smaller than our claimed bound $\lambda = 2^{-20}(2^{q} + 2^{2p_2})$ unless $|\Lambda(\xi)| \sim 2^{q}$ matches the size of $2^{2p_2}$.

We first look at the case where $l_1 + p_1 - 2p_2 \geq (1-\delta)m$. Then,
\begin{align*}
	2^{(1+\beta)l}2^{\beta p} \Vert P_{k,p}R_l \mathcal{B}_{\mathfrak{n}^{res}}(f_1, f_2) \Vert_{L^2}
	& \lesssim 2^{(2+\beta + 2\delta)m} \cdot 2^k \cdot 2^{\frac{3}{2}k + \frac{1}{2}q} \cdot 2^{-(1+\beta)l_1 - \beta p_1} ||f_1||_X \cdot 2^{-(1+\beta)l_2 - \beta p_2} ||f_2||_X \\
	& \lesssim 2^{(2+\beta + 2\delta)m} \cdot 2^{\frac{5}{2}k + p_2} \cdot 2^{p_1 - 2(1+\beta)p_2 - (1+\beta)(1-\delta)m} 2^{-(1+\beta)(p_2 - k_1 + k + (1-\delta)m) - \beta p_2} \epsilon^2 \\
	& \lesssim 2^{(-\beta + (4+2\beta)\delta)m} \cdot 2^{\frac{5}{2}(k_2 + p_2)} \cdot 2^{p_1 - p_2} \cdot 2^{-(2+5\beta)p_2} \epsilon^2 \\
	& \lesssim 2^{(-\beta + (4+2\beta)\delta + 3\delta_0)m} 2^{(\frac{1}{2} - 5\beta)p_2} \epsilon^2,
\end{align*}
gives acceptable bound since $\beta \leq \frac{1}{10}$. Here, we used $\frac{q}{2} \sim p_2$ and $-k_1 + k \sim p_2$.

Finally, we assume $l_2 - p_2 \geq (1-\delta)m$, in which case,
\begin{align*}
	2^{(1+\beta)l}2^{\beta p} \Vert P_{k,p}R_l \mathcal{B}_{\mathfrak{n}^{res}}(f_1, f_2) \Vert_{L^2}
	& \lesssim 2^{(2+\beta + 2\delta)m} \cdot 2^k 2^{k_1+p_1 + \frac{k+q}{2}} \cdot 2^{-(1+\beta)l_1-\beta p_1} \cdot 2^{-(1+\beta)l_2 - \beta p_2} \epsilon^2 \\
	& \lesssim 2^{(2+\beta + 2\delta)m} \cdot 2^{k_1 + \frac{3}{2}k + p_2} \cdot 2^{p_1-l_1} 2^{-\beta(p_1+l_1)} \cdot 2^{-(1+2\beta)p_2 - (1+\beta)(1-\delta)m} \epsilon^2 \\
	& \lesssim 2^{5\delta m} \cdot 2^{2k_1 + \frac{k}{2}} 2^{-\beta(p_1+l_1)} 2^{-2\beta p_2} \epsilon^2 \\
	& \sim 2^{5\delta m} \cdot 2^{\frac{5}{2}k_1} 2^{-\beta(p_1+l_1)} 2^{(\frac{1}{2} - 2\beta)p_2} \epsilon^2.
\end{align*}
Thus, if $p_2 \leq -20\delta m$ or $p_1 + l_1 \geq 400\delta m$, we are done. If not, we perform another normal form using $L = 2^{-300\delta m}$. One can see we are dividing with smaller cutoff, as the former cutoff is of order $2^{2p_2} \geq 2^{-40\delta m}$ with the newest assumption. From $|\partial_3 \Phi| \gtrsim 2^{2p_2 - k_2} \geq 2^{-100\delta m}$ and the set size estimates,
\begin{align*}
	2^{(1+\beta)l}2^{\beta p} \Vert P_{k,p}R_l \mathcal{B}_{\mathfrak{n}^{rr}}(f_1, f_2) \Vert_{L^2}
	& \lesssim 2^{(2+\beta + 2\delta)m} \cdot 2^k \cdot (L 2^{100\delta m})^{\frac{1}{2}} \cdot 2^{k_1+p_1} \cdot 2^{-l_1} ||f_1||_X \cdot 2^{-(1+\beta)l_2 - \beta p_2} ||f_2||_X \\
	& \lesssim 2^{(2+\beta + 2\delta)m} 2^{-100\delta m} 2^{2k_1} 2^{p_1 + k - k_1 - l_1} \cdot 2^{-(1+\beta)(p_2 - l_2) - (1+2\beta)p_2} \epsilon^2 \\
	& \lesssim 2^{(2+\beta + 2\delta)m} 2^{-100\delta m} 2^{2\delta_0 m} 2^{-(1-\delta)m} 2^{-(1+\beta)(1-\delta)m + (1+2\beta)20\delta m} \epsilon^2 \\
	& \leq 2^{-70\delta m} \epsilon^2,
\end{align*}
\begin{align*}
	2^{(1+\beta)l}2^{\beta p} \Vert P_{k,p}R_l \mathcal{Q}_{\mathfrak{n}^{nrr}}(f_1, f_2) \Vert_{L^2}
	& \lesssim 2^{(1+\beta + 2\delta)m} \cdot 2^{k + 300\delta m} 2^{k_1+p_1 + \frac{k+q}{2}} \cdot 2^{-l_1} ||f_1||_X \cdot 2^{-l_2} ||f_2||_X \\
	& \lesssim 2^{(1+\beta + 302\delta)m} \cdot 2^{2k_1} \cdot 2^{p_1+k-k_1-l_1} 2^{p_2-l_2} \epsilon^2 \\
	& \lesssim 2^{(-1 + \beta + 304\delta + 2\delta_0)m} \epsilon^2,
\end{align*}
\begin{align*}
	2^{(1+\beta)l}2^{\beta p} \Vert P_{k,p}R_l \mathcal{B}_{\mathfrak{n}^{nrr}}((\dt-\kap\Delta)f_1, f_2) \Vert_{L^2}
	& \lesssim 2^{(2+\beta + 2\delta)m} \cdot 2^{k+300\delta m} \cdot 2^{k_2+p_2} 2^{-\frac{3}{2}m + \gamma m} \epsilon^2 \cdot 2^{-l_2} ||f_2||_X \\
	& \lesssim 2^{(-\frac{1}{2} + \beta + \gamma + 303\delta + 2\delta_0)m} \epsilon^3,
\end{align*}
\begin{align*}
	2^{(1+\beta)l}2^{\beta p} \Vert P_{k,p}R_l \mathcal{B}_{\mathfrak{n}^{nrr}}(f_1, (\dt-\kap\Delta)f_2) \Vert_{L^2}
	& \lesssim 2^{(2+\beta + 2\delta)m} \cdot 2^{k + 300\delta m} \cdot 2^{k_1+p_1} \cdot 2^{-l_1} ||f_1||_X \cdot 2^{-\frac{3}{2}m + \gamma m} \epsilon^2 \\
	& \lesssim 2^{(-\frac{1}{2} + \beta + \gamma + 303\delta + 2\delta_0)m} \epsilon^3.
\end{align*}
This finishes the case \textcircled{2}, $p_1 \ll p, p_2$.\\

\noindent
\textcircled{3} $p \sim p_1 \ll p_2 \sim 0$ : By Lemma \ref{two_small_sizes}, we get $k_2 \ll k_1 \sim k$ and $k_2+p_2 \lesssim k_1+p_1$. Integration by parts are enough if
\begin{gather*}
	2^{-p_1+2k_1 - k_{\max} - k_{\min}}(1 + 2^{k_2 - k_1 + l_1}) \sim 2^{-p_1}(2^{k_1 - k_2} + 2^{l_1}) \leq 2^{(1-\delta)m}, \quad \text{or,} \\
	2^{-p_2+2k_2 - k_{\max} - k_{\min}}(1 + 2^{k_1 - k_2 + l_2}) \sim 2^{-p_2 + l_2} \leq 2^{(1-\delta)m}.
\end{gather*}
Otherwise, if $k_1 - k_2 \leq l_1$, we have $p_j - l_j \leq -(1-\delta)m$ for $j=1,2$, and hence
\begin{align*}
	2^{(1+\beta)l}2^{\beta p} \Vert P_{k,p}R_l \mathcal{B}_\mathfrak{m}(f_1, f_2) \Vert_{L^2}
	& \lesssim 2^{(2+\beta + 2\delta)m} 2^{\beta p} \cdot 2^k \cdot 2^{\frac{3}{2}k_2 + p_2} \cdot 2^{-(1+\beta)l_1 - \beta p_1} ||f_1||_X \cdot 2^{-(1+\beta)l_2} ||f_2||_X \\
	& \lesssim 2^{(2+\beta + 2\delta)m} 2^k 2^{\frac{3}{2}(k_1+p_1)} 2^{-(1+\beta)l_1} 2^{-(1+\beta)(p_2 - l_2)} \epsilon^2 \\
	& \lesssim 2^{(2+\beta + 2\delta)m} 2^{\frac{5}{2}k + (\frac{1}{2} - \beta)p_1} 2^{-2(1+\beta)(1-\delta)m} \epsilon^2 \\
	& \lesssim 2^{(-\beta + 5\delta + 3\delta_0)m} \epsilon^2.
\end{align*}
If $k_1 - k_2 \geq l_1$, we have $p_1 - k_1 + k_2 \leq -(1-\delta)m$, and hence,
\begin{align*}
	2^{(1+\beta)l}2^{\beta p} \Vert P_{k,p}R_l \mathcal{B}_\mathfrak{m}(f_1, f_2) \Vert_{L^2}
	& \lesssim 2^{(2+\beta + 2\delta)m} 2^{\beta p_1} \cdot 2^k \cdot 2^{\frac{3}{2}k_2 + p_2} \cdot 2^{p_1} ||f_1||_B \cdot 2^{-(1+\beta)l_2} ||f_2||_X \\
	& \lesssim 2^{(2+\beta + 2\delta)m} \cdot 2^{k + (1+\beta)k_1 + (\frac{1}{2} - \beta)k_2} \cdot 2^{(1+\beta)(p_1 - k_1 + k_2)} 2^{(1+\beta)(p_2 - l_2)} \epsilon^2 \\
	& \lesssim 2^{(-\beta + 5\delta + 3\delta_0)m} \epsilon^2.
\end{align*}

\noindent
\textcircled{4} $p_1 \sim p_2 \ll p \sim 0$ : Again by Lemma \ref{two_small_sizes}, we have $k \ll k_1 \sim k_2$ and $k+p \lesssim k_1+p_1 \sim k_2+p_2$. Integration by parts is sufficient if
\begin{gather}
	2^{-p_1+2k_1 - k_{\max} - k_{\min}}(1 + 2^{k_2 - k_1 + l_1}) \sim 2^{-p_1+k_1-k+l_1} \leq 2^{(1-\delta)m}, \quad \text{or,} \label{ibyp-gapinp-last} \\
	2^{-p_2+2k_2 - k_{\max} - k_{\min}}(1 + 2^{k_1 - k_2 + l_2}) \sim 2^{-p_2+k_2-k+l_2} \leq 2^{(1-\delta)m}.
\end{gather}
Hence, we assume $p_j +k - l_j \leq k_1 - (1-\delta)m, j=1,2$ from now on. We do another integration by parts in vertical variant, Lemma \ref{vertical-variant}. From
\begin{equation*}
	|\partial_{\eta_3} \Phi| = |\partial_{\eta_3}\Lambda(\eta) \pm \partial_{\eta_3}\Lambda(\xi-\eta)| = \left\vert \frac{|\eta_h|^2}{|\eta|^3} \mp \frac{|\xi_h-\eta_h|^2}{|\xi-\eta|^3} \right\vert,
\end{equation*}
we can see that the smaller case comes from \eqref{ibyp-gapinp-last}. In this case,
\begin{equation*}
	|\partial_{\eta_3}\Phi|
	= |\partial_{\eta_3}\Lambda(\eta) +\partial_{\eta_3}\Lambda(\xi-\eta)|
	= |\partial_{\eta_3}\Lambda(\eta) - \partial_{\eta_3}\Lambda(\eta - \xi)|
	\sim |\partial_3 \nabla \Lambda(\eta) \cdot \xi|
	\sim 2^{p_2 + k - 2k_2}.
\end{equation*}
Hence, using the vertical variant of integration by parts with respect to $D_3^{\eta} = |\eta|\partial_{\eta_3}$, we get an acceptable contribution if
\begin{equation*}
	2^{-p_2 - k + k_2}(1 + 2^{p_1+l_1} + 2^{p_2 + l_2}) \sim 2^{-k + k_1}(2^{l_1} + 2^{l_2}) \leq 2^{(1-\delta)m}.
\end{equation*}
For the remaining regions, WLOG we now assume $k - k_1 - l_2 \leq -(1-\delta)m$, because the relations on indices are symmetric in $\xi-\eta$ and $\eta$ up to this point.

Now, since $|\partial_{\xi_3} \Phi| \sim |2^{2p-k} \pm 2^{2p_1-k_1}| \sim 2^{-k}$, we can use the set size estimate. Combining this with the normal form with $\lambda = 2^{-\frac{2}{3}m}$
\begin{align*}
	2^{(1+\beta)l}2^{\beta p} \Vert P_{k,p}R_l \mathcal{B}_{\mathfrak{n}^{res}}(f_1, f_2) \Vert_{L^2}
	& \lesssim 2^{(2+\beta + 2\delta)m} \cdot 2^k \cdot (\lambda 2^k)^{\frac{1}{2}} \cdot 2^{k+p} \cdot ||f_1||_2 ||f_2||_2 \\
	& \lesssim 2^{(2+\beta + 2\delta)m} \cdot 2^{\frac{5}{2}k} 2^{-\frac{1}{3}m} \cdot (2^{p_1} ||f_1||_B)^{\frac{1}{4}} (2^{-l_1} ||f_1||_X)^{\frac{3}{4}} \cdot 2^{-l_2} ||f_2||_X \\
	& \lesssim 2^{(\frac{5}{3} + \beta + 2\delta)m} \cdot 2^{\frac{3}{4}(k + k_1+p_1 - l_1)} 2^{\frac{1}{4}p_1} \cdot 2^{k - l_2} \epsilon^2 \\
	& \lesssim 2^{(-\frac{1}{12} + \beta + 4\delta + 3\delta_0)m} \epsilon^2,
\end{align*}
where we used $k \lesssim k_1 + p_1$, and
\begin{align*}
	2^{(1+\beta)l}2^{\beta p} \Vert P_{k,p}R_l \mathcal{Q}_{\mathfrak{n}^{nr}}(f_1, f_2) \Vert_{L^2}
	& \lesssim 2^{(1+\beta + 2\delta)m} \cdot 2^k \cdot (\lambda^{-1} 2^k)^{\frac{1}{2}} \cdot 2^{k + p} \cdot ||f_1||_2 ||f_2||_2 \\
	& \lesssim 2^{(1+\beta + 2\delta)m} \cdot 2^{\frac{5}{2}k} 2^{\frac{1}{3}m} 2^{-\frac{3}{4}l_1 + \frac{1}{4}p_1} \cdot 2^{-l_2} \epsilon^2 \\
	& \lesssim 2^{(\frac{4}{3} + \beta + 2\delta)m} \cdot 2^{\frac{3}{4}(k + k_1+p_1 - l_1)} 2^{\frac{1}{4}p_1} \cdot 2^{k - l_2} \epsilon^2 \\
	& \lesssim 2^{(-\frac{5}{12} + \beta + 4\delta + 3\delta_0)m} \epsilon^2,
\end{align*}
\begin{align*}
	2^{(1+\beta)l}2^{\beta p} \Vert P_{k,p}R_l \mathcal{B}_{\mathfrak{n}^{nr}}((\dt-\kap\Delta)f_1, f_2) \Vert_{L^2}
	& \lesssim 2^{(2+\beta + 2\delta)m} \cdot 2^k \cdot (\lambda^{-1} 2^k)^{\frac{1}{2}} \cdot 2^{k + p} \cdot 2^{-\frac{3}{2}m + \gamma m} \epsilon^2 \cdot ||f_2||_2 \\
	& \lesssim 2^{(\frac{1}{2} + \beta + \gamma + 2\delta)m} \cdot 2^{\frac{5}{2}k} 2^{\frac{1}{3}m} 2^{-l_2} ||f_2||_X \epsilon^2 \\
	& \lesssim 2^{(-\frac{1}{6} + \beta + \gamma + 3\delta + 3\delta_0)m} \epsilon^3,
\end{align*}
\begin{align*}
	2^{(1+\beta)l}2^{\beta p} \Vert P_{k,p}R_l \mathcal{B}_{\mathfrak{n}^{nr}}(f_1, (\dt-\kap\Delta)f_2) \Vert_{L^2}
	& \lesssim 2^{(2+\beta + 2\delta)m} \cdot 2^k \cdot (\lambda^{-1} 2^k)^{\frac{1}{2}} \cdot 2^{k + p} \cdot ||f_1||_2 \cdot 2^{-\frac{3}{2}m + \gamma m} \epsilon^2 \\
	& \lesssim 2^{(\frac{1}{2} + \beta + \gamma + 2\delta)m} \cdot 2^{\frac{5}{2}k} 2^{\frac{1}{3}m} 2^{-l_1} ||f_1||_X \epsilon^2 \\
	& \lesssim 2^{(\frac{5}{6} + \beta + \gamma + 2\delta)m} 2^{\frac{1}{2}k} 2^{k+k_1+p_1-l_1} \epsilon^3 \\
	& \lesssim 2^{(-\frac{1}{6} + \beta + \gamma + 3\delta + 3\delta_0)m} \epsilon^3.
\end{align*}
This finishes the final case for gap in $p$ with $p_{\max} \sim 0$.\\

\noindent \textbf{7.2.2 Gap in $p$ with $p_{\max} \ll 0$} \\
In this case, every $\Lambda(\zeta) \geq \frac{1}{2}, \zeta \in \{\xi,\xi-\eta,\eta\}$, so that we have $|\Phi| \geq \frac{1}{2}$, enabling us to perform normal form with $|\mathfrak{n}| \lesssim 2^k$ every time. The terms with time derivative are easier to bound by
\begin{align*}
	2^{(1+\beta)l}2^{\beta p} \Vert P_{k,p}R_l \mathcal{B}_{\mathfrak{n}^{nr}}((\dt-\kap\Delta)f_1, f_2) \Vert_{L^2}
	& \lesssim 2^{(2+\beta + 2\delta)m} \cdot 2^k \cdot ||(\dt-\kap\Delta)f_1||_2 \cdot ||e^{it\Lambda} f_2||_\infty \\
	& \lesssim 2^{(2+\beta + 2\delta)m} \cdot 2^k \cdot 2^{-\frac{3}{2}m + \gamma m} \epsilon^2 \cdot 2^{-m + \gamma m} \epsilon \\
	& \lesssim 2^{(-\frac{1}{2} + \beta + 2\gamma + 2\delta + \delta_0)m} \epsilon^3,
\end{align*}
where we used Lemma \ref{slower-dispersive-decay}. The symmetric term works the same, and the boundary term is the toughest. WLOG assume $f_1$ has less vector field applied than to $f_2$, and first consider the case when $p_1 \gtrsim p_2$. If $p_2 \leq -\frac{m}{2}$,
\begin{align*}
	2^{(1+\beta)l}2^{\beta p} \Vert P_{k,p}R_l \mathcal{Q}_{\mathfrak{n}^{nr}}(f_1, f_2) \Vert_{L^2}
	& \lesssim 2^{(1+\beta + 2\delta)m} \cdot 2^k \cdot ||e^{it\Lambda} f_1||_\infty \cdot 2^{p_2}||f_2||_B \\
	& \lesssim 2^{(1+\beta + 2\delta)m} \cdot 2^k 2^{-m + \gamma m} \epsilon \cdot 2^{-\frac{m}{2}} \epsilon \\
	& \leq 2^{(-\frac{1}{2} + \beta + \gamma + 2\delta + \delta_0)m} \epsilon^2,
\end{align*}
and when $p_2 \geq -\frac{m}{2}$, we divide $f_1$ into $f_1^I$ and $f_1^{II}$ as in Proposition \ref{dispersive-decay} to achieve
\begin{align*}
	2^{(1+\beta)l}2^{\beta p} \Vert P_{k,p}R_l \mathcal{Q}_{\mathfrak{n}^{nr}}(f_1^I, f_2) \Vert_{L^2}
	& \lesssim 2^{(1+\beta + 2\delta)m} \cdot 2^k \cdot ||e^{it\Lambda} f_1^I||_\infty \cdot ||f_2||_2 \\
	& \lesssim 2^{(1+\beta + 2\delta)m} \cdot 2^k \cdot 2^{-p_1-\frac{3}{2}m} ||f_1^I||_D \cdot 2^{p_2} ||f_2|||_B \\
	& \lesssim 2^{(-\frac{1}{2} + \beta + 2\delta + \delta_0)m} \epsilon^2,
\end{align*}
\begin{align*}
	2^{(1+\beta)l}2^{\beta p} \Vert P_{k,p}R_l \mathcal{Q}_{\mathfrak{n}^{nr}}(f_1^{II}, f_2) \Vert_{L^2}
	& \lesssim 2^{(1+\beta + 2\delta)m} \cdot 2^k \cdot ||e^{it\Lambda} f_1^{II}||_2 \cdot ||e^{it\Lambda}f_2||_\infty \\
	& \lesssim 2^{(1+\beta + 2\delta)m} \cdot 2^k \cdot 2^{(-1-\beta')m - (1-2\beta')p_1} ||f_1^I||_D \cdot 2^{-m+\gamma m} \epsilon \\
	& \lesssim 2^{(\beta + \gamma + 2\delta + \delta_0)m} 2^{(-1-\beta')m + (\frac{1}{2} + \beta')m} \epsilon^2 \\
	& \lesssim 2^{(-\frac{1}{2} + \beta + 2\delta + \delta_0)m} \epsilon^2.
\end{align*}
Now when $p_1 \ll p_2$, from the fact that $\bar{\sigma} \sim 2^{p_{\max} + k_{\max} + k_{\min}}$, integration by parts is enough if
\begin{gather*}
	2^{-p_1 - p_{\max} + 2k_1 - k_{\max} - k_{\min}}(1 + 2^{k_2 - k_1 + l_1}) \leq 2^{(1-\delta)m}
\end{gather*}
Hence, we are left with the case when
\begin{gather*}
	p_1 + p_{\max} + k_{\max} + k_{\min} - 2k_1 \leq -(1-\delta)m, \quad \text{or,} \\
	p_1 + p_{\max} + k_{\max} + k_{\min} - k_1 - k_2 - l_1 \leq -(1-\delta)m.
\end{gather*}
Since $k$ appears from $\mathfrak{m}$ and $k_2$ can be obtained when using Lemma \ref{slower-dispersive-decay}, we can always have one copy of $k_{\min}$.(If $k_1 = k_{\min}$, it will cancel out with $-k_1$) Therefore,
\begin{align*}
	2^{(1+\beta)l}2^{\beta p} \Vert P_{k,p}R_l \mathcal{Q}_{\mathfrak{n}^{nr}}(f_1, f_2) \Vert_{L^2}
	& \lesssim 2^{(1+\beta + 2\delta)m} \cdot 2^{k+p_{\max}} (2^{p_1}||f_1||_B)^{\frac{1}{2}} (2^{-l_1}||f_1||_X)^{\frac{1}{2}} \cdot ||e^{it\Lambda} f_2||_\infty \\
	& \lesssim 2^{(1+\beta + 2\delta)m} \cdot 2^{\frac{k+p_{\max}}{2}} 2^{\frac{1}{2}(p_1 + p_{\max} + k - l_1)} \cdot 2^{\frac{3}{2}k_2 - 3k_2^+ -m + \gamma m} \epsilon^2 \\
	& \lesssim 2^{(-\frac{1}{2} + \beta + \gamma + 3\delta + \delta_0)m} \epsilon^2.
\end{align*}
This finishes the case with $p_{\max} \ll 0$.\\

\noindent \textbf{7.2.3 Gap in $q$} \\
We can assume $p \sim p_1 \sim p_2 \sim 0$ from now on. We further localize in $q$ and write $g_i = P_{k_i,p_i,q_i}\mathcal{U}_i$, $i=1,2$, and look at the case where we have gap in $q$. Using the boostrap assumption on $B$-norm,
\begin{align*}
	2^{(1+\beta)l}2^{\beta p} \Vert P_{k,p}R_l \mathcal{B}_\mathfrak{m}(g_1, g_2) \Vert_{L^2}
	& \lesssim 2^{(1+\beta+2\delta)m} \cdot 2^m \cdot 2^{k+q_{\max}} 2^{\frac{3}{2}k_{\max} + \frac{1}{2}q_{\min}} \cdot 2^{\frac{1}{2}q_1} ||g_1||_B \cdot 2^{\frac{1}{2}q_2} ||g_2||_B,
\end{align*}
we can reduce to the case when $q_{\max} \geq -\frac{14}{15}m$ and $q_{\min} \geq -5m$. This means we are again dealing with only $(\log \langle t \rangle)^6$ many cases, and hence it's enough to prove
\begin{align*}
	2^{(1+\beta)l} \Vert P_{k,p}R_l \mathcal{B}_\mathfrak{m}(f_1, f_2) \Vert_{L^2}
	& \lesssim 2^{-\delta m}.
\end{align*}

\noindent\textcircled{1} $q \ll q_1, q_2$: We divide the cases according to Lemma \ref{two_big_sizes} again.\\
1-1) $k_1 \sim k_2$: This implies $q_1 \sim q_2$ and hence integration by parts is sufficient when $-q_1 + k_1 - k + l_i \leq (1-\delta)m$, $i=1,2$. Now that we can assume $p_1,p_2 \sim 0$, we actively use the precise decay estimates. Without loss of generality, assuming $g_2$ has less copies of $S$ and denoting $g_2 = g_2^{I} + g_2^{II} := \ I + II$ as in Proposition \ref{dispersive-decay},
\begin{align*}
	2^{(1+\beta)l} \Vert P_{k,p}R_l \mathcal{B}_\mathfrak{m}(g_1, g_2^{I}) \Vert_{L^2}
	& \lesssim 2^{(2+\beta + 2\delta)m} \cdot 2^{k+q_{\max}} ||g_1||_2 ||e^{it\Lambda} g_2^{I}||_\infty \\
	& \lesssim 2^{(2+\beta + 2\delta)m} \cdot 2^{k+q_2} (2^{q_1/2} ||g_1||_B)^{1/4} (2^{-(1+\beta)l_1} ||g_1||_X) \cdot 2^{-\frac{q_2}{2} - \frac{3}{2}m} \epsilon \\
	& \lesssim 2^{(\frac{1}{2}+\beta + 2\delta + \delta_0)m} \cdot 2^{\frac{3}{4}(q_2 - l_1)} \epsilon^2 \\
	& \lesssim 2^{(-\frac{1}{4} + \beta + 3\delta)m} \epsilon^2,
\end{align*}
\begin{align*}
	2^{(1+\beta)l} \Vert P_{k,p}R_l \mathcal{B}_\mathfrak{m}(g_1, g_2^{II}) \Vert_{L^2}
	& \lesssim 2^{(2+\beta + 2\delta)m} \cdot 2^{k+q_1} \cdot 2^{k+\frac{1}{2}k_1 + \frac{1}{2}q_1} \cdot 2^{-(1+\beta)l_1} ||g_1||_X \cdot ||e^{it\Lambda} g_2^{II}||_2 \\
	& \lesssim 2^{(2+\beta + 2\delta)m} \cdot 2^{(1-\beta)k+\frac{k_1}{2} + \frac{1-2\beta}{2}q_1} \cdot 2^{(1+\beta)(q_1 + k - l_1)} 2^{-(1+\beta')m} ||g_2||_D \\
	& \lesssim 2^{(-\beta' + 4\delta)m} \epsilon^2,
\end{align*}
give satisfactory bounds. \\

\noindent 1-2) $k_2 \ll k_1$: From $2^{k_1+q_1} \sim 2^{k_2+q_2}$, we also have $2^{q_1 - q_2} \sim 2^{k_2 - k_1} \ll 1$, and $q \ll q_1 \ll q_2$. Integration by parts is sufficient if
\begin{gather*}
	2^{-q_2+2k_1 - k_{\max} - k_{\min}}(1 + 2^{k_2 - k_1}(2^{q_2 - q_1} + 2^{l_1})) \sim 2^{-q_1} + 2^{-q_2+l_1} \leq 2^{(1-\delta)m}, \quad \text{or,} \\
	2^{-q_2+2k_2 - k_{\max} - k_{\min}}(1 + 2^{k_1 - k_2}(2^{q_1 - q_2} + 2^{l_2})) \sim 2^{-q_2 + l_2} \leq 2^{(1-\delta)m}.
\end{gather*}
So we look at the converse situation, and when $q_1 \geq (1+2\beta)q_2$ in addition,
\begin{align*}
	2^{(1+\beta)l} \Vert P_{k,p}R_l \mathcal{B}_\mathfrak{m}(g_1^I, g_2) \Vert_{L^2}
	& \lesssim 2^{(2+\beta + 2\delta)m} \cdot 2^{k+q_2} \cdot 2^{-\frac{3}{2}m - \frac{q_1}{2}} ||e^{it\Lambda} g_1^I||_D \cdot 2^{-(1+\beta)l_2} ||g_2||_X \\
	& \lesssim 2^{(\frac{1}{2}+\beta+2\delta + \delta_0)m} 2^{(1+\beta)(q_2 - l_2)} 2^{-\beta q_2 - \frac{q_1}{2}} \epsilon^2 \\
	& \lesssim 2^{(-\frac{1}{2} + 4\delta)m} 2^{(\frac{1}{2}+2\beta)\frac{14}{15}m} \epsilon^2
\end{align*}
gives sufficient bound since $\beta \leq \frac{1}{60}$, and
\begin{align*}
	2^{(1+\beta)l} \Vert P_{k,p}R_l \mathcal{B}_\mathfrak{m}(g_1^{II}, g_2) \Vert_{L^2}
	& \lesssim 2^{(2+\beta + 2\delta)m} \cdot 2^{k+q_2} \cdot 2^{\frac{3}{2}k_2 + \frac{1}{2}q_2} \cdot 2^{-(1+\beta')m} ||g_1||_D \cdot 2^{-(1+\beta)l_2} ||g_2||_X \\
	& \lesssim 2^{(1-\beta' + \beta + 2\delta + 3\delta_0)m} 2^{(1+\beta)(q_2 - l_2) + \frac{1-2\beta}{2}q_2} \epsilon^2 \\
	& \lesssim 2^{(-\beta' + 4\delta)m} \epsilon^2,
\end{align*}
can be used when $g_1$ has less power of $S$. When $g_2$ has less power of $S$,
\begin{align*}
	2^{(1+\beta)l} \Vert P_{k,p}R_l \mathcal{B}_\mathfrak{m}(g_1, g_2^{I}) \Vert_{L^2}
	& \lesssim 2^{(2+\beta + 2\delta)m} \cdot 2^{k+q_2} \cdot 2^{-(1+\beta)l_1} ||g_1||_X \cdot 2^{\frac{3}{2}k_2 -\frac{3}{2}m - \frac{q_2}{2}} ||g_2||_D \\
	& \lesssim 2^{(\frac{1}{2}+\beta+2\delta + \delta_0)m} \cdot 2^{\frac{1-2\beta}{2}k_2 - (\frac{1}{2}+\beta)q_2} 2^{(1+\beta)(q_2 + k_2 - l_1)} \epsilon^2 \\
	& \lesssim 2^{(-\frac{1}{2} + 4\delta)m} 2^{(\frac{1}{2}+\beta)\frac{14}{15}m} \epsilon^2,
\end{align*}
gives sufficient bound since $q_2+k_2 \sim q_1+k_1$ and $\beta \leq \frac{1}{30}$, and
\begin{align*}
	2^{(1+\beta)l} \Vert P_{k,p}R_l \mathcal{B}_\mathfrak{m}(g_1, g_2^{II}) \Vert_{L^2}
	& \lesssim 2^{(2+\beta + 2\delta)m} \cdot 2^{k+q_2} \cdot 2^{k_2 + \frac{k_1+q_1}{2}} \cdot 2^{-(1+\beta)l_1} ||g_1||_X \cdot 2^{-(1+\beta')m} ||g_2||_D \\
	& \lesssim 2^{(1-\beta'+\beta+2\delta)m} 2^k 2^{\frac{3}{2}(k_2+q_2) - (1+\beta)l_1} \epsilon^2 \\
	& \lesssim 2^{(-\beta' + 4\delta)m} \epsilon^2,
\end{align*}
using $k_1 + q_1 \sim k_2 + q_2$. \\
Now we look at the case when $q_1 \leq (1+2\beta)q_2$. Here, we use $k_2 \sim k_1 + q_1 - q_2 \leq k_1 + 2\beta q_2$, so that we can generate $q_2$ using $k_2$. Using the normal form since $|\Phi| \gtrsim 2^{q_2}$,
\begin{align*}
	2^{(1+\beta)l} \Vert P_{k,p}R_l \mathcal{Q}_\mathfrak{n}(g_1, g_2) \Vert_{L^2}
	& \lesssim 2^{(1+\beta + 2\delta)m} \cdot 2^{k} \cdot 2^{\frac{3}{2}k_2 + \frac{1}{2}q_2} (2^{\frac{1}{2}q_1} ||g_1||_B)^{1 - \frac{4\beta}{3}} (2^{-(1+\beta)l_1} ||g_1||_X)^{\frac{4\beta}{3}} \cdot 2^{-(1+\beta)l_2} ||g_2||_X \\
	& \lesssim 2^{(1+\beta + 2\delta)m} \cdot 2^{\frac{5}{2}k} \cdot 2^{(\frac{1}{2}+3\beta)q_2} \cdot 2^{(\frac{1}{2}-2\beta)q_2 + \frac{4\beta}{3}q_1 - \frac{4\beta}{3}(1+\beta)l_1} \cdot 2^{-(1+\beta)l_2} \epsilon^2 \\
	& \lesssim 2^{(-\beta + 4\delta + 3\delta_0)m} \epsilon^2,
\end{align*}
\begin{align*}
	2^{(1+\beta)l} \Vert P_{k,p}R_l \mathcal{B}_\mathfrak{n}((\dt-\kap\Delta)g_1, g_2) \Vert_{L^2}
	& \lesssim 2^{(2+\beta + 2\delta)m} \cdot 2^k \cdot 2^{\frac{3}{2}k_2 + \frac{1}{2}q_2} \cdot ||(\dt-\kap\Delta)g_1||_2 \cdot 2^{-(1+\beta)l_2} ||g_2||_X \\
	& \lesssim 2^{(2+\beta + 2\delta)m} \cdot 2^k \cdot 2^{\frac{3}{2}k_1 + (\frac{1}{2}+3\beta)q_2} \cdot 2^{-\frac{3}{2}m + \gamma m} \epsilon^2 \cdot 2^{-(1+\beta)l_2} ||g_2||_X \\
	& \lesssim 2^{(-2\beta + \gamma + 3\delta)m} \epsilon^3,
\end{align*}
are acceptable contributions since $\gamma \ll \beta$, and
\begin{align*}
	2^{(1+\beta)l} & \Vert P_{k,p}R_l \mathcal{B}_\mathfrak{m}(g_1, (\dt-\kap\Delta)g_2) \Vert_{L^2} \\
	& \lesssim 2^{(2+\beta + 2\delta)m} \cdot 2^k \cdot 2^{\frac{3}{2}k_2 + \frac{1}{2}q_2} \cdot (2^{\frac{q_1}{2}} ||g_1||_B)^{\frac{1}{3}} (2^{-(1+\beta)l_1} ||g_1||_X)^{\frac{2}{3}} \cdot ||(\dt-\kap\Delta)g_2||_2 \\
	& \lesssim 2^{(2+\beta + 2\delta)m} \cdot 2^{2k_2 + \frac{1}{2}k_1 + \frac{1}{2}q_1} \cdot 2^{\frac{1}{6}q_1 - \frac{2}{3}l_1} \cdot 2^{-\frac{3}{2}m + \gamma m} \epsilon^3 \\
	& \lesssim 2^{(-\frac{1}{6} + \beta + \gamma + 3\delta)m} \epsilon^3
\end{align*}
hence we get acceptable contributions for the all terms after normal form transformation.\\

\noindent 1-3) $k_1 \ll k_2$: This forces $q_2 \ll q_1$ which is against the assumption, and hence excluded. \\

\noindent \textcircled{2} $q_1 \ll q_2, q$ : Again, we divide the cases by Lemma \ref{two_big_sizes}.\\
2-1) $k \sim k_2$ : This implies $q \sim q_2$ and hence integration by parts is enough if
\begin{gather}
	2^{-q_2+k_1 - k_2}(1 + 2^{k_2 - k_1}(2^{q_2 - q_1} + 2^{l_1})) \sim 2^{-q_1} + s^{-q_2 + l_1} \leq 2^{(1-\delta)m}, \quad \text{or,} \label{ibyp2-1-1} \\
	2^{-q_2+k_2 - k_1}(1 + 2^{k_1 - k_2}(2^{q_1 - q_2} + 2^{l_2})) \sim 2^{-q_2 + k_2 - k_1} + 2^{-q_2 + l_2} \leq 2^{(1-\delta)m} \label{ibyp2-1-2}.
\end{gather}
So we look at the other cases. Assume first that \eqref{ibyp2-1-2} does not hold because $q_2 + k_1 - k_2 \leq -(1-\delta)m$. Then,
\begin{align*}
	2^{(1+\beta)l} \Vert P_{k,p}R_l \mathcal{B}_\mathfrak{m}(g_1, g_2) \Vert_{L^2}
	& \lesssim 2^{(2+\beta + 2\delta)m} \cdot 2^{k+q_2} \cdot 2^{\frac{3}{2}k_1 + \frac{1}{2}q_1} \cdot (2^{k_1} \Vert g_1 \Vert_{\dot{H}^{-1}})^{\frac{2}{3}} (2^{\frac{q_1}{2}} \Vert g_1 \Vert_B)^{\frac{1}{3}} \cdot 2^{\frac{1}{2}q_2} ||g_2||_B \\
	& \lesssim 2^{(2+\beta + 2\delta)m} \cdot 2^{k} \cdot 2^{(2+\frac{1}{6})(k_1 + q_2)} \epsilon^2 \\
	& \lesssim 2^{(-\frac{1}{6} + \beta + 5\delta)m} \epsilon^2,
\end{align*}
so from here, we assume $q_2 - l_2 \leq - (1-\delta)m \ll q_2 - k_2 + k_1$.

Assume that the inequality in \eqref{ibyp2-1-1} fails due to $q_1 \leq -(1-\delta)m$. Together with the above line we have $q_1 \ll q_2 - k_2 + k_1$. Hence, from $|\nabla_{\eta_h} \Phi| \gtrsim |2^{q_2 - k_2} \pm 2^{q_1 - k_1}| \sim 2^{q_2 - k_2}$, we can use set size estimate Lemma \ref{arbitrary-normal-form} along with normal form. For $\lambda = 2^{(1+2\beta)q_2 - 2\beta m}$, the resonant term can be bounded by
\begin{align*}
	2^{(1+\beta)l} \Vert P_{k,p}R_l \mathcal{B}_{\mathfrak{m}^{res}}(g_1, g_2) \Vert_{L^2}
	& \lesssim 2^{(2+\beta + 2\delta)m} \cdot 2^{\frac{k_1+q_1}{2}} \cdot 2^{\frac{k_2}{2}} \cdot (\lambda 2^{k_2-q_2})^{\frac{1}{2}} \cdot 2^{k+q_2} \cdot ||g_1||_2 \cdot ||g_2||_2 \\
	& \lesssim 2^{(2+\beta+2\delta)m} \cdot 2^{\frac{k_1}{2}+2k_2} \lambda^{\frac{1}{2}} 2^{\frac{1}{2}(q_1 + q_2)} \cdot 2^{\frac{q_1}{2}}||g_1||_B \cdot 2^{-(1+\beta)l_2} ||g_2||_X \\
	& \leq 2^{(2+2\delta)m} \cdot 2^{\frac{k_1}{2}+2k_2} \cdot 2^{q_1} \cdot 2^{(1+\beta)(q_2 - l_2)} \epsilon^2 \\
	& \leq 2^{(-\beta + 5\delta + 3\delta_0)m} \epsilon^2.
\end{align*}
Now we turn to nonresonant terms. First, the boundary term can be bounded by
\begin{align*}
	2^{(1+\beta)l} \Vert P_{k,p}R_l \mathcal{Q}_{\mathfrak{n}^{nr}}(g_1, g_2) \Vert_{L^2}
	& \lesssim 2^{(1+\beta + 2\delta)m} \cdot 2^{\frac{k_1+q_1}{2}} \cdot 2^{\frac{k_2}{2}} \cdot (\lambda 2^{k_2-q_2})^{\frac{1}{2}} \cdot 2^{k + q_2} \lambda^{-1} \cdot 2^{\frac{q_1}{2}} ||g_1||_B \cdot 2^{-(1+\beta)l_2} ||g_2||_X \\
	& \lesssim 2^{(1+2\beta + 2\delta)m} \cdot 2^{\frac{k_1}{2}+2k_2} \cdot 2^{q_1} 2^{-\beta q_2} \cdot 2^{-(1+\beta)l_2} \epsilon^2 \\
	& \lesssim 2^{(1+2\beta + 2\delta)m} \cdot 2^{3\delta_0 m} \cdot 2^{-(1-\delta)m} 2^{\frac{14}{15}\beta m} \cdot 2^{-(1+\beta)(\frac{1}{15} - \delta)m} \epsilon^2 \\
	& \lesssim 2^{(-\frac{1}{15} + \frac{43}{15}\beta + 4\delta + 3\delta_0)m} \epsilon^2
\end{align*}
since $q_2 \sim q_{\max} \gtrsim -\frac{14}{15}m$ and hence $-l_2 \leq -(1-\delta)m - q_2 \leq -(\frac{1}{15} - \delta)m$. As $\beta < \frac{1}{43}$, this gives acceptable contribution. The rest are bounded by
\begin{align*}
	2^{(1+\beta)l} \Vert P_{k,p}R_l \mathcal{B}_{\mathfrak{n}^{nr}}((\dt - \kap\Delta)g_1, g_2) \Vert_{L^2}
	& \lesssim 2^{(2+\beta + 2\delta)m} \cdot 2^{\frac{k_1+q_1}{2}} \cdot 2^{\frac{k_2}{2}} \cdot (\lambda 2^{k_2-q_2})^{\frac{1}{2}} \cdot 2^{k + q_2} \lambda^{-1} \cdot 2^{-\frac{3}{2}m + \gamma m} \epsilon^2 \cdot ||g_2||_2 \\
	& \lesssim 2^{(\frac{1}{2} + 2\beta + \gamma + 2\delta)m} \cdot 2^{\frac{k_1}{2}+2k_2} \cdot 2^{-\beta q_2 + \frac{q_1}{2}} \cdot (2^{\frac{q_2}{2}}||g_2||_B)^{\frac{2}{3}} (2^{-(1+\beta)l_2}||g_2||_X)^{\frac{1}{3}} \epsilon^2 \\
	& \lesssim 2^{(3\beta + \gamma + 3\delta + 3\delta_0)m} \cdot 2^{\frac{1}{3}q_2 - \frac{1+\beta}{3}l_2} \epsilon^3 \\
	& \lesssim 2^{(-\frac{1}{3} + 3\beta + \gamma + 3\delta + 3\delta_0)m} \epsilon^3,
\end{align*}
\begin{align*}
	2^{(1+\beta)l} \Vert P_{k,p}R_l \mathcal{B}_{\mathfrak{n}^{nr}}(g_1, (\dt - \kap\Delta)g_2) \Vert_{L^2}
	& \lesssim 2^{(2+\beta + 2\delta)m} \cdot 2^{\frac{k_1+q_1}{2}} \cdot 2^{\frac{k_2}{2}} \cdot (\lambda 2^{k_2-q_2})^{\frac{1}{2}} \cdot 2^{k + q_2} \lambda^{-1} \cdot 2^{\frac{q_1}{2}} ||g_1||_B \cdot 2^{-\frac{3}{2}m + \gamma m} \epsilon^2 \\
	& \lesssim  2^{(\frac{1}{2} + 2\beta + \gamma + 2\delta)m} \cdot 2^{\frac{k_1}{2}+2k_2} \cdot 2^{-\beta q_2 + q_1} \epsilon^3 \\
	& \lesssim 2^{(-\frac{1}{2} + 3\beta + \gamma + 3\delta_0)m} \epsilon^3.
\end{align*}
So now we assume $q_2 - l_1 \leq -(1-\delta)m \ll q_2$ for the remaining complement case of \eqref{ibyp2-1-1}. We use the precise decay estimate again and use $g_j^{I}, g_j^{II}$ notations again. It's easier when $g_1$ has more $S$-vector field applied. In these cases,
\begin{align*}
	2^{(1+\beta)l} \Vert P_{k,p}R_l \mathcal{B}_\mathfrak{m}(g_1, g_2^I) \Vert_{L^2}
	& \lesssim 2^{(2+\beta + 2\delta)m} \cdot 2^{k+q_2} \cdot (2^{\frac{q_1}{2}} ||g_1||_B)^{4\beta} (2^{-(1+\beta)l_1})^{1-4\beta} \cdot 2^{-\frac{3}{2}m - \frac{q_2}{2}} \epsilon \\
	& \lesssim 2^{(\frac{1}{2}+\beta+2\delta)m} \cdot 2^k \cdot 2^{\frac{1}{2}(1+4\beta)(q_2-l_2)} \epsilon^2 \\
	& \lesssim 2^{(-\beta + 3\delta + \delta_0)m} \epsilon^2,
\end{align*}
\begin{align*}
	2^{(1+\beta)l} \Vert P_{k,p}R_l \mathcal{B}_\mathfrak{m}(g_1, g_2^{II}) \Vert_{L^2}
	& \lesssim 2^{(2+\beta + 2\delta)m} \cdot 2^{k+q_2} \cdot 2^{\frac{3}{2}k_1 + \frac{1}{2}q_1} \cdot 2^{-(1+\beta)l_1} ||g_1||_B \cdot 2^{-(1+\beta')m} \epsilon \\
	& \lesssim 2^{(-\beta' + 4\delta)m} \epsilon^2.
\end{align*}
When $g_2$ has more copies of $S$,
\begin{align*}
	2^{(1+\beta)l} \Vert P_{k,p}R_l \mathcal{B}_\mathfrak{m}(g_1^{II}, g_2) \Vert_{L^2}
	& \lesssim 2^{(2+\beta + 2\delta)m} \cdot 2^{k+q_2} \cdot 2^{\frac{3}{2}k_1 + \frac{1}{2}q_1} \cdot 2^{-(1+\beta')m} \epsilon \cdot 2^{-(1+\beta)l_2} ||g_2||_B \\
	& \lesssim 2^{(-\beta' + 4\delta)m} \epsilon^2,
\end{align*}
follows in the same pattern. To deal with the term involving $g_1^I$, we have to subdivide the case, and we first consider when $q_1 \geq -(1-\beta)m - 2\beta q_2$. Then,
\begin{align*}
	2^{(1+\beta)l} \Vert P_{k,p}R_l \mathcal{B}_\mathfrak{m}(g_1^I, g_2) \Vert_{L^2}
	& \lesssim 2^{(2+\beta + 2\delta)m} \cdot 2^{k+q_2} \cdot 2^{-\frac{3}{2}m - \frac{q_1}{2}} \epsilon \cdot 2^{-(1+\beta)l_2} ||g_2||_X \\
	& \lesssim 2^{(\frac{1}{2}+\beta+2\delta)m} \cdot 2^k \cdot 2^{q_2 + \frac{1-\beta}{2}m + \beta q_2} \cdot 2^{-(1+\beta)l_2} \epsilon^2 \\
	& \lesssim 2^{(1+\frac{\beta}{2}+2\delta+\delta_0)m} 2^{(1+\beta)(q_2-l_2)} \epsilon^2 \\
	& \lesssim 2^{(-\frac{\beta}{2} + 4\delta + \delta_0)m} \epsilon^2,
\end{align*}
and we are done. Hence, now we assume $q_1 \leq -(1-\beta)m - 2\beta q_2$. Here, we compare $q_1$ and $q_2-k_2+k_1$. When $q_1 \gtrsim q_2-k_2+k_1$, the same remaining term can be bounded by
\begin{align*}
	2^{(1+\beta)l} \Vert P_{k,p}R_l \mathcal{B}_\mathfrak{m}(g_1^I, g_2) \Vert_{L^2}
	& \lesssim 2^{(2+\beta + 2\delta)m} \cdot 2^{k+q_2} \cdot 2^{\frac{3}{2}k_1 -\frac{3}{2}m - \frac{q_1}{2}} \epsilon \cdot ||g_2||_2 \\
	& \lesssim 2^{(\frac{1}{2}+\beta+2\delta)m} \cdot 2^{k+k_1} \cdot 2^{q_2 + \frac{k_2-q_2}{2}} \epsilon \cdot (2^{\frac{q_2}{2}} ||g_2||_B)^{\frac{1}{3}} (2^{-(1+\beta)l_2} ||g_2||_X) \\
	& \lesssim 2^{(\frac{1}{2} + \beta + 2\delta + 3\delta_0)m} 2^{\frac{2}{3}(q_2-l_2)} \epsilon^2 \\
	& \lesssim 2^{(-\frac{1}{6} + \beta + 3\delta + 3\delta_0)m} \epsilon^2,
\end{align*}
and hence we are finished again. The last case when $q_1 \ll q_2-k_2+k_1$, we can use the set size estimate again using the fact that in this case, $|\nabla_{\eta_h}\Phi| \gtrsim 2^{k_2-q_2}$. The difference comes from the fact that we are no longer in the range $q_1 \leq -(1-\delta)m$, and the additional decay comes from $q_1 \leq -(1-\beta)m - 2\beta q_2$. Setting $\lambda = 2^{(1+6\beta)q_2 - 4\beta m}$,
\begin{align*}
	2^{(1+\beta)l} \Vert P_{k,p}R_l \mathcal{B}_{\mathfrak{n}^{res}}(g_1, g_2) \Vert_{L^2}
	& \lesssim 2^{(2+\beta + 2\delta)m} \cdot 2^{\frac{k_1+q_1}{2}} \cdot 2^{\frac{k_2}{2}} \cdot (\lambda 2^{k_2-q_2})^{\frac{1}{2}} \cdot 2^{k+q_2} \cdot ||g_1||_2 \cdot ||g_2||_2 \\
	& \lesssim 2^{(2+\beta+2\delta)m} \cdot 2^{\frac{k_1}{2}+2k_2} \cdot 2^{(\frac{1}{2}+3\beta)q_2 - 2\beta m} \cdot 2^{\frac{1}{2}(q_1 + q_2)} \cdot 2^{\frac{q_1}{2}}||g_1||_B \cdot 2^{-(1+\beta)l_2} ||g_2||_X \\
	& \leq 2^{(2-\beta + 2\delta)m} \cdot 2^{\frac{k_1}{2}+2k_2} \cdot 2^{-(1-\beta)m - 2\beta q_2} \cdot 2^{(1+3\beta)q_2} \cdot 2^{-(1+\beta)l_2} \epsilon^2 \\
	& \leq 2^{(-\beta + 4\delta + 3\delta_0)m} \epsilon^2.
\end{align*}
Nonresonant terms use normal form along with the set size estimates. First, the boundary term can be bounded by
\begin{align*}
	2^{(1+\beta)l} \Vert P_{k,p}R_l \mathcal{Q}_{\mathfrak{n}^{nr}}(g_1, g_2) \Vert_{L^2}
	& \lesssim 2^{(1+\beta + 2\delta)m} \cdot 2^{\frac{k_1+q_1}{2}} \cdot 2^{\frac{k_2}{2}} \cdot (\lambda 2^{k_2-q_2})^{\frac{1}{2}} \cdot 2^{k + q_2} \lambda^{-1} \cdot 2^{\frac{q_1}{2}} ||g_1||_B \cdot ||g_2||_2 \\
	& \lesssim 2^{(1+3\beta + 2\delta)m} \cdot 2^{\frac{k_1}{2}+2k_2} \cdot 2^{q_1} 2^{-3\beta q_2} \cdot (2^{\frac{1}{2}q_2} ||g_2||_B)^{\frac{2}{3}} (2^{-(1+\beta)l_2} ||g_2||_X)^{\frac{1}{3}} \epsilon \\
	& \lesssim 2^{(4\beta + 2\delta)m} \cdot 2^{3\delta_0 m} \cdot 2^{(\frac{1}{3} - 5\beta)q_2 - \frac{1+\beta}{3}l_2} \epsilon^2 \\
	& \lesssim 2^{(-\frac{1}{3} + 9\beta + 3\delta + 3\delta_0)m} \epsilon^2,
\end{align*}
\begin{align*}
	2^{(1+\beta)l} \Vert P_{k,p}R_l \mathcal{B}_{\mathfrak{n}^{nr}}((\dt - \kap\Delta)g_1, g_2) \Vert_{L^2}
	& \lesssim 2^{(2+\beta + 2\delta)m} \cdot 2^{\frac{k_1+q_1}{2}} \cdot 2^{\frac{k_2}{2}} \cdot (\lambda 2^{k_2-q_2})^{\frac{1}{2}} \cdot 2^{k + q_2} \lambda^{-1} \cdot 2^{-\frac{3}{2}m + \gamma m} \epsilon^2 \cdot ||g_2||_2 \\
	& \lesssim 2^{(\frac{1}{2} + 3\beta + \gamma + 2\delta)m} \cdot 2^{\frac{k_1}{2}+2k_2} \cdot 2^{-3\beta q_2 + \frac{q_1}{2}} \cdot (2^{\frac{q_2}{2}}||g_2||_B)^{\frac{2}{3}} (2^{-(1+\beta)l_2}||g_2||_X)^{\frac{1}{3}} \epsilon^2 \\
	& \lesssim 2^{(\frac{7}{2}\beta + \gamma + 2\delta + 3\delta_0)m} \cdot 2^{(\frac{1}{3}-4\beta)q_2 - \frac{1+\beta}{3}l_2} \epsilon^3 \\
	& \lesssim 2^{(-\frac{1}{3} + \frac{15}{2}\beta + \gamma + 3\delta + 3\delta_0)m} \epsilon^3,
\end{align*}
\begin{align*}
	2^{(1+\beta)l} \Vert P_{k,p}R_l \mathcal{B}_{\mathfrak{n}^{nr}}(g_1, (\dt - \kap\Delta)g_2) \Vert_{L^2}
	& \lesssim 2^{(2+\beta + 2\delta)m} \cdot 2^{\frac{k_1+q_1}{2}} \cdot 2^{\frac{k_2}{2}} \cdot (\lambda 2^{k_2-q_2})^{\frac{1}{2}} \cdot 2^{k + q_2} \lambda^{-1} \cdot ||g_1||_2 \cdot 2^{-\frac{3}{2}m + \gamma m} \epsilon^2 \\
	& \lesssim 2^{(\frac{1}{2} + 3\beta + \gamma + 2\delta)m} \cdot 2^{\frac{k_1}{2}+2k_2} \cdot 2^{-3\beta q_2 + \frac{q_1}{2}} \cdot (2^{\frac{q_1}{2}}||g_1||_B)^{\frac{2}{3}} (2^{-(1+\beta)l_1}||g_1||_X)^{\frac{1}{3}} \epsilon^2 \\
	& \lesssim 2^{(\frac{7}{2}\beta + \gamma + 2\delta + 3\delta_0)m} \cdot 2^{(\frac{1}{3}-4\beta)q_2 - \frac{1+\beta}{3}l_1} \epsilon^3 \\
	& \lesssim 2^{(-\frac{1}{3} + \frac{15}{2}\beta + \gamma + 3\delta + 3\delta_0)m} \epsilon^3.
\end{align*}
This finishes the case 2-1.\\

\noindent 2-2) $k \ll k_1 \sim k_2$: In this case, we have $k+q \sim k_2+q_2$ and $q_2-q \sim k-k_2 \ll 0$, so that $q_1 \ll q_2 \ll q$. The integration by parts already gives enough result if
\begin{gather*}
	2^{-q + k_1 - k}(1 + (2^{q_2 - q_1} + 2^{l_1})) \sim 2^{-q_1} + 2^{-q_2 + l_1} \leq 2^{(1-\delta)m}, \quad \text{or,} \label{ibyp2-2-1} \\
	2^{-q + k_2 - k}(1 + (2^{q_1 - q_2} + 2^{l_2})) \sim 2^{-q_2 + l_2} \leq 2^{(1-\delta)m} \label{ibyp2-2-2}.
\end{gather*}
In the complement cases, we use the normal form using the fact that $|\Phi| \gtrsim 2^q = 2^{q_{\max}}$. The boundary term can be bounded as
\begin{align*}
	2^{(1+\beta)l} \Vert P_{k,p}R_l \mathcal{Q}_\mathfrak{n}(g_1, g_2) \Vert_{L^2}
	& \lesssim 2^{(1+\beta + 2\delta)m} \cdot 2^{k} \cdot ||e^{it\Lambda} g_1||_\infty ||g_2||_2 \\
	& \lesssim 2^{(1+\beta + 2\delta)m} \cdot 2^{k} \cdot 2^{-m+\gamma m} ||g_1||_D \cdot (2^{\frac{q_2}{2}} ||g_2||_B)^{\frac{2}{3}} (2^{-(1+\beta)l_2} ||g_2||_X)^{\frac{1}{3}} \\
	& \lesssim 2^{(-\frac{1}{3} + \beta + \gamma + 3\delta + \delta_0)m} \epsilon^2
\end{align*}
using Lemma \ref{slower-dispersive-decay}. Lastly, the time derivative terms can be bounded easily as
\begin{align*}
	2^{(1+\beta)l} \Vert P_{k,p}R_l \mathcal{B}_\mathfrak{n}((\dt-\kap\Delta)g_1, g_2) \Vert_{L^2}
	& \lesssim 2^{(2+\beta + 2\delta)m} \cdot 2^{k} \cdot 2^{\frac{3}{2}k_1 + \frac{1}{2}q_1} \cdot 2^{-\frac{3}{2}m + \gamma m} \epsilon^2 \cdot (2^{\frac{q_2}{2}} ||g_2||_B)^{\frac{1}{3}} (2^{-(1+\beta)l_2} ||g_2||_X)^{\frac{2}{3}} \\
	& \lesssim 2^{(\frac{1}{2} + \beta+\gamma +2\delta + 3\delta_0)m} \cdot 2^{\frac{2}{3}(q_2 - l_2)} \epsilon^3 \\
	& \leq 2^{(-\frac{1}{6} + \beta+\gamma + 3\delta)m} \epsilon^3,
\end{align*}
\begin{align*}
	2^{(1+\beta)l} \Vert P_{k,p}R_l \mathcal{B}_\mathfrak{n}(g_1, (\dt-\kap\Delta)g_2) \Vert_{L^2}
	& \lesssim 2^{(2+\beta + 2\delta)m} \cdot 2^{k} \cdot 2^{\frac{3}{2}k_1 + \frac{1}{2}q_1} \cdot (2^{\frac{q_1}{2}} ||g_1||_B)^{\frac{1}{3}} (2^{-(1+\beta)l_1} ||g_1||_X)^{\frac{2}{3}} \cdot 2^{-\frac{3}{2}m + \gamma m} \epsilon^2 \\
	& \lesssim 2^{(\frac{1}{2} + \beta+\gamma +2\delta + 3\delta_0)m} \cdot 2^{\frac{2}{3}(q_1 - l_2)} \epsilon^3 \\
	& \leq 2^{(-\frac{1}{6} + \beta+\gamma + 3\delta)m} \epsilon^3,
\end{align*}
since $q_1 - l_1 \leq \min\{q_1, q_2-l_1\}$.\\

\noindent 2-3) $k_2 \ll k_1 \sim k$. In this case, we have $k+q \sim k_2+q_2$ and $q-q_2 \sim k_2-k \ll 0$, so that $q_1 \ll q \ll q_2$. The case is symmetric to case 2-2 and the estimates follow similarly. The integration parts are enough if
\begin{gather*}
	2^{-q_2 + k_1 - k_2}(1 + 2^{k_2-k_1}(2^{q_2 - q_1} + 2^{l_1})) \sim 2^{-q_1} + s^{-q_2 + l_1} \leq 2^{(1-\delta)m}, \quad \text{or,} \\
	2^{-q_2 + k_2 - k_1}(1 + 2^{k_1-k_2}(2^{q_1 - q_2} + 2^{l_2})) \sim 2^{-q_2 + l_2} \leq 2^{(1-\delta)m}.
\end{gather*}
The conditions on the indices are same as in 2-2, and we can still use normal form since $|\Phi| \gtrsim 2^{q_{\max}} = 2^{q_2}$. Hence, the result follows exactly the same way as in the case 2-2.\\

\noindent \textcircled{3} $q \sim q_1 \ll q_2$ : By Lemma \ref{two_small_sizes}, we have $k \sim k_1 \gg k_2$, $k+q \sim k_1+q_1 \gtrsim k_2+q_2$ in this case. Integration by parts is sufficient if
\begin{gather*}
	2^{-q_2 + 2k_1 - k_{\max} - k_{\min}}(1 + 2^{k_2 - k_1}(2^{q_2 - q_1} + 2^{l_1})) \sim 2^{-q_2}(2^{k_1 - k_2} + 2^{l_1}) \leq 2^{(1-\delta)m}, \quad \text{or,} \\
	2^{-q_2+2k_2 - k_{\max} - k_{\min}}(1 + 2^{k_1 - k_2}(2^{q_1 - q_2} + 2^{l_2})) \sim 2^{-q_2 + l_2} \leq 2^{(1-\delta)m}.
\end{gather*}
For the remaining case, we use normal form since $|\Phi| \sim 2^{q_2}$.
\begin{align*}
	2^{(1+\beta)l} \Vert P_{k,p}R_l \mathcal{Q}_\mathfrak{n}(g_1, g_2) \Vert_{L^2}
	& \lesssim 2^{(1 + \beta + 2\delta)m} \cdot 2^k \cdot ||e^{it\Lambda} g_1 ||_\infty \cdot (2^{\frac{q_2}{2}} ||g_2||_B)^{\frac{2}{3}} (2^{-(1+\beta)l_2} ||g_2||_X)^{\frac{1}{3}} \\
	& \lesssim 2^{(1 + \beta + 2\delta)m} \cdot 2^k \cdot 2^{-m+\gamma m} \epsilon \cdot 2^{\frac{q_2 - l_2}{3}} \epsilon \\
	& \lesssim 2^{(-\frac{1}{3} + \beta + \gamma + 3\delta)m} \epsilon^2,
\end{align*}
where the second inequality holds by Lemma \ref{slower-dispersive-decay}. For the derivative terms,
\begin{align*}
	2^{(1+\beta)l} \Vert P_{k,p}R_l \mathcal{B}_\mathfrak{n}((\dt-\kap\Delta)g_1, g_2) \Vert_{L^2}
	& \lesssim 2^{(2 + \beta + 2\delta)m} \cdot 2^k \cdot 2^{\frac{3}{2}k_2 + \frac{1}{2}q_2} \cdot ||(\dt - \kap\Delta)g_1||_2 \cdot (2^{\frac{q_2}{2}} ||g_2||_B)^{\frac{1}{3}} (2^{-(1+\beta)l_2} ||g_2||_X)^{\frac{2}{3}} \\
	& \lesssim 2^{(2 + \beta + 2\delta)m} \cdot 2^{\frac{5}{2}k} \dot 2^{-\frac{3}{2}m + \gamma m} \epsilon^2 \cdot 2^{\frac{2}{3}(q_2 - l_2)} \epsilon \\
	& \lesssim 2^{(-\frac{1}{6} + \beta + \gamma + 3\delta + 3\delta_0)m} \epsilon^3,
\end{align*}
\begin{align*}
	2^{(1+\beta)l} \Vert P_{k,p}R_l \mathcal{B}_\mathfrak{n}(g_1, (\dt-\kap\Delta)g_2) \Vert_{L^2}
	& \lesssim 2^{(2 + \beta + 2\delta)m} \cdot 2^k \cdot 2^{\frac{3}{2}k_2 + \frac{1}{2}q_2} \cdot (2^{\frac{q_1}{2}} ||g_2||_B)^{\frac{1}{3}} (2^{-(1+\beta)l_1} ||g_2||_X)^{\frac{2}{3}} \cdot ||(\dt-\kap\Delta)g_2||_2 \\
	& \lesssim 2^{(2 + \beta + 2\delta)m} \cdot 2^k \cdot 2^{\frac{2}{3}(q_2 - l_1) + \frac{3}{2}k_2} \epsilon \cdot 2^{-\frac{3}{2}m + \gamma m} \epsilon^2 \\
	& \lesssim 2^{(-\frac{1}{6} + \beta + \gamma + 3\delta + 3\delta_0)m} \epsilon^3.
\end{align*}

\noindent\textcircled{4} $q_1 \sim q_2 \ll q$ : Again by Lemma \ref{two_small_sizes}, we have $k \ll k_1 \sim k_2$ and $k+q \lesssim k_1+q_1 \sim k_2+q_2$. Integration by parts is sufficient if
\begin{align*}
	& 2^{-q + 2k_1 - k_{\max} - k_{\min}}(1 + 2^{k_2 - k_1}(2^{q_2 - q_1} + 2^{l_1})) \sim 2^{-q + k_1 - k + l_1} \leq 2^{(1-\delta)m}, \quad \text{or,} \\
	& 2^{-q + 2k_2 - k_{\max} - k_{\min}}(1 + 2^{k_1 - k_2}(2^{q_1 - q_2} + 2^{l_2})) \sim 2^{-q + k_2 - k + l_2} \leq 2^{(1-\delta)m}.
\end{align*}
Hence, now we assume $q + k - k_1 - l_j \leq -(1-\delta)m$, $j=1,2$. Since the conditions are symmetric on $\eta$ and $\xi-\eta$, assume WLOG $g_1$ has more copies of $S$. Using the normal form,
\begin{align*}
	2^{(1+\beta)l} \Vert P_{k,p}R_l \mathcal{Q}_\mathfrak{n}(g_1, g_2) \Vert_{L^2}
	& \lesssim 2^{(1 + \beta + 2\delta)m} \cdot 2^k \cdot (2^{\frac{q_1}{2}} ||g_1||_B)^{\frac{2}{3}} (2^{-(1+\beta)l_1} ||g_1||_X)^{\frac{1}{3}} \cdot ||e^{it\Lambda} g_2||_\infty \\
	& \lesssim 2^{(1 + \beta + 2\delta)m} \cdot 2^{\frac{2}{3}k + \frac{1}{3}k_1} \cdot 2^{\frac{1}{3}(q_1 + k - k_1 - l_1)} \cdot 2^{-m} \epsilon^2 \\
	& \lesssim 2^{(-\frac{1}{3} + \beta + 3\delta + \delta_0)m} \epsilon^2,
\end{align*}
\begin{align*}
	2^{(1+\beta)l} \Vert P_{k,p}R_l \mathcal{B}_\mathfrak{n}((\dt-\kap\Delta)g_1, g_2) \Vert_{L^2}
	& \lesssim 2^{(2 + \beta + 2\delta)m} \cdot 2^k \cdot 2^{\frac{3}{2}k + \frac{1}{2}q} \cdot 2^{-\frac{3}{2}m + \gamma m} \epsilon^2 \cdot (2^{\frac{q_2}{2}} ||g_2||_B)^{\frac{1}{3}} (2^{-(1+\beta)l_2} ||g_2||_X)^{\frac{2}{3}} \\
	& \lesssim 2^{(\frac{1}{2} + \beta + \gamma + 2\delta)m} \cdot 2^{\frac{11}{6}k + \frac{4}{6}k_2} \cdot 2^{\frac{2}{3}(q_2 + k - k_2 - l_2)} \epsilon^3 \\
	& \lesssim 2^{(-\frac{1}{6} + \beta + \gamma + 3\delta + 3\delta_0)m} \epsilon^3,
\end{align*}
\begin{align*}
	2^{(1+\beta)l} \Vert P_{k,p}R_l \mathcal{B}_\mathfrak{n}(g_1, (\dt-\kap\Delta)g_2) \Vert_{L^2}
	& \lesssim 2^{(2 + \beta + 2\delta)m} \cdot 2^k \cdot 2^{\frac{3}{2}k + \frac{1}{2}q} \cdot (2^{\frac{q_1}{2}} ||g_1||_B)^{\frac{2}{3}} (2^{-(1+\beta)l_1} ||g_1||_X)^{\frac{1}{3}} \cdot 2^{-\frac{3}{2}m + \gamma m} \epsilon^2 \\
	& \lesssim 2^{(\frac{1}{2} + \beta + \gamma + 2\delta)m} \cdot 2^{\frac{11}{6}k + \frac{4}{6}k_1} \cdot 2^{\frac{2}{3}(q_1 + k - k_1 - l_1)} \epsilon^3 \\
	& \lesssim 2^{(-\frac{1}{6} + \beta + \gamma + 3\delta + 3\delta_0)m} \epsilon^3.
\end{align*}
This finishes the cases with gap in $q$. \\

\noindent \textbf{7.2.4 No Gaps} \\
Now we can assume $p \sim p_1 \sim p_2 \sim 0$, and $-\frac{14}{15}m \leq q_{\max} \sim q_{\min}$. WLOG we assume $g_2$ has more copies of $S$ applied, and divide into resonant and nonresonant terms according to $\lambda = 2^{q-20}$.

We consider the nonresonant term first and perform normal form. Then,
\begin{align*}
	2^{(1+\beta)l} \Vert P_{k,p}R_l \mathcal{Q}_\mathfrak{n}(g_1^I, g_2) \Vert_{L^2}
	& \lesssim 2^{(1 + \beta + 2\delta)m} \cdot 2^k \cdot ||e^{it\Lambda} g_1^I||_\infty \cdot 2^{\frac{q_2}{2}} ||g_2||_B \\
	& \lesssim 2^{(1 + \beta + 2\delta)m} \cdot 2^k \cdot 2^{-\frac{3}{2}m - \frac{q_1}{2}} ||g_1||_D \cdot 2^{\frac{q_2}{2}} \epsilon \\
	& \lesssim 2^{(-\frac{1}{2} + \beta + 2\delta + \delta_0)m} \epsilon^2,
\end{align*}
\begin{align*}
	2^{(1+\beta)l} \Vert P_{k,p}R_l \mathcal{Q}_\mathfrak{n}(g_1^{II}, g_2) \Vert_{L^2}
	& \lesssim 2^{(1 + \beta + 2\delta)m} \cdot 2^k \cdot ||e^{it\Lambda} g_1^{II}||_2 \cdot ||e^{it\Lambda}g_2||_\infty \\
	& \lesssim 2^{(1 + \beta + 2\delta)m} \cdot 2^k \cdot 2^{-m - \beta'm} ||g_1||_D \cdot 2^{-m + \gamma m} \epsilon \\
	& \lesssim 2^{(-\frac{1}{2} + \beta + \gamma - \beta' + 2\delta + \delta_0)m} \epsilon^2,
\end{align*}
\begin{align*}
	2^{(1+\beta)l} \Vert P_{k,p}R_l \mathcal{B}_\mathfrak{n}((\dt-\kap\Delta)g_1, g_2) \Vert_{L^2}
	& \lesssim 2^{(2 + \beta + 2\delta)m} \cdot 2^k \cdot ||(\dt-\kap\Delta)g_1||_2 \cdot ||e^{it\Lambda} g_2||_2 \\
	& \lesssim 2^{(2 + \beta + 2\delta)m} \cdot 2^k \cdot 2^{-\frac{3}{2}m + \gamma m} \epsilon^2 \cdot 2^{-m + \gamma m} \epsilon \\
	& \lesssim 2^{(-\frac{1}{2} + \beta + 2\gamma + \delta + \delta_0)m} \epsilon^3,
\end{align*}
\begin{align*}
	2^{(1+\beta)l} \Vert P_{k,p}R_l \mathcal{B}_\mathfrak{n}(g_1, (\dt-\kap\Delta)g_2) \Vert_{L^2}
	& \lesssim 2^{(2 + \beta + 2\delta)m} \cdot 2^k \cdot ||e^{it\Lambda} g_1||_\infty \cdot ||(\dt-\kap\Delta)g_2||_2 \\
	& \lesssim 2^{(2 + \beta + 2\delta)m} \cdot 2^k \cdot 2^{-m} ||g_1||_D \cdot 2^{-\frac{3}{2}m + \gamma m} \epsilon^2 \\
	& \lesssim 2^{(-\frac{1}{2} + \beta + \gamma + \delta + \delta_0)m} \epsilon^3.
\end{align*}

For resonant term, trivial estimate gives
\begin{align*}
	2^{(1+\beta)l} \Vert P_{k,p}R_l \mathcal{B}_\mathfrak{m}(g_1, g_2) \Vert_{L^2}
	& \lesssim 2^{(2 + \beta + 2\delta)m} \cdot 2^{k+q} \cdot 2^{\frac{3}{2}k_{\min} + \frac{1}{2}q} \cdot 2^{\min\{ \frac{q_1}{2}, k_1\}} \cdot 2^{\min\{ \frac{q_2}{2}, k_2\}} \epsilon^2 \\
	& \lesssim 2^{(2 + \beta + 2\delta)m} \cdot 2^{k_{\max}} \cdot 2^{\frac{3}{2}k_{\min} + 2q} \cdot 2^{\min\{ \frac{q}{2}, k_{\min}\}} \epsilon^2 \\
	& \lesssim 2^{(2 + \beta + 2\delta)m} \cdot 2^{k_{\max}} \cdot 2^{\frac{13}{6}(k_{\min} + q)} \epsilon^2,
\end{align*}
which is enough if $k_{\min} + q \leq -(1-\beta)m$. Hence, we assume $k_{\min} + q \geq -(1-\beta)m$. Now, since $|\Phi| \leq 2^{q-20}$ for the resonant term, from Proposition \ref{no-spacetime-resonance} we can use integration by parts, and this is sufficient if
\begin{align*}
	& 2^{2k_1 - q - k_{\max} - k_{\min}}(1 + 2^{k_2 - k_1}(2^{q_2 - q_1} + 2^{l_1})) \sim 2^{2k_1 - q - k_{\max} - k_{\min}}(1 + 2^{k_2 - k_1 + l_1}) \leq 2^{(1-\delta)m}, \quad \text{or,} \\
	& 2^{2k_2 - q - k_{\max} - k_{\min}}(1 + 2^{k_1 - k_2}(2^{q_1 - q_2} + 2^{l_2})) \sim 2^{2k_2 - q - k_{\max} - k_{\min}}(1 + 2^{k_1 - k_2 + l_2}) \leq 2^{(1-\delta)m}.
\end{align*}
In other words, if $\max\{2k_j, k_1+k_2+l_j\} \leq k_{\max} + k_{\min} + q + (1-\delta)m$, we are done. Thus we are in the case
\begin{equation*}
	\max\{2k_j, k_1+k_2+l_j\} \geq k_{\max} + k_{\min} + q + (1-\delta)m \geq k_{\max} - (1-\beta)m + (1-\delta)m = k_{\max} + (\beta - \delta)m.
\end{equation*}
If $2k_j = \max\{2k_j, k_1+k_2+l_j\}$ for any $j=1,2$, we get a contradiction by $\delta_0 m \geq k_{\max} \geq 2k_j - k_{\max} \geq (\beta - \delta)m$. Therefore, $k_1+k_2+l_j$'s have to be the maximums. Plugging this back to get $q-l_j + k_{\max} + k_{\min} - k_1 - k_2 \leq -(1-\delta)m$,
\begin{align*}
	2^{(1+\beta)l} \Vert P_{k,p}R_l \mathcal{B}_\mathfrak{m}(g_1^I, g_2) \Vert_{L^2}
	& \lesssim 2^{(2 + \beta + 2\delta)m} \cdot 2^{k+q} \cdot ||e^{it\Lambda} g_1^I||_\infty \cdot (2^{\frac{q_2}{2}} ||g_2||_B)^{\frac{1}{3}} (2^{-(1+\beta)l_2} ||g_2||_X)^{\frac{2}{3}} \\
	& \lesssim 2^{(2 + \beta + 2\delta)m} \cdot 2^{k+q} \cdot 2^{-\frac{3}{2}m - \frac{q}{2}} \epsilon \cdot 2^{\frac{1}{6}q - \frac{2}{3}(1+\beta)l_2} \epsilon \\
	& \lesssim 2^{(\frac{1}{2} + \beta + 2\delta)m} \cdot 2^{k + \frac{2}{3}(k_1 + k_2 - k_{\max} - k_{\min} -(1-\delta)m)} \epsilon^2 \\
	& \lesssim 2^{(-\frac{1}{2} + \beta + 3\delta + \delta_0)m} \epsilon^2,
\end{align*}
\begin{align*}
	2^{(1+\beta)l} \Vert P_{k,p}R_l \mathcal{B}_\mathfrak{m}(g_1^{II}, g_2) \Vert_{L^2}
	& \lesssim 2^{(2 + \beta + 2\delta)m} \cdot 2^{k+q} \cdot 2^{\frac{3}{2}k + \frac{1}{2}q} \cdot ||e^{it\Lambda} g_1^{II}||_2 \cdot 2^{-(1+\beta)l_2} ||g_2||_X \\
	& \lesssim 2^{(2 + \beta + 2\delta)m} \cdot 2^{-(1+\beta')m} \epsilon \cdot 2^{\frac{5}{2}k + (1+\beta)(q - l_2)} \epsilon \\
	& \lesssim 2^{(-\beta' + 4\delta + 3\delta_0)m} \epsilon^2,
\end{align*}
and this completes the proof.

\bigskip
\bigskip

\textbf{Acknowledgments.} The author would like to thank Beno\^{i}t Pausader for his support and discussions during the process of this work. The author was partially supported through his NSF grant DMS-2154162 while writing this paper.

\newpage

\appendix
\section{Additional lemmas}\label{apdx}
\subsection{Properties of angular localization}
\begin{proposition}
	For any $l \in \mathbb{N}$, the angular localization operator $\bar{R}_l$ satisfies the following.\\
	(1) $\bar{R}_l$ commmutes with the regular Littlewood-Paley projectors, vector fields $S$ and $\Omega_{a b} = x_z \partial_{x_b} - x_b \partial_{x_a}$, the Fourier transform, and the Laplacian. In other words,
	\begin{equation*}
		[P_k, \bar{R}_l] = [S, \bar{R}_l] = [\Omega_{a b}, \bar{R}_l] = [\mathcal{F}, \bar{R}_l] = [\Delta, \bar{R}_l] = 0.
	\end{equation*}	
	(2) $\bar{R}_l$ constitutes an orthogonal partition of unity in the sense that
	\begin{equation*}
		f = \sum_{l \geq 0} \bar{R}_l f, \quad ||f||_2^2 \simeq \sum_{l \geq 0} ||\bar{R}_l f||_2^2, \quad \bar{R}_l \bar{R}_{l'} = 0 \text{ if } |l-l'| > 3.
	\end{equation*}
	(3) Bernstein property holds, i.e.
	\begin{equation*}
		\sum_{1 \leq a < b \leq 3} ||\Omega_{a b} \bar{R}_l f ||_{L^r} \simeq 2^l ||\bar{R}_l f||_{L^r}.
	\end{equation*}
\end{proposition}
\begin{proof}
	Everything is proved in \cite[Proposition 3.1]{GPW} except $[\Delta, \bar{R}_l] = 0$. However, this is evident since $\bar{R}_l$ commutes with the Fourier transform, hence
	\begin{equation*}
		\mathcal{F}\{\bar{R}_l (-\Delta f)\} (\xi)
		= \sum_{n \geq 0} \varphi(2^{-l}n) \int_{\mathbb{S}^2} ||\xi|\theta|^2 f(|\xi|\theta) \mathfrak{Z}_n(\langle \theta, \frac{x}{|x|} \rangle) dS(\theta)
		= |\xi|^2(\bar{R}_l \hat{f})(\xi)
		= \mathcal{F}\{(-\Delta \bar{R}_l f)\}(\xi).
	\end{equation*}
\end{proof}

\subsection{Set size estimates}
The following lemmas are used to estimate the bilinear terms after reducing the indices to certain cases, by using the set size in Fourier space.
\begin{lemma}\label{set-size-estimate}
	(1) For a bilinear term $\mathcal{Q}_{\mathfrak{m}}$ with a multiplier $\mathfrak{m}$ and localized by $P_{k,p,q}$, i.e.
	\begin{equation*}
		P_{k,p,q}\mathcal{Q}_\mathfrak{m}(P_{k_1,p_1,q_1} f_1, P_{k_2,p_2,q_2} f_2) = \mathcal{F}^{-1} \left\{\int_{\mathds{R}^3_\eta} e^{is\Phi} \chi(\xi,\eta) \mathfrak{m}(\xi,\eta) \hat{f_1}(\xi-\eta) \hat{f_2}(\eta) d\eta \right\},
	\end{equation*}
	then with
	\begin{equation*}
		|S| := \min\{2^{k+p}, 2^{k_1+p_1}, 2^{k_2+p_2}\} \cdot \min\{2^{\frac{k+q}{2}}, 2^{\frac{k_1+q_1}{2}},2^{\frac{k_2+q_2}{2}}\},
	\end{equation*}
	the following bound holds:
	\begin{equation*}
		||P_{k,p,q}\mathcal{Q}_\mathfrak{m}(P_{k_1,p_1,q_1} f_1, P_{k_2,p_2,q_2} f_2)||_2 \lesssim |S| \cdot ||\mathfrak{m} \chi||_\infty \cdot ||P_{k_1,p_1,q_1} f_1||_2 \cdot ||P_{k_2,p_2,q_2} f_2)||_2.
	\end{equation*}
	(2) If the bilinear term is localized where the phase is small in addition, i.e.
	\begin{equation*}
		\mathcal{Q}_{\mathfrak{m}}(f_1, f_2) = \mathcal{F}^{-1} \left\{\int_{\mathds{R}^3_\eta} e^{is\Phi} \varphi(\lambda^{-1}\Phi) \chi(\xi,\eta) \mathfrak{m}(\xi,\eta) \hat{f_1}(\xi-\eta) \hat{f_2}(\eta) d\eta \right\},
	\end{equation*}
	and we have a lower bound $|\nabla_{\eta_h} \Phi| \gtrsim L > 0$ on the support of $\chi$, then
	\begin{equation*}
		||\mathcal{Q}_\mathfrak{m}(f_1, f_2)||_2 \lesssim \min\{2^{\frac{k_1+q_1}{2}},2^{\frac{k_2+q_2}{2}}\} \cdot 2^{\frac{k_2+p_2}{2}} \cdot (\lambda L^{-1})^{\frac{1}{2}} \cdot ||\mathfrak{m} \chi||_\infty \cdot ||P_{k_1,p_1,q_1} f_1||_2 \cdot ||P_{k_2,p_2,q_2} f_2)||_2.
	\end{equation*}
\end{lemma}
\begin{proof}
	The estimates use the fact that the localizations, as well as the conditions on $\Phi$ in case of the estimate on (2), give the set size related terms. See \cite[Lemma A.3, A.4]{GPW}.
\end{proof}

\subsection{Slower dispersive decay}
\begin{lemma}\label{slower-dispersive-decay}
	Under the bootstrap assumption \eqref{BootstrapAssumption}, for some $0 < \gamma \ll \beta$, we have
	\begin{equation*}
		||P_k e^{it\Lambda} S^b \U||_\infty \lesssim 2^{\frac{3}{2}k - 3k^+} t^{-1+\gamma} \epsilon, \quad 0 \leq b \leq N.
	\end{equation*}
\end{lemma}
\begin{proof}
	See \cite[Corollary A.7]{GPW}.
\end{proof}

\bibliographystyle{siam}
\bibliography{Navier-Stokes-Coriolis-Submission}

\end{document}